%% file: 2005-32.tex
\documentclass{gtart_h}  
\input gtoutput

\lognumber{545}
\received{18 December 2004}
\volumenumber{9}\papernumber{32}\volumeyear{2005}
\pagenumbers{1381}{1441}   
\revised{2 August 2005}
\published{5 August 2005}
\accepted{4 July 2005}
\proposed{Joan Birman}
\seconded{Walter Neumann, Martin Bridson}

\usepackage{amsmath, amssymb, amsthm}
\usepackage{graphicx, psfrag} 
\def\psfraga <#1,#2> #3#4{%
\psfrag {#3}{\smash{\rlap{\kern #1 \raise #2\hbox{#4}}}}}
\def\figref#1{\hyperlink{#1anchor}{Figure~\ref*{#1}}}
\def\anchor#1{\noindent\hypertarget{#1anchor}{\smash{$\phantom{99}$}}\newline}

\theoremstyle{plain}

\newtheorem{thm}{Theorem}
\newtheorem{prop}[thm]{Proposition}
\newtheorem{lemma}[thm]{Lemma}
\newtheorem*{claim}{Claim}
\newtheorem*{sublemma}{Sublemma}

\theoremstyle{definition}
\newtheorem*{defn}{Definition}
\newtheorem*{remark}{Remark}
\newtheorem*{example}{Example}
\newtheorem*{notation}{Notation}
\newtheorem*{ackn}{Acknowledgement}
\newtheorem{defnnum}[thm]{Definition}
\newtheorem{remarknum}[thm]{Remark}

\def\ubrace#1#2{\underbrace{#1}\limits_{\displaystyle{#2}}}
\def\Cal#1{{\cal#1}}
\def\V{{\cal V}}\def\W{{\cal W}}\def\F{{\cal F}}

\def\<{\langle}\def\>{\rangle}\def\inv{^{-1}}
\def\what{\widehat}\def\wtil{\widetilde}\def\ov{\overline}

\def\Z{{\mathbb Z}}\def\N{{\mathbb N}} 
\def\R{{\mathbb R}} \def\E{{\mathbb E}}
\def\D{{\mathbb D}}\def\mF{{\mathbb F}}\def\X{{\bf X}}

\def\Lk{\text{\sl Lk}} \def\pro{\text{\sl prod}}
\def\Aut{\text{\sl Aut}}\def\Out{\text{\sl Out}}\def\Fix{\text{\sl Fix}}
\def\Iso{\text{\sl Iso}}\def\Inv{\text{\sl Inv}}\def\Stab{\text{\sl Stab}}
\def\Comm{\text{\sl Comm}}\def\Sym{\text{\sl Sym}}\def\Inn{\text{\sl Inn}}
\def\Pure{\text{\sl Pure}}\def\Biject{\text{\sl Biject}}\def\Twist{\text{\sl Twist}}

\def\al{\alpha}\def\be{\beta}\def\ga{\gamma}\def\sig{\sigma}
\def\Ga{\Gamma}\def\Sig{\Sigma}\def\ep{\epsilon}
\begin{document}

\title{Automorphisms and abstract commensurators\\of 2--dimensional Artin groups}
\shorttitle{Automorphisms and abstract commensurators}
\asciititle{Automorphisms and abstract commensurators of 2-dimensional Artin groups}                    
\authors{John Crisp}                  
\address{IMB(UMR 5584 du CNRS), Universit\'e de Bourgogne\\BP 47 870, 
21078 Dijon, France}                  
\asciiaddress{IMB(UMR 5584 du CNRS), Universite de Bourgogne\\BP 47 870, 
21078 Dijon, France}                  
            
\email{jcrisp@u-bourgogne.fr}                     
\urladdr{http://math.u-bourgogne.fr/IMB/crisp}                       

\begin{abstract}  
In this paper we consider the class of 2--dimensional Artin groups
with connected, large type, triangle-free defining graphs (type
CLTTF).  We classify these groups up to isomorphism, and describe a
generating set for the automorphism group of each such Artin group. In
the case where the defining graph has no separating edge or vertex we
show that the Artin group is not abstractly commensurable to any other
CLTTF Artin group.  If, moreover, the defining graph satisfies a
further ``vertex rigidity'' condition, then the abstract commensurator
group of the Artin group is isomorphic to its automorphism group and
generated by inner automorphisms, graph automorphisms (induced from
automorphisms of the defining graph), and the involution which maps
each standard generator to its inverse.

We observe that the techniques used here to study automorphisms 
carry over easily to the Coxeter group situation.
We thus obtain a classification of the CLTTF type Coxeter groups up to
isomorphism and a description of their automorphism groups analogous to that 
given for the Artin groups.
\end{abstract}

\asciiabstract{%
In this paper we consider the class of 2-dimensional Artin groups with
connected, large type, triangle-free defining graphs (type CLTTF).  We
classify these groups up to isomorphism, and describe a generating set
for the automorphism group of each such Artin group. In the case where
the defining graph has no separating edge or vertex we show that the
Artin group is not abstractly commensurable to any other CLTTF Artin
group.  If, moreover, the defining graph satisfies a further `vertex
rigidity' condition, then the abstract commensurator group of the
Artin group is isomorphic to its automorphism group and generated by
inner automorphisms, graph automorphisms (induced from automorphisms
of the defining graph), and the involution which maps each standard
generator to its inverse.  We observe that the techniques used here to
study automorphisms carry over easily to the Coxeter group situation.
We thus obtain a classification of the CLTTF type Coxeter groups up to
isomorphism and a description of their automorphism groups analogous
to that given for the Artin groups.}

\asciikeywords{2-dimensional Artin group, Coxeter group, commensurator group,
graph automorphisms, triangle free}                    

\maketitle

\section*{Introduction and statement of results}
\addcontentsline{toc}{section}{Introduction and statement of results}

Let $\Delta$ denote a simplicial graph with vertex set $V(\Delta)$ and 
edge set $E(\Delta)\subset V(\Delta)\times V(\Delta)$.
Suppose also that every edge $e=\{ s,t\}\in E(\Delta)$ 
carries a label $m_e=m_{st}\in \N_{\geq 2}$. 
We define the \emph{Artin group} $G(\Delta)$
associated to the (labelled) \emph{defining graph} $\Delta$
to be the group given by the presentation
\footnote{Our notion of defining graph differs from the frequently used
``Coxeter graph'' where, by contrast, the absence
of an edge between $s$ and $t$  indicates a commuting
relation ($m_{st}=2$) and the label $m_{st}=\infty$
is used to designate the absence of a relation between $s$ and $t$.
In our convention the label $\infty$ is never used.} 
\[
G(\Delta)=\<\ V(\Delta)\  \mid\  
\ubrace{ststs\cdots}{m_{st}}=\ubrace{tstst\cdots}{m_{st}}
\ \text{ for all } \{ s,t\}\in E(\Delta)\ \>\,.
\]
Adding the relations $s^2=1$ for each $s\in V(\Delta)$ yields a presentation 
of the associated \emph{Coxeter group} $W(\Delta)$ of type $\Delta$. We denote 
$\rho_\Delta\co G(\Delta)\to W(\Delta)$ the canonical quotient map obtained by
the addition of these relations.

The following observations are true for all Artin groups and were proved in \cite{vdL}. 
If $T$ is a full subgraph of $\Delta$ then the subgroup of $G(\Delta)$
generated by the vertices of $T$ is canonically isomorphic to $G(T)$. Such subgroups
shall be called \emph{standard parabolic}. Moreover, the intersection of two 
standard parabolic subgroups of an Artin group is again a standard parabolic subgroup.
Thus, for example, if $e,f\in E(\Delta)$ then $G(e)\cap G(f)=G(e\cap f)$, which is either 
the cyclic group $\<s\>$ in the case that $e$ and $f$
share a common vertex $s$, or the trivial group (in the case $e$ and $f$ are disjoint).
The analogous statements also hold for Coxeter groups.

\begin{defn}[{\rm(}CLTTF Artin group\/{\rm)}]
The main Theorems in this paper shall apply to Artin (and Coxeter) groups whose 
defining graph satisfies  the following conditions: 
\begin{description}
\item[\rm(C)] $\Delta$ is connected and has at least 3 vertices; 
\item[\rm(LT)] all labels $m_e$, for $e\in E(\Delta)$, are at least 3; and
\item[\rm(TF)] $\Delta$ has no triangles (no simple circuits of length 3).
\end{description}
If $\Delta$ satisfies all three of the above conditions then we refer to it as a 
\emph{CLTTF defining graph} and we refer to $G(\Delta)$ as a \emph{CLTTF Artin group}, 
and to $W(\Delta)$ as a \emph{CLTTF Coxeter group}.
\end{defn}

Conditions (LT) and (TF) correspond to 
what are known as the \emph{large type} and \emph{triangle free} conditions, either of 
which implies that the Artin group has cohomological (or geometric) dimension 2.
The triangle free Artin groups are exactly the 2--dimensional, so-called, ``FC type''
Artin groups.
The condition (C) simply serves to rule out the $2$--generator or ``dihedral 
type'' Artin groups which are best treated as a separate case 
(see \cite{GHMR} for a treatment of their automorphism groups), as well as 
those Artin groups which are proper free products. (Using the Kurosh Subgroup Theorem
it can be shown that an arbitrary Artin group $G(\Delta)$ is freely 
indecomposable if and only if $\Delta$ is connected).

\begin{thm}\label{Thm1}
Let $\Cal G$ denote the set of all CLTTF defining graphs (up to labelled graph
isomorphism) and write $\Iso(\Cal G)$ for the category (a groupoid) with objects $\Cal G$ and 
morphisms the set of all isomorphisms $G(\Delta)\to G(\Delta')$ where
$\Delta,\Delta'\in\Cal G$.
Then $\Iso(\Cal G)$ is generated by the isomorphisms of type (1)--(4) listed 
below.
\end{thm}

For the following definitions we  make no assumptions on the defining
graph $\Delta$. We first describe three classes of automorphisms.

\begin{itemize}
\item[\rm(1)] {{\bf Graph automorphisms -- $\Aut(\Delta)\,$}\hfill\break 
Any label preserving graph automorphism
of $\Delta$ induces in an obvious way an automorphism of $G(\Delta)$.
We denote by $\Aut(\Delta)$ the group of all such automorphisms.}

\item[\rm(2)] {{\bf Inversion automorphisms -- $\Inv(\Delta)\,$}\hfill\break  
These include the involution 
$\epsilon\co G(\Delta)\to G(\Delta)$ such that $\epsilon(s)=s^{-1}$
for all $s\in V(\Delta)$, which we shall refer to as the \emph{global inversion} of 
$G(\Delta)$, as well as the following involutions which we shall refer to as
\emph{leaf inversions}. For any edge $e=\{s,t\}\in E(\Delta)$ where $t$ is
a terminal vertex and $m_e$ is even, we define the involution 
$\mu_e\co G(\Delta)\to G(\Delta)$ by setting $\mu_e(t)=(sts)^{-1}$ 
and $\mu_e(v)=v$ for all $v\in V(\Delta)\setminus\{ t\}$. 
The global and leaf inversions together generate
a subgroup of $\Aut(G(\Delta))$ isomorphic to $(\Z/2\Z)^{l+1}$, 
where $l$ denotes the number of even labelled 
terminal edges in $\Delta$. We shall denote this subgroup by $\Inv(\Delta)$.}

\item[\rm(3)]{{\bf Inner and Dehn twist automorphisms -- $\Pure(\Delta)\,$}\hfill\break 
Let $T$ denote an edge or vertex of $\Delta$
and suppose that $\Delta=\Delta_1\cup_T\Delta_2$ (by which we imply that 
$\Delta_1$, $\Delta_2$ are full subgraphs of $\Delta$ such that 
$\Delta_1\cup\Delta_2=\Delta$ and $\Delta_1\cap\Delta_2=T$).
Let $g\in C_G(G(T))$ be an element of the centralizer of
$G(T)$. Then we may define an automorphism of $G$ by setting 
\[
\varphi(v)=gvg^{-1} \text{ if }v\in V(\Delta_1),
\ \text{ and } \varphi(v)=v \text{ if }v\in V(\Delta_2)\,.  
\]
Such automorphisms shall be called \emph{Dehn twist}
automorphisms (along $T$). We define $\Pure(\Delta)$ to be the 
subgroup of $\Aut(G(\Delta))$ generated by the Dehn twist automorphisms.
Note that putting $\Delta_2=T=\{s\}$ we obtain the inner automorphism
`conjugation by $s$' as a Dehn twist automorphism. Thus, $\Pure(\Delta)$
contains the group $\Inn(G(\Delta))$ of inner automorphisms of $G(\Delta)$.
By a \emph{nondegenerate} Dehn twist we mean one which is not just 
an inner automorphism, namely a Dehn twist along a separating edge or vertex.}
\end{itemize}

\noindent Note that each nondegenerate Dehn twist is defined in terms of
a ``visual splitting'' of the Artin group, a decomposition as an amalgamated free
product of standard parabolic subgroups, namely 
\[
G(\Delta)=G(\Delta_1)\,\star_{G(T)}\,G(\Delta_2)\,.
\]
 The global and leaf inversions respect any (proper) visual 
splitting of the group while the graph automorphisms carry any visual splitting to a
similar one. Thus graph automorphisms and inversions of $G(\Delta)$ 
conjugate Dehn twist automorphisms to Dehn twist automorphisms. 
Moreover, the graph automorphisms preserve the 
set of even labelled terminal edges and therefore act by conjugation on the inversions. 
Thus $\Aut(G(\Delta))$ contains
a subgroup of the form
\[
\Pure(\Delta)\rtimes\Inv(\Delta) \rtimes \Aut(\Delta)\,.
\]

\begin{remark}
If $e=\{s,t\}\in E(\Delta)$ and $m_e\geq 3$ then the group $G(e)$ has 
infinite cyclic centre
generated by the element $z_e=(st)^k$ where $k=\text{lcm}(m_e,2)/2$.
We also define the element
\[
x_e=\ubrace{ststs\cdots}{m_e}\,.
\]
This element generates the \emph{quasi-centre} of $G(e)$, the subgroup 
of elements which leave the 
generating set $\{s,t\}$ invariant by conjugation. 
We have $z_e=x_e^2$ if $m_e$ is odd and $z_e=x_e$ if $m_e$ is even.

In the case where $G=G(\Delta)$ is a large type (LT) Artin group we can 
explicitly describe the centralizers of separating edges and vertices.
If $e\in E(\Delta)$ then $C_G(G(e))=Z(G(e))=\< z_e\>$. The centralizer of 
a generator $s\in V(\Delta)$ is the direct product of $\< s\>$ with a 
(typically non-cyclic) free group of finite rank. A generating 
set for this free group 
may be obtained by observing that $C_G(\< s\>)/\< s\>$ is isomorphic to 
the vertex group at $s\in V(\Delta)$ in the groupoid with object set $V(\Delta)$
and generated by arrows $x_e\co  r\to r'$ where $e=\{r,t\}$ and $r'=x_e r x_e^{-1}$ 
($r'=r$ if $m_e$ is even, and $t$ otherwise).
We refer the reader to \cite{God1}, or \cite{God2}, for a more detailed description.
\end{remark}

\begin{itemize}
\item[\rm(4)]{{\bf Edge twist isomorphisms}\hfill\break
Suppose that $\Delta=\Delta_1\cup_e\Delta_2$ where $e$ is a
separating edge whose label $m_e$ is odd. Let $\Delta'$ denote
the labelled graph obtained by gluing $\Delta_1$ and $\Delta_2$ together
along the edge $e$ where the identification map reverses the edge.
Then we may define an isomorphism 
\[
\varphi\co G(\Delta)\to G(\Delta')
\]   
by setting 
\[
\varphi(v)=x_evx_e^{-1} \text{ if }v\in V(\Delta_1),
\ \text{ and } \varphi(v)=v \text{ if }v\in V(\Delta_2)\,. 
\]
We shall call such an isomorphism an \emph{edge twist},
and we say that $\Delta$ and $\Delta'$ are \emph{twist equivalent} graphs.
This generates an equivalence relation 
on the set of all  defining graphs. (The collection $\Cal G$ of all CLTTF
defining graphs is invariant under twist equivalence).
Note that in the case where $e=\{ s,t\}$ and $t$ is a terminal vertex of 
$\Delta_1$, then $s$ is a separating vertex and we may think of 
$\Delta$ as the union of $\Delta'_1:=\Delta_1\setminus e$ and 
$\Delta_2$ joined at the vertex $s$. 
In this case the edge twist $\varphi$ modifies the graph $\Delta$ 
by sliding the component $\Delta'_1$ along the edge $e$ so that 
it is attached to $\Delta_2$ at the vertex $t$, instead of at $s$.}
\end{itemize}

\begin{remark} 
The edge twist isomorphism described here is a special case of the ``diagram
twist'' isomorphisms between Artin (and Coxeter) groups first described 
by Brady, McCammond, M\"uhlherr and Neumann in \cite{BMMN}. 
The notion of (diagram) twist 
equivalence as introduced in \cite{BMMN} is defined, more generally, over the
family of all defining graphs and there is considerable evidence for the conjecture
that it is essentially this equivalence relation which classifies all Coxeter groups 
up to isomorphism. A recent survey of the isomorphism problem for Coxeter groups
has been written by M\"uhlherr \cite{Muh}.
\end{remark}

\begin{defn}[{\rm(}Twist equivalence groupoid\/{\rm)}]
Denote $\Biject(\Cal G)$ the groupoid with object set 
$\Cal G =\{\,\text{CLTTF defining graphs}\,\}$ and a morphism $f\co \Delta\to\Delta'$
for each bijection $f\co E(\Delta)\to E(\Delta')$ of the edge sets. 
Observe that every edge twist and every graph automorphism is naturally associated
with a morphism in $\Biject(\Cal G)$. We define $\Twist(\Cal G)$ to be the subgroupoid
of $\Biject(\Cal G)$ generated by the edge twists and graph automorphisms.
\end{defn}

It is known that in any 2--dimensional Artin group the subgroups 
$\< z_e\>$, for $e\in E(\Delta)$, are mutually non-conjugate (this may be readily seen from 
the action of $G(\Delta)$ on its Deligne complex, as described in Section \ref{sect:Deligne}).
Thus the bijection $\ov\varphi\in\Biject(\Cal G)$ induced by any edge twist or graph 
automorphism $\varphi$ is determined by the action of $\varphi$ on the set of conjugacy
classes of the cyclic subgroups $\< z_e\>$ for $e\in E(\Delta)$. Note also that any element
of $\Pure(\Delta)\rtimes\Inv(\Delta)$ acts trivially on this set.
The following statement is largely a consequence 
of Theorem \ref{Thm1} and the above discussion.

\begin{thm}\label{Thm2}
There exists a unique well-defined groupoid homomorphism
\[
\pi\co \Iso(\Cal G)\to \Twist(\Cal G)\,
\]
such that, writing $\pi(\varphi)=\ov\varphi$, we have 
$\<z_{\ov{\varphi}(e)}\>\sim \varphi(\< z_e\>)$ for all $e\in E(\Delta)$.
The image of $\pi$ is $\Twist(\Cal G)$ and the kernel at $\Delta\in \Cal G$ is given by
\[
\ker(\pi,\Delta) = \Pure(\Delta)\rtimes \Inv(\Delta)\,.
\] 
In particular, for fixed $\Delta\in\Cal G$, the automorphism group of $G(\Delta)$
is a (finite) extension of $\Pure(\Delta)\rtimes \Inv(\Delta)$ by a subgroup of
$\Sym(E(\Delta))$ which 
consists of those permutations of $E(\Delta)$ obtained by composing edge twists
and label preserving graph automorphisms. Moreover, two CLTTF Artin groups are
isomorphic if and only if their defining graphs lie in the same connected component
of $\Twist(\Cal G)$, ie, if and only if their defining graphs are twist equivalent.
\end{thm}

Note that the connected components of the groupoids $\Iso(\Cal G)$, 
and $\Twist(\Cal G)$ alike, correspond to the isomorphism classes of 
CLTTF Artin groups. Moreover, the connected components of $\Twist(\Cal G)$ 
are \emph{finite}, and easily computable. Thus, as well as 
determining the automorphism group of any CLTTF Artin group, the above Theorem also 
solves the problem of classifying these groups up to isomorphism. 
In the language of \cite{BMMN}, CLTTF Artin groups are ``rigid up to 
diagram twisting''. Note that spherical type Artin groups (those whose associated
Coxeter groups are finite) are also known to be diagram rigid. This was recently shown 
by Paris in \cite{Paris}. Diagram rigidity is also known for right-angled Artin groups
(the case where all edge labels in $\Delta$ are equal to $2$) by the 
work of Droms \cite{Droms}.
Other partial results on diagram rigidity appear in \cite{BMMN}.

\begin{example}[{\rm(}No separating edges or vertices\/{\rm)}]
Restricting our attention to those CLTTF Artin groups $G=G(\Delta)$ where
$\Delta$ has no separating edge or vertex, we see that two such groups
are isomorphic if and only if their defining graphs are isomorphic, and that
\[
\Aut(G) = \Inn(G)\rtimes (\<\epsilon\> \times \Aut(\Delta))\,.
\]
This is simply because, with no separating edges or vertices, there are no
leaf inversions, nondegenerate Dehn twists or edge twist isomorphisms.
Note that we also have $\Inn(G)\cong G$, since any CLTTF Artin group $G$ has
trivial centre.

A simple example of the above type is where $\Delta$ is the 1-skeleton of a 3-cube 
and all edge labels are $3$. This defining graph also satisfies the vertex rigidity
condition (VR) required by part (ii) of Theorem \ref{Thm3} below. 
\end{example}

\begin{example}[{\rm(}No separating vertices\/{\rm)}]
When $\Delta$ has separating edges but no separating vertices, then the group $\Pure(\Delta)$
is generated by the inner automorphisms and the Dehn twists along separating edges.

A \emph{chunk} of $\Delta$ is a maximal connected full subgraph of $\Delta$ which is not 
separated by the removal of any  edge or vertex which is separating in $\Delta$
(see Section \ref{sect:ChunkEquiv} for a more detailed definition).
Thus if $\Delta$ has no separating vertices it is the union of, say, $N$ distinct chunks 
glued along separating edges.  Fixing a ``base" chunk $B$,
we may suppose that, up to an inner automorphism, each Dehn twist restricts to the
identity on $G(B)$. It can be easily checked that the
 Dehn twists fixing $G(B)$ are mutually commuting elements. 
In this case we therefore have $\Pure(\Delta)\cong G\rtimes \Z^{N-1}$. 
\end{example}

\begin{example}[{\rm(}$\Delta$ a star graph\/{\rm)}]
On the other hand, when there are separating vertices in $\Delta$ we expect the structure of 
$\Aut(G)$ to be somewhat more complicated. For example, one can check that when $\Delta$
is the star graph of $n+1$ vertices ($n$ edges adjoined along a common vertex), and all edge 
labels are $3$ say, then $\Aut(G)$ contains a subgroup isomorphic to 
the $n$--string braid group $B_n$. Let $e_1,..,e_n$ denote the edges of $\Delta$
and, for $i=1,..,n-1$, let $\sig_i$ denote the automorphism of $G$ which is the 
product of the graph automorphism exchanging the edges $e_i$ and $e_{i+1}$ and the Dehn twist 
which conjugates the subgroup $G(e_i)$ by the element $z_{e_{i+1}}$. 
These automorphisms leave invariant the subgroup $F_n$ of $G$ which is 
freely generated by the set $\{ z_e:e\in E(\Delta)\}$ (see Proposition \ref{freegroup}),  
and they describe precisely the standard generators for Artin's representation of the 
braid group as a subgroup of $\Aut(F_n)$. (Moreover, one can check that elements 
of $B_n$ are represented by inner automorphism of $G$ if and only if they are
central in the braid group. Thus $\Out(G)$ is not virtually abelian in this case).   
\end{example}

\subsection*{Abstract commensurators of Artin groups}

We recall that the \emph{abstract commensurator group} $\Comm(\Ga)$ of a group $\Ga$
is defined to be the group of equivalence classes of isomorphisms between finite 
index subgroups of $\Ga$, where two isomorphisms are considered equivalent if 
they agree on common finite index subgroup of their domains. Moreover, two groups 
$\Ga$, $\Ga'$ are said to \emph{abstractly commensurable} if they possess
finite index subgroups $H<\Ga$ and $H'<\Ga'$ which are isomorphic.

\begin{thm}\label{Thm3}
Let $\Delta$ be a CLTTF defining graph with no separating edge or vertex. 
\begin{itemize}
\item[\rm(i)] If $G(\Delta)$ is abstractly commensurable to any CLTTF Artin group
$G(\Delta')$ then $\Delta$ and $\Delta'$ are label isomorphic.
\item[\rm(ii)] Suppose moreover that $\Delta$ satisfies
the following \emph{vertex rigidity} condition:
\begin{description}
\item[\rm(VR)] Any label preserving automorphism of $\Delta$ which
fixes the neighbourhood of a vertex is the identity automorphism.
\end{description}
Then we have 
$\hskip3mm \Comm(G) = \Aut(G) \cong G \rtimes (\<\ep\> \times \Aut(\Delta))\,.$
\end{itemize}
\end{thm}

With regard to part (i) of the above Theorem, we note that a 2--dimensional Artin
group is not commensurable to any other Artin group which is not also 2--dimensional 
(since,
for an Artin group, being 2--dimensional is equivalent to having $\Z\times\Z$ as a 
maximal rank abelian subgroup). We do not know whether the smaller class of CLTTF 
Artin groups is rigid in this sense.

Part (ii) of this Theorem should be compared with \cite{ChCr} where it is shown that
$G$ is commensurable with its abstract commensurator group when $G$ belongs to
one of the two infinite families of 
Artin groups of affine type $\wtil A_n$ and $\wtil C_n$, with $n\geq 2$. 
(The same holds for $G/Z$ where $G$ is an Artin group of finite type $A_n$ or $B_n$, 
with $n\geq 3$, and $Z$ denotes the 
infinite cyclic centre of $G$). 
In Section \ref{sect:Comms} we give an 
example of an abstract commensurator of a CLTTF Artin group $G(\Delta)$ 
which is \emph{not} equivalent to an automorphism in the case where $\Delta$
has no separating edge or vertex, but fails to satisfy the condition (VR).
This hypothesis is therefore necessary. Examples are also given of  
abstractly commensurable but non-isomorphic CLTTF Artin groups.

\subsection*{Isomorphisms  of Coxeter groups}

Finally we consider isomorphisms between Coxeter groups of CLTTF type.
Let $\Iso_W(\Cal G)$ denote the category (a groupoid) with objects $\Cal G$
and morphisms the isomorphisms $W(\Delta)\to W(\Delta')$ for 
$\Delta,\Delta'\in\Cal G$. We note (by inspection of the isomorphisms of type (1)--(4))
that every isomorphism $\varphi\co G(\Delta)\to G(\Delta')$ induces  an isomorphism
$\varphi_W\co W(\Delta)\to W(\Delta')$. This is natural in the sense that
$\varphi_W\circ\rho_\Delta=\rho_{\Delta'}\circ\varphi$,
where $\rho_\Delta\co G(\Delta)\to W(\Delta):g\mapsto \ov g$ denotes the canonical surjection. 
Thus the mapping $\varphi\mapsto\varphi_W$ defines a groupoid homomorphism 
\[
\rho\co  \Iso(\Cal G)\to \Iso_W(\Cal G)\,.
\]

\begin{remark}
The above remarks imply, in particular, that the \emph{pure Artin group} $PG(\Delta)$, 
which is defined as the kernel of the canonical quotient 
$\rho_\Delta\co G(\Delta)\to W(\Delta)$, is a characteristic subgroup of $G(\Delta)$
for CLTTF type Artin groups. This agrees with results already known for irreducible 
finite type Artin groups by Cohen and Paris \cite{CohPar} which generalised  
a much earlier Theorem of Artin \cite{Artin} in the case of the braid groups.
\end{remark}

There is a further source of Coxeter group automorphisms \emph{not} 
induced from automorphisms of the associated Artin groups. These shall be thought of as 
``pure" automorphisms since, as with the inner and Dehn twist automorphism 
(induced from $\Pure(\Delta)$), they respect the conjugacy class of the element 
$\ov x_e=\rho_\Delta(x_e)$, for each $e\in E(\Delta)$. 

\paragraph{Pure automorphisms of $W(\Delta)$}
Let $e=\{ s,t\}\in E(\Delta)$ denote a \emph{cut edge}: every edge path in $\Delta$ 
from $s$ to $t$ passes through $e$. Then there are disjoint connected full subgraphs
$\Delta_1,\Delta_2$ of $\Delta$ such that $\Delta = \Delta_1\cup e \cup \Delta_2$ 
with $\Delta_1\cap e=\{s\}$ and $\Delta_2\cap e=\{t\}$. Let $m=m_e\geq 3$, 
and let $r\in\N$ such that $2r+1$ is congruent (mod $m$) to a unit in the 
ring $\Z/m\Z$. Then we may define an automorphism of $W(\Delta)$ by setting 
\[
\varphi(v)=(st)^rv(st)^{-r} \text{ if }v\in V(\Delta_1),
\ \text{ and } \varphi(v)=v \text{ if }v\in V(\Delta_2)\,.  
\]
Such automorphisms shall be called \emph{dihedral twist}
automorphisms. We define $\Pure_W(\Delta)$ to be the 
subgroup of $\Aut(W(\Delta))$ generated by all dihedral twists, 
Dehn twists and inner automorphisms. In particular, $\Pure_W(\Delta)$ 
contains all automorphisms induced from $\Pure(\Delta)$.\\

The following Theorem gives a solution to the ``classification'' and ``automorphism''
problems for CLTTF Coxeter groups. We remark that the classification up to isomorphism
is already contained in the work of M\"uhlherr and Weidmann \cite{MW} 
on reflection rigidity and 
reflection independance in  large type (what they call ``skew-angled") Coxeter groups. 
Also, the automorphism groups have already been determined in many of the cases covered 
here (and some besides) by Bahls \cite{Bahls}.
The proof of Theorem \ref{ThmCoxeter} which we give 
consists in repeating the same sequence of 
arguments used to establish Theorems \ref{Thm1} and \ref{Thm2}, 
with appropriate slight modification, 
in the context of Coxeter groups.
  
\begin{thm}\label{ThmCoxeter}
The groupoid $\Iso_W(\Cal G)$ is generated
by pure automorphisms (elements of $\Pure_W(\Delta)$, for $\Delta\in \Cal G$),
graph automorphisms and edge twist isomorphisms.
More precisely, there is a surjective groupoid homomorphism
\[
\pi_W\co \Iso_W(\Cal G)\to \Twist(\Cal G)\,
\]
with $\pi_W\circ\rho=\pi$, and for each $\Delta\in\Cal G$ we have 
\[
\ker(\pi,\Delta) = \Pure_W(\Delta)\,.
\]
In particular, for fixed $\Delta\in\Cal G$, the automorphism group of $W(\Delta)$
is a (finite) extension of $\Pure_W(\Delta)$ by the subgroup of
$\Sym(E(\Delta))$ appearing as a vertex group in $\Twist(\Cal G)$.
Moreover, two CLTTF Coxeter groups are
isomorphic if and only if their defining graphs lie in the same connected component
of $\Twist(\Cal G)$, ie, if and only if their defining graphs are twist equivalent.
\end{thm}

The automorphism group of a CLTTF Coxeter group has previously been described
by Patrick Bahls \cite{Bahls} under the added hypotheses that 
all edge labels are even and the defining graph cannot be separated into more than
2 components by removal of a single edge. In fact, in his work, Bahls does not suppose
that the defining graph is triangle free, and so treats many cases which are not 
covered here. He also gives several statements (see Corollaries 1.2, 1.3, 1.4 of \cite{Bahls})
giving further details on the size and structure of $\Out(W)$ which probably extend
to the CLTTF case. 

As an example, consider the case where $\Delta$ has no separating 
vertices. In this case there are no dihedral twists 
and $\Pure_W(\Delta)\cong W(\Delta)\rtimes (\Z/2\Z)^{R-1}$, where
$R$ is the number of distinct maximal full subgraphs of $\Delta$ not separated by any 
\emph{even} labelled edge (compare with Corollary 1.3 in \cite{Bahls}). 
In particular, $\Out(W)$ is finite in this case.  Note, however, that the corresponding 
Artin groups have typically infinite outer automorphism groups. 
In the case of no separating vertices we have already seen that 
$\Pure(\Delta)\cong G(\Delta)\rtimes (\Z)^{N-1}$ with $N\geq R$. 

 Recently, M\"uhlherr and Weidmann \cite{MW} have proved results on reflection
rigidity and reflection independance in the wider class of large type (LT) 
Coxeter groups which give the same solution to the classification problem 
as given by Theorem \ref{ThmCoxeter} above. We note  that Bahls \cite{Bahls2} has
also obtained a similar classification for those
Coxeter groups having 2--dimensional Davis complex (equivalently, 
those associated to 2--dimensional Artin groups).
Several other results in this direction are discussed in the survey
by M\"uhlherr \cite{Muh}.
It seems reasonable to conjecture that Theorems \ref{Thm1}, \ref{Thm2}
and  \ref{ThmCoxeter} all hold unchanged over the class of connected 
large type (CLT) defining graphs, and that similar results might also hold 
for all 2--dimensional Artin groups, or for general Coxeter groups.

\begin{ackn}
This work has benefitted from discussions with many people. In particular,
I would like to thank Benson Farb, Luisa Paoluzzi, 
Bernhard M\"uhlherr, Patrick Bahls, Gilbert Levitt and the referee
for a variety of helpful suggestions and comments.  
\end{ackn}


\section{The Deligne complex $\D$}\label{sect:Deligne}

For simplicity, we formulate the following definitions only in the case where
the Artin group $G=G(\Delta)$ is \emph{2--dimensional}, equivalently, where
every triangle in $\Delta$ with edge labels $m,n,p$ satisfies $1/m+1/n+1/p\leq 1$. 
See \cite{CD} for details of the general construction.
 
\begin{figure}[ht!]\small\anchor{Fig:delignecx}
\cl{
\psfrag {e}{$e$}
\psfrag {f}{$f$}
\psfrag {t}{$t$}
\psfrag {s}{$s$}
\psfrag {r}{$r$}
\psfrag {D}{$\Delta$}
\psfrag {K}{$K$}
\psfrag {Ve}{$V_e$}
\psfrag {Vf}{$V_f$}
\psfrag {Vs}{$V_s$}
\psfrag {Vt}{$V_t$}
\psfrag {Vr}{$V_r$}
\psfrag {V0}{$V_\emptyset$}
\includegraphics[width=10cm]{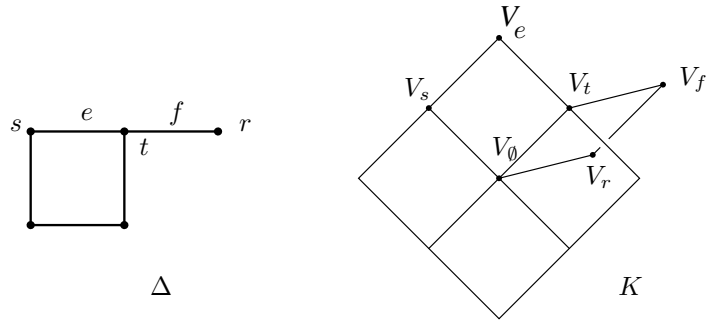}}
\caption{Defining graph $\Delta$ and squared complex $K$ for 
a 2--dimensional Artin group}
\label{Fig:delignecx}
\end{figure}

\paragraph{Definition of the Deligne complex $\D\,$}
Let $K$ denote the geometric realisation of the derived complex of 
the partially ordered set 
\[
\{ V_\emptyset\}\cup\{V_s:s\in V(\Delta)\}\cup\{ V_e:e\in E(\Delta)\}\,,
\]
where the partial order is given by setting $V_\emptyset <V_s$ 
for all $s\in V(\Delta)$, and $V_s < V_e$ whenever $s$ is a vertex of the edge $e$. 
Thus $K$ is a finite 2--dimensional simplicial complex. We may also view $K$ as a 
squared complex with one square cell for each edge of $\Delta$.
If $e=\{s,t\}\in E(\Delta)$ then the corresponding square cell has vertices
$V_\emptyset,V_s,V_t,V_e$. We note that, viewing $K$ as a squared complex 
in this way we have $\Lk(V_\emptyset,K)\cong\Delta$. See \figref{Fig:delignecx}.

Let $\Cal K$ denote the complex of groups with underlying complex $K$
and vertex groups $G(V_\emptyset)=\{ 1\}$, $G(V_s)=\<s\>=G(s)$, 
for $s\in V(\Delta)$, and
$G(V_e)=G(e)$, for $e\in E(\Delta)$. Then $\Cal K$ is a developable
complex of groups (cf \cite{CD}) whose
fundamental group is the Artin group: $\pi_1(\Cal K)=G(\Delta)$.

\begin{defn}[{\rm(}Deligne complex\/{\rm)}]
Let $G=G(\Delta)$ be a 2--dimensional Artin group.
We define the \emph{Deligne complex} $\D$, of type $\Delta$,
to be  the universal covering $\wtil{\Cal K}$ of the complex of groups
$\Cal K$ just described, equipped with the action of $G$ by covering 
transformations.
\end{defn}

The Artin group acts by simplicial isomorphisms of $\D$
with vertex stabilizers
either trivial or conjugate to one of the standard parabolic subgroups
$G(s)$, for $s\in V(\Delta)$, or $G(e)$, for $e\in E(\Delta)$.
We classify the vertices of $\D$ into three kinds 
according to their stabilizers:

\begin{description}
\item[Rank 0] 
vertices of the form $gV_\emptyset$ for $g\in G$. These have trivial stabilizer.
\item[Rank 1] 
vertices $gV_s$ for $s\in V(\Delta)$ and $g\in G$ --- $Stab(gV_s)=g\<s\>g^{-1}$.
\item[Rank 2] 
vertices $gV_e$ for $e\in E(\Delta)$ and $g\in G$ --- $Stab(gV_e)=gG(e)g^{-1}$.
\end{description}
 
Note that every point in the open neighbourhood of a rank 0 vertex represents
a free orbit of the group action (since the group action is strictly cellular).

We also note that an analogous construction replacing the vertex groups of 
$\Cal K$ with the corresponding finite standard  parabolic subgroups of the 
Coxeter group $W$ results in a description of the \emph{Davis complex},
which we shall denote by $\D_W$. There
is a natural simplicial map $p_W\co \D\to \D_W$ induced by the canonical projection
$G\to W$ and an inclusion $i_W\co \D_W\hookrightarrow \D$ induced by the Tits section
$W\hookrightarrow G$. We have $p_W\circ i_W$ equal to the identity on $\D_W$.

\begin{defn}[{\rm(}Metrics on $\D$\/{\rm)}]
There are two natural choices of $G$--equivariant piecewise Euclidean metric
for the complex $\D$. The first, and perhaps most natural, is known as
the \emph{Moussong metric} and is defined such that, for $e=\{s,t\}\in E(\Delta)$, 
the simplex $(V_\emptyset,V_s,V_e)$ is a Euclidean triangle with 
angles $\frac{\pi}{2}$ at
$V_s$ and $\frac{\pi}{2m_e}$ at $V_e$. (See \cite{Mou}, also \cite{CD}).
The Moussong metric on $\D$ is known to be CAT(0) for all 2--dimensional Artin groups.
This property will be used in Section \ref{sect:CNVA}. 

The second is the \emph{cubical metric} obtained by viewing $\D$ as a 
squared complex (as in \figref{Fig:delignecx}) built from unit Euclidean squares.
For $G(\Delta)$ 2--dimensional, the cubical metric on $\D$ is known to be CAT(0) 
if and only if $\Delta$ is triangle free (see \cite{CD}). In particular, this metric is
CAT(0) in the CLTTF case. The cubical metric shall be used in Section \ref{sect:TwoProps}.

We note that each of these metrics induces a unique metric on the Davis complex $\D_W$ 
such that the map $i_W\co \D_W\hookrightarrow \D$ is an isometric embedding.
\end{defn}

The following definition and lemma will be relevant in Section \ref{sect:CNVA}.

\begin{defn}[{\rm(}Hyperbolic type\/{\rm)}]
We shall say that a defining graph $\Delta$, or the associated Artin group $G(\Delta)$,
is of \emph{hyperbolic type} if the Coxeter group $W(\Delta)$ is a Gromov 
hyperbolic group. Equivalently, $\Delta$ is of hyperbolic type if and only if the Davis 
complex $\D_W$ is a $\delta$--hyperbolic metric space with respect to either the 
Moussong metric or the cubical metric. (This is because the Coxeter group acts properly 
and co-compactly by isometries with respect to either metric on the Davis complex
and so is quasi-isometric to  $\D_W$).
\end{defn}

\begin{lemma}\label{Dhyp}
Let $G(\Delta)$ be a 2--dimensional Artin group. Then the following are equivalent
\begin{itemize}
\item[\rm(1)] $G(\Delta)$ (or $\Delta$) is of hyperbolic type;
\item[\rm(2)] the Deligne complex $\D$ equipped with the Moussong metric 
is a $\delta$--hyperbolic metric space;
\item[\rm(3)] $\Delta$ contains no triangle having edge labels $m,n,p$ with 
$1/m+1/n+1/p =1$ and no square with all edge labels equal to $2$. 
\end{itemize}
\end{lemma}

\begin{proof}
We suppose throughout that the Deligne complex $\D$ is equipped with the Moussong metric.
Since $G(\Delta)$ is 2--dimensional, this implies that $\D$ is a CAT(0) space.
By the Flat Plane Theorem (see \cite{BH}) this space is $\delta$--hyperbolic 
if and only if it contains no isometrically embedded flat plane $\E^2$. 
If such a plane existed in $\D$, it would necessarily be a simplicial subcomplex 
and so contain at least one rank 0 vertex. Morever, it would contribute a simple
circuit of length exactly $2\pi$ to the link of any such vertex.
We note that the link of a rank 0 vertex of $\D$ contains a circuit of
length exactly $2\pi$ if and only if there exists either a triangle in $\Delta$ 
with labels $m,n,p$ such that $1/m+1/n+1/p=1$, or a square in $\Delta$ with 
all labels $2$. Thus, condition (3) implies that no embedded flat plane can occur
in $\D$, and hence that (2) holds. On the other hand if (3) fails then $W(\Delta)$
contains either a Euclidean triangle group, or $D_\infty\times D_\infty$. In either
case $W(\Delta)$ contains a subgroup $\Z\times \Z$, and so cannot be Gromov hyperbolic.
Thus, we have shown (1) implies (3), as well as (3) implies (2).

Finally, we use the fact that the Davis complex $\D_W$ (with the Moussong metric)  
embeds isometrically in $\D$. Any flat plane in $\D_W$ is therefore also
a flat plane in $\D$. Thus, by the Flat Plane Theorem, $\D_W$ is $\delta$--hyperbolic 
if $\D$ is, and so (2) implies (1).
\end{proof}
 
We note that any CLTTF Artin group is necessarily a 2--dimensional Artin group
of hyperbolic type. Similarly, any CLTTF Coxeter group has 2--dimen\-sional 
$\delta$--hyperbolic Davis complex.

In  Section \ref{sect:CNVA} we shall also use the following statement which is a 
consequence of a quite general result  due to Bridson \cite{Bri}. We recall
that an isometry $\ga$ of a geodesic metric space $X$ 
is said to be \emph{semi-simple} if it
attains its translation length: $|\ga|:=\textsl{inf}\{d(x,\ga x):x\in X\}$ 
is realised at some point in $X$. Semi-simple elements are 
classified into two classes: \emph{elliptic} if $|\ga|=0$; and \emph{hyperbolic} if
$|\ga|\neq 0$. Bridson's result in  \cite{Bri} states that any isometry of a geodesic
metric simplicial complex having finitely many isometry types of cells is 
necessarily semi-simple. As a consequence we have:

\begin{lemma}\label{Daction}
Let $G$ be a 2--dimensional Artin group. Then the action of $G$ on $\D$ is 
semi-simple (with respect to either the Moussong metric, or the cubical metric).
\qed
\end{lemma}


\section{Structure of vertex stabilisers and fixed sets in $\D$}

We consider the 2--generator Artin groups which appear as the stabilizers of
rank 2 vertices of the Deligne complex $\D$ associated to a 2--dimensional 
Artin group and derive some basic properties which will be useful in the sequel.
Using one of these properties, we also give a classification of the fixed sets 
in $\D$ for arbitrary elements of a 2--dimensional Artin group.

Recall that if $e=\{s,t\}\in E(\Delta)$ with label $m_e$ then the group $G(e)$
is given by the presentation
\[
G(e)=\<\ s,t\  \mid\ \ubrace{ststs\cdots}{m_e}=\ubrace{tstst\cdots}{m_e}\ \>
\]
When $m_e\geq 3$ the centre of $G(e)$ is infinite cyclic generated by the 
element  $z_e:=(st)^k$ where $k=m_e$ if $m_e$ is odd, and $k=m_e/2$
if $m_e$ is even.  (Alternatively $k=\text{lcm}(m_e,2)/2$.) We wish to consider
the quotient of $G(e)$ by its centre, which we shall denote by
\[
\Gamma\, =\, G(e)/\<z_e\>\,.
\]
We shall systematically write $\ov x$ for the image in $\Gamma$ of an element 
$x\in G(e)$. 

We note that $\Gamma$ is a virtually free group (and virtually 
cyclic if and only if $m_e=2)$. In fact, 
\[
\Gamma\,\cong\,\left\{ \begin{aligned}
\Z_2\star \Z_k\quad &\text{if } m_e \text{ odd,}\\
\Z\star \Z_k\quad &\text{if } m_e \text{ even.}\\
\end{aligned}
\right.
\]
The free factors here are generated by the elements 
$\ov{st}$ of order $k=\text{lcm}(m_e,2)/2$ and either 
$\ov x_e$ of order $2$ when $m_e$ is odd,
or $\ov s$ of infinite order in the case $m_e$ even.

It is clear from the the above description that, in each case, 
$\Gamma$ admits a proper co-compact action on a regular $k$--valent metric 
tree $T$ (with edge lengths equal to $1$) where 
the fixed set of any elliptic element consists 
of a single point -- elements conjugate to $\ov{st}$ fix the vertices and, in
the case $m_e$ odd, elements conjugate to $\ov x_e$ fix the midpoints of edges.
On the other hand, both generators $s$ and $t$ of $G(e)$  
act by hyperbolic isometries of $T$ of translation length $1$. 
(These actions are described in more detail in Section 2 of \cite{BrCr},
for example).

Finally, we note that any Artin group $G(\Delta)$ admits a standard \emph{length} 
homomorphism $\ell\co G(\Delta)\to\Z$ defined by setting $\ell(s)=1$ 
for all $s\in V(\Delta)$. 
  
\begin{lemma}\label{structureGe}
Let $G(e)$ be the rank 2 Artin group associated to an edge $e=\{ s,t\}$, with label
$m_e\geq 2$. Let $R$ denote the set of all elements conjugate in $G(e)$ 
into the generating set $\{ s,t\}$, and let $x\in G(e)$. Then
\begin{itemize}
\item[\rm(i)] $C_{G(e)}(\<x\>)$ is virtually abelian if and only if 
either $m_e=2$ or $m_e\geq 3$ and $x$ is not central.
\item[\rm(ii)] Let $u\in R$ and $k\in\Z\setminus\{ 0\}$. Then 
$C_{G(e)}(\<u^k\>)=\< u,z_e\>\cong\Z\times\Z$, and if $x^k=u^k$ then $x=u$.
\item[\rm(iii)] Suppose $m_e\geq 3$. Let $u,v\in R$ 
and $k,l\in\Z\setminus\{ 0\}$. If $u^k$ and $v^l$ commute then $u=v$. 
\end{itemize}
\end{lemma} 

\begin{proof}
We shall suppose throughout that $m_e\geq 3$, the case  where
$m_e=2$ and $G(e)\cong \Z^2$ being easily checked.
\medskip

\noindent\textbf{(i)}\qua
If $x$ lies in the centre $Z(G(e))$ then $C_{G(e)}(\< x\>)=G(e)$ which
is not virtually abelian (since $m_e\geq 3$).
On the other hand, if $x\notin Z(G(e))$ then its image
$\ov x$ in $\Ga$ is nontrivial. We consider the action of $\ov x$ on the tree $T$.
If $\ov x$ is elliptic then its fixed set consists of a single point $p$. But
then $C_\Ga(\ov x)$ fixes $p$, so must be finite. If $\ov x$ is hyperbolic then
$C_\Ga(\ov x)$ leaves invariant its axis. In either case $\ov x$ generates 
a finite index subgroup of $C_\Ga(\ov x)$. Therefore $x$ and $z_e$ generate a finite
index abelian subgroup of $C_{G(e)}(x)$.
\medskip

\noindent\textbf{(ii)}\qua
Since $u$ is conjugate to a generator, $\ov{u}$ is hyperbolic on $T$ with
translation length $|\ov u|=1$. Let $A\subset T$ denote the translation axis
for $\ov u$. This is also the unique translation axis for each power of $\ov u$.
Let $x\in C_{G(e)}(\<u^k\>)$. Then, since $\ov x$  commutes with $\ov u^k$, 
it leaves invariant the axis $A$ (without reversing its direction).
Since $\ov u$ has unit translation length we can find $n\in\Z$ such that $\ov x$ and $\ov u^n$ 
differ by an elliptic fixing the whole axis $A$ and, since the fixed set of any elliptic
in $\Gamma$ is a single point in $T$, we have that $\ov x=\ov u^n$. It follows that
$x\in \< u,z_e\>$. Thus $C_{G(e)}(\< u^k\>)=\< u,z_e\>\cong\Z\times\Z$.
 
If $x^k=u^k$ then $x$ centralizes $u^k$, and so $x\in \< u,z_e\>\cong\Z\times\Z$.
But now we have $x=u$ since uniqueness of roots holds in a free abelian group. 
\medskip

\noindent\textbf{(iii)}\qua
Since they are conjugate to generators, $u$ and $v$ project to 
hyperbolic isometries $\ov u$ and $\ov v$ of $T$ with translation length $1$ in 
each case. If $u^k$ and $v^l$ commute for nonzero $k$ and $l$ then $\ov u$ 
and $\ov v$ must also share an axis in $T$. But then $\ov{u} =\ov{v}^{\pm 1}$,
and one of $uv^{-1}$ or $uv$ lies in the centre $\< z_e\>$. 
But since $\ell(z_e) =\text{lcm}(m_e,2)\geq 3$, while $\ell(u)=\ell(v)=1$ 
we deduce that $u=v$.
\end{proof}

We now consider the action of a 2--dimensional Artin group $G=G(\Delta)$ on
its Deligne complex $\D$.

\begin{defn}[{\rm(}Fixed sets and $F_s$\/{\rm)}]
For $g\in G$ we write $\Fix(g)$ for the (possibly empty) set of points in $\D$
left fixed by $g$. If $s\in V(\Delta)$ we write $F_s$ for the fixed 
set $\Fix(s)$ of $s$.
\end{defn}

Note that $F_s$ is necessarily a  geodesically convex subcomplex
of $\D$. Moreover, since rank 0 vertices have trivial 
stabilizer, $F_s$  lies in that part of the 1-skeleton of $\D$ which is
spanned by rank 1 and 2 vertices. Consequently $F_s$ is a tree 
(since it is geodesically convex) 
whose vertices are alternately rank 1 and 2 vertices of $\D$.

\begin{lemma}\label{fixedsets}
Suppose that $G=G(\Delta)$ is a 2--dimensional Artin group, 
and let $x\in G\setminus\{ 1\}$.
\begin{itemize}
\item[\rm(i)] If $x\in \< s\>$, for $s\in V(\Delta)$, then $\Fix(x)=F_s$.
\item[\rm(ii)] If $x\in G(e)$, for $e=\{ s,t\}\in E(\Delta)$, but $x$ is 
not conjugate in $G(e)$ into $\<s\>$ or $\<t\>$, then $\Fix(x)=\{V_e\}$.
\item[\rm(iii)] If $x$ is not conjugate in $G$ to any of the elements covered by
cases (i) and (ii) above, then $\Fix(x)=\emptyset$. 
\end{itemize}  
\end{lemma}

\begin{proof}
\noindent\textbf{(i)}\qua
Let $x=s^k$ for some $k\neq 0$. Clearly $F_s\subset \Fix(s^k)$. 
If $\Fix(s^k)\neq F_s$ then there must be some edge $g[V_t,V_e]$ of $\D$
($g\in G$, $e=\{ t,t'\}\in E(\Delta)$) which is fixed by $s^k$ but only one 
of whose vertices is fixed by $s$. If $s$ fixes $gV_t$ then it also fixes $gV_e$ 
(since $G(t)<G(e)$) so we may suppose that $s$ fixes $gV_e$ but not $gV_t$. 
Then, writing $y= g^{-1}sg$, we have $y\in G(e)$. On the other hand, since $s^k$ fixes 
$gV_t$ we have that $y^k\in\< t\>$. Comparing lengths, we must have $y^k=t^k$
and therefore $y=t$, by Lemma \ref{structureGe} (ii). But then $s$ fixes the vertex
$gV_t$ contrary to the choice of edge. Thus $\Fix(s^k)=F_s$.
\medskip

\noindent\textbf{(ii)}\qua
If $x\in G(e)$ then it clearly fixes the point $V_e$ in $\D$.
If, however, $\Fix(x)$ contains any other vertex of $\D$ then it contains 
a neighbouring vertex, that is $gV_t$ or $gV_s$ for some $g\in G(e)$. 
But that is to say that $x$ is conjugate, in $G(e)$, into one
of the subgroups $\<s\>$ or $\<t\>$.
\medskip

\noindent\textbf{(iii)}\qua
 If $\Fix(x)\neq\emptyset$ then $x$ must fix some rank 2 vertex
(if it fixes a rank 1 vertex then it fixes every neighbouring rank 2 vertex).
But then $x$ is conjugate to $x'\in G(e)$ for some edge $e\in E(\Delta)$ and if 
$x'$ is not covered by case (ii) it is conjugate to an element covered by case (i).
\end{proof}


\section{CNVA  subgroups and their fixed sets in $\D$}\label{sect:CNVA}

\begin{defn}
Let $C$ denote a nontrivial (necessarily infinite) cyclic subgroup of $G$.
We say that $C$ is \emph{CNVA  (``centralizer not virtually abelian'') in $G$}
if its centralizer $C_G(C)$ is not virtually abelian.
\end{defn}

Note that if $H$ is a finite index subgroup of $G$ and $C<H$,
then $C$ is  CNVA  in $G$ if and only if it is  CNVA  in $H$ (ie $C_H(C)$ 
is not virtually abelian). The property of being  CNVA  is also inherited by
subgroups of $C$, for if $C'<C$ then the centralizer
$C_G(C')$ contains $C_G(C)$ and so fails to be virtually
abelian unless $C_G(C)$ is virtually abelian.

\begin{defn}[{\rm(}Internal vertex\/{\rm)}] Let $\Delta$ be an Artin defining graph. 
By an \emph{internal vertex} of $\Delta$ we mean a vertex of valence at least two.
\end{defn}

\begin{lemma}\label{star1}
Let $G=G(\Delta)$ be a 2--dimensional Artin group.
\begin{itemize}
\item[\rm(i)] If $e\in E(\Delta)$ with $m_e\geq 3$ then each nontrivial subgroup
of $\< z_e\>$ is CNVA.
\item[\rm(ii)] If $s\in V(\Delta)$ is an internal vertex then each nontrivial subgroup
of $\< s\>$ is CNVA.
\item[\rm(iii)] Suppose that $s\in V(\Delta)$ is not conjugate in $G(\Delta)$ to any generator 
corresponding to an internal vertex of $\Delta$. Then NO nontrivial subgroup of $\< s\>$ is CNVA.
\end{itemize}
\end{lemma}

\begin{proof}

{\bf (i)}\qua The fact that $C_G(\< z_e\>)=G(e)$ is virtually nonabelian free by cyclic
when $m_e\geq 3$ ensures that $\< z_e\>$ (and each of its subgroups)
is  CNVA  for all $e\in E(\Delta)$ with $m_e\geq 3$.
\medskip

\noindent{\bf (ii)}\qua Let $s\in V(\Delta)$. We consider the tree $F_s$ lying in
the 1-skeleton of $\D$ which is the fixed
point set of $s$. This tree is left invariant by $C_G(\< s\>)$, and we may
therefore consider the action of the centralizer on $F_s$.
We note that the rank 1 vertex $V_s$ lies in $F_s$, and that
$\Stab(V_s)=G(s)=\< s\>$. Any vertex of $\D$ which is adjacent to $V_s$
is $V_e$ for some $e\in E(\Delta)$ such that $s$ is a vertex of $e$.

Suppose now that $s$ is internal.
Then  there are rank 2 vertices $V_e$ and $V_d$
which lie in $F_s$, for distinct edges $e,d$ adjacent to $s$.
(Note that the vertex $V_s$ lies midway between $V_e$ and $V_d$). 
The element $z_e$ (resp. $z_d$) centralizes $s$ and fixes $V_e$ (resp. $V_d$)
but does not fix the point $V_s$. Since the elements $z_e$ and $z_d$
are acting in this way on a tree they necessarily generate a free group
of rank 2 inside $C_G(\< s\>)$, implying that $\< s\>$ 
(and hence $\< s^k\>$ for any $k\neq 0$) is CNVA. 
\medskip

\noindent{\bf (iii)}\qua
We note that if $s$ belongs to an \emph{odd} labelled edge $e=\{s,t\}$ 
then $s$ is conjugate to $t$ (by the element $x_e$).
It  follows that there are exactly three ways that $s$ can fail to be conjugate to
an internal vertex generator (we have not supposed here that $\Delta$ is connected). Either
\begin{itemize}
\item[(a)] $s$ is an isolated vertex of $\Delta$, or
\item[(b)] $s$ lies in a component of $\Delta$ which consists of a single edge $e$, or 
\item[(c)] $s$ lies in a unique edge $e$, and $m_e$ is even.  
\end{itemize}
We recall that, in general, the fixed set $F_s=\Fix(s)$ is a connected 1-dimensional 
subcomplex of $\D$, in fact a tree, whose vertices are alternately vertices of 
rank 1 and 2. 
Recall also that $F_s=\Fix(s^k)$, for all $k\geq 1$, by Lemma \ref{fixedsets}(i).
We use the basic fact that the centralizer of any element $g$ 
must leave invariant the set $\Fix(g)$. 
Thus $C_G(\< s^k\>)$ leaves $F_s$ invariant, for all $k\geq 1$.

In case (a), $F_s$ consists solely of the vertex $V_s$, since this vertex is not 
adjacent in $\D$ to any rank 2 vertex at all. In this case, 
$C_G(\< s^k\>)$ must fix $V_s$ and is therefore an infinite 
cyclic group (since $\Stab(V_s)=\<s\>$). 
Thus, in case (a), $\< s^k\>$ fails to be CNVA, for all $k\geq 1$.

In cases (b) and (c) we claim that $F_s$ is a bounded (but still infinite) 
tree containing exactly one rank 2 vertex, namely the vertex $V_e$. 
First note that any rank 1 vertex of $F_s$ can be written $hV_t$ where $h\in G$
and $t\in V(\Delta)$ is a generator which is conjugate to $s$ 
(in fact we must have $s=hth^{-1}$ because 
$\< s\>\leq \Stab(hV_t)=h\< t\> h^{-1}$ and $\ell(s)=\ell(t)=1$). 
Moreover, any edge of $F_s$ may be written $h[V_t,V_f]$ 
for some $h\in G$, some  $t$ conjugate to $s$,
and some $f\in E(\Delta)$ such that $\ t\in f$. 

Now observe that, in both cases (b) and (c), there exists a homomorphism 
$\nu\co G\to\Z$ such that $\nu(t)=0$  if $t$ lies in some edge different from $e$, 
and $\nu(t)=1$ otherwise. In particular $\nu(s)=1$, and clearly $\nu(t)=1$ 
for any generator $t$ which is conjugate to $s$.
This shows that $e$ is the only edge which can possibly contain a vertex $t$
such that $s$ and $t$ are conjugate. It follows that every edge of $F_s$ is a translate of 
$[V_t,V_e]$ for some $t\in e$, and in particular that every rank 2 vertex is a translate of 
$V_e$. However, since each rank 1 vertex in $\D$ can be adjacent to at most one translate 
of a given rank 2 vertex, it now follows (by connectedness) that $F_s$ lies entirely in the 
neighbourhood of the vertex $V_e$.  

Since it leaves $F_s$ invariant, we deduce in cases (b) and (c) that $C_G(\< s^k\>)$ 
must fix $V_e$ (the unique rank 2 vertex of $F_s$), and hence is a subgroup of $G(e)$,
for all $k\geq 1$. 
But then $\< s^k\>$ is not  CNVA  since, by Lemma \ref{structureGe}(i), it 
has virtually abelian centralizer in $G(e)$.
\end{proof}

\begin{remarknum}\label{rem:unboundedtrees}
It is implicit in the above proof that, for $s\in V(\Delta)$, the
cyclic group $\<s\>$ is CNVA if and only if its fixed set $F_s$ is an unbounded 
tree. 
\end{remarknum}

\begin{lemma}\label{starlike}
Let $G=G(\Delta)$ be a 2--dimensional Artin group of hyperbolic type.
A cyclic subgroup of $G$ is  CNVA  if and only if it is conjugate in $G$ to either 
a subgroup of $\< z_e\>$ for some $e\in E(\Delta)$ with $m_e\geq 3$, or a subgroup
of $\< s\>$ for some internal vertex $s\in V(\Delta)$.
\end{lemma}

\begin{proof}
By Lemma \ref{star1}, it will suffice to show that any CNVA subgroup 
is conjugate into a subgroup of  either $\< z_e\>$, for some $e\in E(\Delta)$ 
with $m_e\geq 3$, or $\< s\>$, for some $s\in V(\Delta)$.
 
We suppose for the purposes of this proof that $\D$ is equipped with the
Moussong metric, and so is CAT(0) by \cite{CD}.
Suppose that $C$ is a  CNVA  cyclic subgroup of $G$ generated by
the element $\ga$.
By Lemma \ref{Daction} and the classification of semi-simple
isometries, this element is either elliptic or hyperbolic.

Assume firstly that $\ga$ is elliptic. By Lemma \ref{fixedsets}, 
either $\ga$ fixes some rank 1 vertex, and so is conjugate 
into $\< s\>$ for some $s\in V(\Delta)$ as required, or
$\Fix(\ga)=\{ gV_e\}$ for some $g\in G$ and $e\in E(\Delta)$.
In the latter case the centralizer $C_G(\< \ga\>)$  must also fix
the vertex $gV_e$ and so is a subgroup of $\Stab(gV_e)=gG(e)g^{-1}$.
But then, by Lemma \ref{structureGe}(i), it follows that
 $\ga$ is an element of $g\< z_e\>g^{-1}$ and $m_e\geq 3$,
since otherwise it would have virtually abelian centralizer.
 
We now assume that $\ga$ is hyperbolic. Let $M$ denote the minset of $\ga$.
Then by, Theorem II.6.8 of \cite{BH},  $M\cong T\times\R$ where $T$ is, 
in our case, a metric tree. However, $T$ must be a bounded tree, since otherwise we
would have a flat plane $\E^2$ isometrically embedded in $\D$, contradicting
Lemma \ref{Dhyp} (with the hypothesis that $\Delta$ is hyperbolic type).
Therefore, $T$ has a fixed point $c$ under the action of $C_G(C)$
(cf Corollary II.2.8 of \cite{BH}). Thus $C_G(C)$
leaves invariant the $\ga$--axis $A_c =\{c\}\times \R$. 
Note that $A_c$ has a metric simplicial structure (induced from the structure
on $\D$) with discrete automorphism group $\Aut(A_c)$. The group $C_G(C)$ acts 
via a homomorphism to $\Aut(A_c)$ whose kernel we denote $H$. 
Moreover the translation $\ga$ acts co-compactly on the axis, so
generates a finite index subgroup of $\Aut(A_c)$. 
It follows that $C_G(C)$ contains the product $H\times C$ with finite index. 
Note also that the only points in $\D$
which have non-abelian stabilizer are the rank 2 vertices. Since these form
a discrete set, while the fixed set of $H$ contains a whole real line $A_c$, it
follows that  $H$ must be abelian (either trivial or infinite cyclic). 
But then $C_G(C)$ is virtually abelian, a contradiction.
\end{proof}


\section{Abstract commensurators and the graph $\Theta$ of fixed sets in $\D$}

We recall briefly the definition of an abstract commensurator of groups.
Given groups $\Ga_1,\Ga_2$, we define 
\[
\Comm(\Ga_1,\Ga_2) = \{\,\varphi\co H_1\buildrel\cong\over\rightarrow H_2\ :\ 
H_i<\Ga_i \text{ finite index, } i=1,2\,\}\,/\lower.75ex\hbox{$\sim$}\,,
\]
where isomorphisms $\varphi$ and $\psi$ are equivalent, $\varphi\sim\psi$, 
if they agree on restriction to a common finite 
index subgroup of their domains. Elements of $\Comm(\Ga_1,\Ga_2)$ shall be called
\emph{abstract commensurators} from $\Ga_1$ to $\Ga_2$, and when this set is nonempty
we shall say that $\Ga_1$ and $\Ga_2$ are \emph{abstractly commensurable}.
Note that when $\Ga_1$ and $\Ga_2$ are the same group this set has the
structure of a group (under composition of isomorphisms after passing to
appropriate finite index subgroups). We shall write $\Comm(\Ga)=\Comm(\Ga,\Ga)$
and refer to this as the \emph{abstract commensurator group} of $\Ga$.  
Note that there is a natural homomorphism $\Aut(\Ga)\to\Comm(\Ga)$ whose kernel
consists of those automorphisms which fix a finite index subgroup $\Ga$. 
 
Before continuing, we make some general observations concerning the relationships
between a 2--dimensional Artin group, its automorphism group and its abstract
commensurator group. If $\Delta$ is a 2--dimensional defining graph with at least 3 
vertices then $G=G(\Delta)$ has a trivial centre and so is isomorphic to $\Inn(G)$.
Moreover, consideration of Lemma \ref{fixedsets}(i) shows that $s$ is the unique 
$N$th root of $s^N$ for any generator $s\in V(\Delta)$ and any $N\in \N$. It follows
that any automorphism of $G$ which restricts to the identity on a finite index 
subgroup of $G$ is the identity on all of $G$. 
Thus the natural homomorphism $\Aut(G)\to\Comm(G)$
is injective.  Identifying $\Aut(G)$ with its image, we have
\[
G\cong \Inn(G)< \Aut(G)< \Comm(G)\,.
\]
We now turn to the class of CLTTF Artin groups. 
Our principal tool for studying abstract commensurators
between these groups is the following structure:

\begin{defn}[{\rm(}The fixed set graph $\Theta$\/{\rm)}]
Let $\Delta$ denote a CLTTF defining graph. We define the following sets
of subsets of the Deligne complex $\D$ of type $\Delta$: 
\[
\begin{aligned}
\V &=\{ \text{ singletons }\{gV_e\} : g\in G\,,\ e\in E(\Delta)\}\,,\text{ and }\\
\F &= \{ \text{ unbounded trees } gF_s :  g\in G\,,\ s\in V(\Delta)\}\,.
\end{aligned}
\]
We define the \emph{fixed set graph} $\Theta=\Theta(\Delta)$ 
to be the bipartite graph with the following vertex and edge sets: 
\[
\begin{aligned}
\text{Vert}(\Theta)&:=\V\cup \F\\
\text{Edge}(\Theta)&:=\{\,(V,F)\,:\, V\in\V,\  F\in\F\, \text{ and }\, V\subset F\,\}\,.
\end{aligned}
\]
\end{defn}

Observe that, by Lemma \ref{fixedsets}, Remark \ref{rem:unboundedtrees}, 
and Lemma \ref{starlike}, and since we are supposing large type (LT), we have that
\[
\V\cup \F = \{\, \Fix(C)\, :\, C \text{ is a  CNVA  subgroup of } G\,\}\,,
\] 
where $\Fix(C)\in\V$ if $C$ is conjugate to a subgroup of
$\< z_e\>$ for some $e\in E(\Delta)$, and $\Fix(C)\in\F$ if 
$C$ is conjugate to a subgroup of $\< s\>$ for some internal vertex 
$s\in V(\Delta)$.

\begin{lemma}\label{commute}
Let $C,C'$ be  CNVA  subgroups of $G$. Then 
\begin{itemize}
\item[\rm(i)] $C\cap C'\neq\{ 1\}$ if and only if $\Fix(C)=\Fix(C')$.
\item[\rm(ii)] $\< C,C'\>\cong\Z\times\Z$ if and only if 
$(\,\Fix(C),\Fix(C')\,)\in\text{Edge}(\Theta)$.
\end{itemize} 
\end{lemma}

\begin{proof} 
\noindent\textbf{(i)}\qua If $\Fix(C)=gF_s$ for some $g\in G$ and $s\in V(\Delta)$, then 
$C<g\< s\>g^{-1}$ since $gV_s\in gF_s$. On the other hand, if $\Fix(C)=\{gV_e\}$,
for some $g\in G$ and $e\in E(\Delta)$, then $C_G(C)<gG(e)g^{-1}$ and, by 
Lemma \ref{structureGe}(i), $C<g\<z_e\>g^{-1}$ (else it fails to be  CNVA).
Therefore, if $\Fix(C)=\Fix(C')$ then $C$ and $C'$ lie in a common infinite cyclic
subgroup, so must intersect nontrivially. 

On the other hand, it follows from Lemma \ref{fixedsets} that
a cyclic subgroup of $G$ has the same fixed set as any of 
its nontrivial subgroups. Thus, if  $C''=C\cap C'$ is nontrivial we have
$\Fix(C)=\Fix(C'')=\Fix(C')$.
\medskip

\noindent\textbf{(ii)}\qua
Suppose that $\Fix(C)=V\in\V$ and $\Fix(C')=F\in\F$ such that $V\subset F$.
Up to conjugation of $C,C'$ in $G$ we may suppose that $F$ contains the 
edge $[V_e,V_s]$, for some $e=\{s,t\}\in E(\Delta)$, and that $V=\{ V_e\}$.
(This is because any edge of $F$ emanating from $V$ can be viewed as the translate
of some edge in the fundamental region $K$).
But then  we have $C<\< z_e\>$ (since $\Fix(C)=V_e$) and $C'<\Stab(V_s)=\< s\>$. 
Since $\< s, z_e\>\cong\Z\times\Z$ it follows that $\< C,C'\>\cong\Z\times\Z$.
 
Suppose now that $\< C,C'\>\cong\Z\times\Z$.
It follows, since they commute, that $C$ and $C'$ have a common 
fixed point in $\D$, for $C$ must leave $\Fix(C')$ invariant and 
so fixes the orthogonal projection $p'\in\Fix(C')$ of any point 
$p\in\Fix(C)$. However, a rank 2 abelian subgroup can only fix
a rank 2 vertex. So $\Fix(C)\cap\Fix(C')$ consists of a single vertex 
$V\in\V$, say. Up to conjugation by an element of $G$ we may suppose that
$C,C'<\Stab(V_e)=G(e)$, for some $e=\{s,t\}\in E(\Delta)$. 
Each of the two  CNVA  subgroups is then either 
a subgroup of $Z(G(e))=\< z_e\>$ or conjugate in $G(e)$ to a subgroup of
$\< s\>$ or of $\< t\>$. By Lemma \ref{structureGe}(iii) they cannot both be
of the latter kind unless they lie in a common cyclic subgroup. Similarly,
they cannot both lie in the centre. But then one is central and one is
conjugate into $\< s\>$ say. That is to say that, up to conjugacy in $G$, we
have $\{\Fix(C),\Fix(C')\}=\{V_e,F_s\}$. 
\end{proof}

\begin{prop}\label{ThetaIsom}
Let $\Delta,\Delta'$ denote CLTTF defining graphs, and suppose that 
$\varphi \co H\to H'$ is an abstract commensurator from $G(\Delta)$ to $G(\Delta')$.
Then $\varphi$ determines a unique
well-defined graph isomorphism $\Phi \co \Theta(\Delta)\to\Theta(\Delta')$ 
such that $\Phi(\Fix(C)) = \Fix(\varphi(C\cap H))$ for any  CNVA  
subgroup $C$ of $G(\Delta)$.
\end{prop}  

\begin{proof}
This is a consequence of Lemma \ref{commute} above and
the fact that the properties ``$C$ is  CNVA'', ``$C\cap C'\neq\{ 1\}$'' and 
``$\< C,C'\>\cong\Z\times\Z$'' are all preserved by isomorphism and 
passage to finite index subgroups.
\end{proof}

\begin{remark} 
Consider a fixed $\Delta$ of type CLTTF, and write $G=G(\Delta)$ 
and $\Theta=\Theta(\Delta)$.
Note that the action of $G$ on $\D$ induces an action of $G$ by graph
automorphisms of $\Theta$. We remark that the action of
$\Comm(G)$ on $\Theta$ given by the above Proposition extends this action
of $G$ when $G$ is identified with the subgroup of $\Comm(G)$
consisting of inner automorphisms. 
\end{remark}


\section{Circuits in the graph $\Theta$}\label{sect:Circuits}

In this and subsequent sections we analyse the structure of the graph 
of fixed sets  associated to a CLTTF defining graph $\Delta$. For simplicity
we shall write $\Theta=\Theta(\Delta)$ and $G=G(\Delta)$.

\begin{lemma}
If $\Delta_2$ denotes the first subdivision of $\Delta$
(so $\Delta_2$ has vertex set $E(\Delta)\cup V(\Delta)$)
we let $\what\Delta$ denote the full subgraph of $\Delta_2$ 
spanned by the non-terminal vertices. Then there is a graph embedding
$f\co \what\Delta\hookrightarrow \Theta$ defined by $f(s)=F_s$, if $s\in V(\Delta)$ 
is an internal vertex, and $f(e)=V_e$, if $e\in E(\Delta)$. 
\end{lemma}

\begin{proof}
It is clear that, as written, $f$ is a well-defined graph morphism. Clearly, also, 
$f$ is injective on $E(\Delta)$. Suppose $s,t\in V(\Delta)$ and $s\neq t$. 
Suppose that $F_s=F_t$. Then, by convexity, $F_s$ must contain the
geodesic segment $[V_s, V_t]$. However, $[V_s, V_t]$ intersects the interior of
the fundamental region of $\D$, while $F_s$ does not, a contradiction.
Therefore $f$ is injective on $V(\Delta)$.
\end{proof} 

From now on we shall identify $\what\Delta$ with its image $f(\what\Delta)$ in $\Theta$.
We also observe that $\Theta$ is the union of translates of the subgraph
$\what\Delta$ by elements of $G$. In particular, since we suppose that $\Delta$
connected (C), we deduce that $\Theta$ is also connected ($G$ is generated by
elements which individually fix some part of $\what\Delta\,$).
A particular consequence of this is that any automorphism
of $\Theta$ respects the given bi-partite structure. 
However, we have not ruled out the possibility that
$\Phi(\V)=\F$ and $\Phi(\F)=\V$ for some $\Phi\in\Aut(\Theta)$.

In order to understand which structural properties of the graph $\Theta$ are 
respected by graph isomorphisms (coming from abstract commensurators of $G$)
we shall study the properties of simple closed circuits in $\Theta$.

Let $\Sigma = (V_1, F_1, V_2,F_2,\ldots , V_k, F_k)$ denote a simple circuit
of length $2k$ in $\Theta$, where $V_i\in\V$, $F_i\in\F$, and
$V_i,V_{i+1}\subset F_i$ for each $i=1,2,..,k$, with indices taken mod $k$
(so that $V_1\subset F_k$). Note that $F_{i-1}\cap F_i=V_i$ and consists of a
single rank 2 vertex of $\D$.  For each $i=1,..,k$, let $\ga_i$ denote the
geodesic segment in $F_i$ from $V_i$ to $V_{i+1}$. To the simple 
circuit $\Sigma$ we associate the  closed polygonal curve 
$\ov\Sigma = (\ga_1,\ldots ,\ga_k)$ in $\D$.
Note that each segment $\ga_i$ is an edge path in the 1-skeleton of $\D$
and is geodesic in $\D$ (by convexity of $F_i$).

\begin{defn}[{\rm(}Basic circuit\/{\rm)}]
A simple circuit $\Sigma$ in $\Theta$, and its associated polygon
$\ov\Sig$ in $\D$, are said to be \emph{basic} if $\Sigma$ is the translate
by an element of $G$ of a simple circuit in the subgraph $\what\Delta$, 
equivalently if the polygon $\ov\Sigma$ lies wholly in (the boundary of) 
a single translate of the fundamental region $K$ in $\D$.    
\end{defn}

Note that the property of being a basic circuit depends upon the structure
of the Deligne complex (rather than just the structure of $\Theta$). We wish to
characterize certain basic circuits purely in terms of the graph theoretic
properties of $\Theta$.

\begin{defn}[{\rm(}Minimal circuit\/{\rm)}]
Let $\Sigma$ denote a simple circuit in $\Theta$. A \emph{short-circuit} of $\Sigma$
is any simple path $A$ in $\Theta$ which intersects $\Sig$ only in its endpoints
$P,Q$, and which is strictly shorter than any path in $\Sigma$ from $P$ to $Q$.
We say that $\Sig$ is a \emph{minimal} circuit if it is a simple circuit and admits
no short-circuit. (More succinctly, a circuit is minimal if and only
if it is isometrically embedded  
when $\Theta$ is viewed as a metric graph with edges of constant length). Note that 
if the circuit $\Sigma$ admits a short circuit $A$ then we may decompose
$\Sig$ into a pair of  simple circuits each of length strictly smaller than $\Sig$,
namely:
\[
\Sig_1=A_1A \hskip3mm\text{ and }\hskip3mm\Sig_2=A^{-1}A_2
\,,\hskip3mm\text{ where }\hskip3mm\Sig=A_1A_2\,.
\]
This provides an inductive procedure for reducing an arbitrary simple circuit
into (a finite collection of) minimal circuits.
\end{defn}

We devote the next section to proving the following two Propositions.

\begin{prop}\label{first}
Any minimal circuit of $\Theta$ is a basic circuit.
\end{prop}

\begin{prop}\label{second}
Any  minimal circuit of $\what\Delta$ is minimal as a circuit of $\Theta$.
\end{prop}

\begin{remark}
While the family of all minimal circuits of the graph $\Theta$ is easily seen to be
preserved by any abstract commensurator of $G$, the above 
Propositions show that this structure in the graph $\Theta$
is also closely related to the combinatorial definition of $G$, and hence to the
structure of the Deligne complex $\D$. Namely, the minimal circuits are precisely
the translates in $\Theta$ of the minimal circuits of $\what\Delta$. 
This connection to the Deligne complex shall be developed further 
in subsequent sections, and will ultimately lead to the proof of Theorem~\ref{Thm3}.
\end{remark}


\section{On minimal circuits --  Propositions 
\ref{first} and \ref{second}}\label{sect:TwoProps}

Let $\Delta$ be a CLTTF defining graph. 
Throughout this section we shall regard the associated Deligne complex $\D$
as a squared complex equipped with the cubical metric $d_C$.
Since $\Delta$ is of type CLTTF, the metric space $(\D,d_C)$ is 
a CAT(0) squared complex.
 
We begin with a useful lemma which reflects the $\delta$--hyperbolicity 
of the Deligne complex.

\eject

\begin{lemma}\label{supergeodesic}
Let $F\in\F$, and let $\ga$ be any geodesic segment in $F$ which passes through a 
rank 2 vertex $p$. Then $\ga$ is ``super-geodesic" at $p$, by which we mean
that $\ga$ enters and leaves 
$p$ through points separated in $\Lk(p,\D)$ by a path distance strictly greater 
than $\pi$, in fact at least $3\pi/2$.  
\end{lemma}

\begin{proof}
Let $E,E'$ denote the edges of $F$ along which $\ga$ enters and leaves the point $p$.
If $E$ and $E'$ define points in $\Lk(p,\D)$ which are joined by a path of length $\pi$
then there exist squares $Q,Q'$ adjacent to $E,E'$ respectively, which share a common 
edge $E''$ such that $E''\cap F=\{p\}$. Let $g$ denote a nontrivial element of $\Stab(F)$,
and observe that $g(E'')\neq E''$. It follows that the squares $Q,Q',g(Q)$ and $g(Q')$
form a larger square with the vertex $p$ at its centre. In particular 
we see that $\Lk(p,\D)$ contains a simple circuit of length exactly $2\pi$. 
Since $G(\Delta)$ is assumed to be large type (LT), the shortest simple circuit 
in the link of a rank 2 vertex of $\D$ has length  at least $3\pi$ 
(see \cite{CD}, also Lemma \ref{ShortCycles} of Section \ref{sect:LinkTheta} below),
a contradiction. Thus, any path in $\Lk(p,\D)$ 
from $E$ to $E'$ has length strictly greater than $\pi$, and therefore at least 
$3\pi/2$ since all edges of the link graph are of length $\pi/2$.   
\end{proof}
 
In the following arguments we shall use the properties of walls (or hyperplanes)
in a CAT(0) cubed complex. The notion appears frequently in the literature. See
for example \cite{Sageev}, or \cite{NR}. Two edges in a CAT(0) squared 
complex $X$ may be said to be parallel if they are opposite edges of the 
same square in the complex (more generally, if they are parallel edges of the 
same $n$--cube in the case of a higher dimensional cube complex).
This generates an equivalence relation on the set of all edges. 
By a \emph{wall} in $X$ we shall mean the convex subspace spanned by 
the midpoints of the edges lying in a single parallelism class.
Since $X$ is CAT(0) and 2--dimensional
 this defines a tree which is isometrically embedded in $X$.
Moreover, a wall in $X$ separates $X$ into exactly two components, usually called 
\emph{half-spaces}. 

\begin{defn}[{\rm(}$\mF\subset\D$; $W^+$, $W^-$ and $\partial W^+$ for a wall $W$\/{\rm)}]
It will be convenient to write $\mF$ for the subcompex of $\D$ which is the 
union of the sets $F\in\F$. This is the largest subcomplex of $\D$
which contains no rank 0 vertices, equivalently the set of all points in $\D$
with nontrivial stabilizer in $G(\Delta)$. 
We note that any wall $W$ in $\D$ may be naturally 
oriented: we thus denote the connected components of $\D\setminus W$ by 
$W^+,W^-$ in such a way that every edge of $\mF$ which crosses $W$ has a 
rank 1 vertex in $W^-$ and a rank 2 vertex in $W^+$. (All other edges of $\D$ 
which cross $W$ have a rank 0 vertex in $W^-$ and a rank 1 vertex in $W^+$).
Let $\partial W^+$ denote the subcomplex of $\D$ spanned by those vertices in $W^+$
which belong to edges crossing $W$.
Thus $\partial W^+$ is a parallel copy of $W$ spanned by vertices of rank 1 and 2.
In particular, $\partial W^+$ is a convex subtree of $\D$, and also a subcomplex
of $\mF$.
\end{defn}

\begin{defn}[{\rm(}Orthogonality\/{\rm)}]
We shall say that two convex subsets of $\D$ \emph{intersect orthogonally} at a
point $p$ if the orthogonal projection of either one to the other contains only the point $p$.
Note that convex subsets of the 1-skeleton of the squared complex $\D$ 
which intersect in a single point always intersect orthogonally, since all angles 
are multiples of $\frac{\pi}{2}$.
\end{defn}

\begin{defn}[{\rm(}$\V$--paths\/{\rm)}]
By a \emph{$\V$--path} we shall mean an edge path
in $\Theta$ whose initial and terminal vertices lie in $\V$. 
Given a $\V$--path $A$ we shall write $\ov A$ to denote the 
piecewise geodesic path in $\D$ induced by $A$ (in the manner
described previously for simple circuits).
Thus, if $A=(V_1,F_1,V_2,F_2,\ldots,V_k,F_k,V_{k+1})$ denotes a $\V$--path of
length $2k$ where $V_i\in\V$, $F_i\in\F$ and $V_i,V_{i+1}\subset F_i$ for $i=1,..,k$,
then $\ov A$ is simply the union of the geodesic segments joining
$V_i$ to $V_{i+1}$ in $F_i$, for $i=1,..,k$.

Given a $\V$--path $A$ in $\Theta$ we shall write $\W(A)$ to denote
the set of walls of the squared complex $\D$ which are traversed
by the induced path $\ov A$. In particular, if $\ov A$ happens 
to be geodesic in $\D$ then $\W(A)$ is exactly the set of walls 
which separate the endpoints of $\ov A$. 

We shall use $L(A)$ to denote the edge length of a path $A$ in $\Theta$. 
For $\V$--paths this length is always even.
\end{defn}

\begin{lemma}\label{Vpaths}
Let $A,B$ denote $\V$--paths in $\Theta$, and write $\al=\ov A$ and $\be=\ov B$.
Suppose that $\al$ is a nontrivial geodesic in $\D$ and $\W(A)\subset\W(B)$. Then
\begin{enumerate}
\item[\rm(i)] $L(A)\leq L(B)$, and
\item[\rm(ii)] if, moreover, $\al$ and $\be$ share a common endpoint, $p$, but
do not both leave the vertex $p$ along the same edge of $\D$, 
then the inequality is strict: $L(A)< L(B)$.
\end{enumerate}
\end{lemma}

\begin{proof} {\bf (i)}\qua
We shall compare the number of vertices of type $\V$ appearing along each path.
Note that, since $\al$ is geodesic and intersects orthogonally with each element of
$\W(A)$, the walls $\W(A)$ are mutually disjoint.
Let $V$ denote a type $\V$ vertex lying along the path $A$ and let 
$W_1,W_2$ denote the  walls of $\D$ traversed by 
$\al$ immediately before and after passing through $V$. 
Since $V$ is a rank 2 vertex, we have $W_1\subset W_2^+$ and $W_2\subset W_1^+$. 
Since it traverses both walls, the path $\be$ must pass across the
region $W_1^+\cap W_2^+$ between the two walls.
We now claim that $B$ also has a vertex of type $\V$ 
 lying in the region $W_1^+\cap W_2^+$ containing $V$.
Since these regions for the different type $\V$ vertices of $A$ are disjoint,
it follows that $B$ has at least as many type $\V$
vertices as $A$ does, and hence that $L(A)\leq L(B)$.

To see the claim, observe that, if the path $B$ has no type $\V$ vertex falling between
$W_1$ and $W_2$ then it must contain a subpath $(V,F,V')$, with $V,V'\in\V$
and $F\in\F$, such that the geodesic segment of $\be$ joining $V$ to $V'$ 
in $F$ crosses both $W_1$ and $W_2$. 
Let $a,b,c,d$ denote the four points of intersection between 
convex sets $\al$, $W_1$, $W_2$, and the geodesic segment of $\be$ just mentioned.
 Since all intersections are orthogonal, 
the Flat Quadrilateral Theorem, II.2.11 of \cite{BH}, implies that the four points $a,b,c,d$
and the geodesic segments connecting them in $\al$, $\be$, $W_1$ and $W_2$
form the boundary of a flat Euclidean rectangle isometrically embedded in $\D$.
We conclude, by Lemma \ref{supergeodesic}, that $\be$ does not pass through any
rank 2 vertex between $W_1$ and $W_2$, 
since $\be$ would have to be super-geodesic at any such vertex.
  This is of course a contradiction, since any
edge of $\mF$ which crosses $W_1$ in the direction of $W_2$ immediately encounters
a rank 2 vertex (because $W_2\subset W_1^+$).
\medskip

\noindent{\bf (ii)}\qua
Suppose without loss of generality that $p$ is the initial vertex of both 
$\al$ and $\be$. Note that, since $\al$ is nontrivial, $0\neq L(A)\leq L(B)$ and 
we may write $B= (V,F,B')$ where $V=\{p\}$, $F\in\F$,
 and $B'$ is a $\V$--path of strictly smaller length: $L(B')=L(B)-2$. 
Since $\al$ and $\be$ set off along different edges
of the Deligne complex, the convex sets $\be\cap F$ and $\al$ must 
intersect orthogonally at $p$. It follows that no wall of $\W(A)$ 
can cross $\be\cap F$, so $\W(A)\subset\W(B')$.
Applying part (i) of the Lemma, we conclude that $L(A)\leq L(B')<L(B)$. 
\end{proof}

\begin{defn}[{\rm(}Chords\/{\rm)}]
Let $\Sigma$ denote a simple circuit in $\Theta$, and $\ov\Sigma$ the
corresponding piecewise geodesic closed curve in $\D$. 
By a \emph{chord} of $\ov\Sigma$ we mean a path $\alpha$ from $p$ to $q$ 
such that 
\begin{description}
\item[(C1)]{$\alpha$ is a geodesic path in $\D$ and is contained in the subcomplex $\mF$,}
\item[(C2)]{$\alpha\cap\ov\Sigma =\{ p,q\}$ where the endpoints $p$ and $q$ are vertices of $\D$.}
\end{description}
Observe that the endpoints $p,q$ of a chord $\al$ serve to cut $\ov\Sigma$ into 
a concatenation of two paths $\sig_1$ and $\sig_2$ each joining $p$ to $q$.
We shall say that $\al$ is \emph{aligned with $\sig_1$} if $\al$ and $\sig_1$ 
form a right angle at each endpoint.

A rank 2 vertex lying
on the path $\ov\Sig$ shall be termed \emph{essential} if it corresponds to
a type $\V$ vertex of $\Sigma$, and \emph{inessential} otherwise.
\end{defn}

\begin{lemma}\label{sigalign}
Let $\Sigma$ denote a minimal circuit of $\Theta$. Then
\begin{itemize}
\item [\rm(i)] $\ov\Sig$ is a simple closed curve 
in $\D$, and 
\item [\rm(ii)] if $\al$ is a chord of $\ov\Sigma$, cutting $\ov\Sigma$ 
into subpaths $\sig_1,\sig_2$, then both endpoints of $\al$ are
inessential rank 2 vertices and $\al$ is aligned with either $\sig_1$ or $\sig_2$
(but not with both).
\end{itemize}
\end{lemma}

\begin{proof} {\bf (i)}\qua
The essential vertices of $\ov\Sig$ 
are clearly mutually distinct, since $\Sigma$ is a simple circuit of $\Theta$.
Suppose then that $p$ is a point which lies on the interior of more than one
``side'' of $\ov\Sigma$, and let $P,Q$ denote fixed trees which appear as
distinct vertices of $\Sigma$ (of type $\F$) and which both contain $p$.
Then $P\cap Q=V$, where $V=\{p\}\in\V$.
In particular, $V$ is the unique element of
$\V$ which is adjacent, in $\Theta$, to both $P$ and $Q$. 
It follows that $(P,V,Q)$ defines a short circuit for $\Sig$,
contradicting the minimality of $\Sig$.
\medskip

\noindent{\bf (ii)}\qua
Let $P$ (resp. $Q$) denote the smallest fixed set which appears as a 
vertex of $\Sigma$ and which contains the endpoint $p$ (resp. $q$) of $\al$.
Either $P=\{ p\}\in\V$, or $P\in \F$ is an unbounded tree containing $p$ (similarly for $Q$).
For any $\V$--path $X$ in $\Theta$ with endpoints contained in $P$ and $Q$ respectively
we shall denote by $PXQ$ the path in $\Theta$ obtained by appending the vertices
$P$ and/or $Q$ whenever they belong to $\F$ (ie, whenever they
 are not already endpoints of $X$).   

Let $\sig'_1,\sig'_2$ and $\al'$ denote the
shortest subpaths of $\sig_1$, $\sig_2$, and $\al$ respectively, which have an endpoint
in each of the sets $P$ and $Q$. 
Note that the endpoints of these subpaths are all necessarily  rank 2 vertices (since
the paths all lie in $\mF$, which is the union of the elements of $\F$, and any two elements
of $\F$ meet, if at all, in a vertex of rank 2). 
Now $\Sigma$ may be expressed as the union of paths
$P\Sig_1Q$ and $P\Sig_2Q$ where $\Sig_1$ and $\Sig_2$ are $\V$--paths with
$\ov\Sig_i=\sig'_i$, for $i=1,2$.
One may also easily construct a $\V$--path $A$, with endpoints in $P$ and $Q$,
such that $\ov A=\al'$.
We observe that $\al'$ intersects the sets
$P$ and $Q$ orthogonally at its endpoints. This ensures that
 $\W(A)\subset\W(\Sig_i)$ for each $i=1,2$.

If either $P\in\V$ or $Q\in\V$ then it follows from Lemma \ref{Vpaths}(ii) that
$L(A)<L(\Sig_i)$ for $i=1,2$. But then $PAQ$ is a short circuit for $\Sigma$, 
contradicting minimality.
 
Now suppose that both $P,Q\in\F$. In particular, each of the endpoints $p,q$ of $\al$
is either rank 1 or an inessential rank 2 vertex of $\ov\Sig$.
Note that if $p$ is a rank 1 vertex then $\sig_1\cup\al$
is geodesic at the point $p$. This is also the case if $p$ is rank 2, as long as
$\al$ and $\sig_1$ do not form a right angle at $p$.
In either case $\al'$ extends (through $P$) to a geodesic having an
endpoint, $v$ say, in common with $\sig_1'$.
Applying Lemma \ref{Vpaths} to the $\V$--paths $(\{v\},P,A)$ 
and $\Sig_1$ now shows that $L(A)<L(\Sig_1)$,
except possibly when $p$ is inessential, rank 2 and 
$\al$ forms a right angle at $p$ with $\sig_1$.
By applying this argument at each end of $\al$ and with respect to both $\Sig_1$ 
and $\Sig_2$, we conclude that $PAQ$ defines a short-circuit for $\Sig$ 
unless both $p$ and $q$ are inessential rank 2 vertices of $\ov\Sig$  
and $\al$ is aligned with one of $\sig_1$ or $\sig_2$.
Finally, note that $\al$ cannot be aligned with \emph{both} $\sig_1$ and $\sig_2$ 
since $\ov\Sig$ is super-geodesic at $p$ and $q$, by Lemma \ref{supergeodesic}.
\end{proof}

\begin{lemma}\label{sigma-basic}
Let $\Sigma$ be a minimal circuit of $\Theta$, and suppose that all rank 2 vertices
of $\ov\Sig$ are essential. Then $\Sigma$ is a basic circuit.
\end{lemma}

\begin{proof}
The fact that all rank 2 vertices  are essential implies, 
by Lemma \ref{sigalign}(ii), that $\ov\Sig$ has no chords. It follows that, for every
wall $W$ which intersects $\ov\Sig$, the set $\partial W^+\cap \ov\Sig$ is connected.
For, otherwise, any shortest length path in $\partial W^+$ joining distinct components
of $\partial W^+\cap \ov\Sig$ would be a chord of $\ov\Sig$.

Given a wall $W$, we shall define 
$N(W^-)$ to be the convex hull in $\D$ of the subsets $W^-$ and $\partial W^+$.
It follows from the connectivity statement above and the fact that $\ov\Sig$ is
a simple closed curve (Lemma \ref{sigalign}(i)) that, for any wall $W$ of $\D$, 
either  $\ov\Sig\subset W^+$ or $\ov\Sig\subset N(W^-)$.

Consider the fundamental region $K$ of $\D$, with rank 0 vertex $V_\emptyset$.
For each $t\in V(\Delta)$ we define $W_t$ to be the unique wall of $\D$ which cuts
the edge $[V_\emptyset,V_t]$. These are exactly the walls of $\D$ which intersect 
nontrivially with $K$. Moreover, we have
\[
K=\bigcap_{t\in V(\Delta)} N(W_t^-)\,.
\]
Let $W_0$ be any wall crossed by $\ov\Sig$. Note that, since $\ov\Sig$ is a simple 
closed curve  it must cross $W_0$ in at least two different places.
 The intersection $\ov\Sig\cap\partial W_0^+$ must therefore contain at least two rank 2 
vertices and, by connectedness, at least one rank 1 vertex.  Let $v$ denote such a
rank 1 vertex in $\ov\Sig\cap\partial W_0^+$, and let $u$ be the unique rank 0 vertex such that
$u$ is adjacent to $v$ and the edge $[u,v]$ crosses $W_0$. Up to an isometry of $\D$
we may as well suppose that $[u,v]=[V_\emptyset,V_s]\subset K$ for some $s\in V(\Delta)$,
and hence that $W_0=W_s$.
Note that $W_s$ is the unique wall which separates $V_s$ from $V_\emptyset$.
It follows that, for any $t\in V(\Delta)$ different from $s$, we have $V_s\subset W_t^-$.
Thus  $\ov\Sig\subset N(W_t^-)$ (as $\ov\Sig$ is clearly not contained in $W_t^+$).
Since $\ov\Sig$ crosses the wall $W_s$, we must also have $\ov\Sig\subset N(W_s^-)$.
It follows that $\ov\Sig\subset K$, and so $\Sig$ is a basic circuit.  
\end{proof}

\subsection{Proof of Proposition \ref{first}} 

Let $\Sig$ denote a minimal circuit in $\Theta$. In order to show that
$\Sig$ is basic it will suffice, by Lemma \ref{sigma-basic},
to show that every rank 2 vertex of $\ov\Sigma$ is essential 
(ie, corresponds to a vertex of $\Sig$ of type $\V$). 

\begin{figure}[ht!]\small\anchor{ChordsFig}
\psfrag {a1}{$\alpha_1$}
\psfrag {a2}{$\alpha_2$}
\psfrag {s1}{$\sigma_1$}
\psfrag {s2}{$\sigma_2$}
\psfrag {s3}{$\sigma_3$}
\psfrag {q1}{$q_1$}
\psfrag {q2}{$q_2$}
\psfrag {p}{$p$}
\cl{\includegraphics[width=7cm]{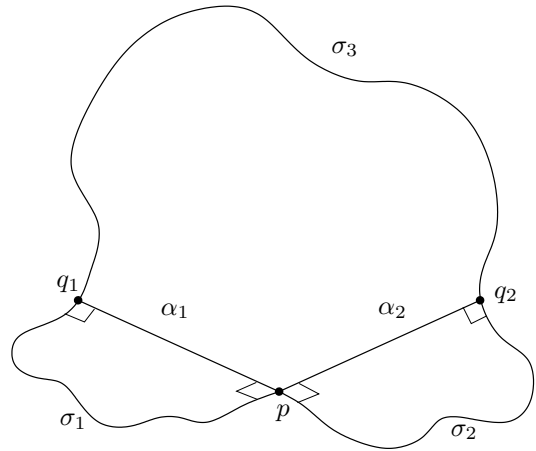}}
\caption{Chords $\al_i$ aligned with $\sig_i$ for $i=1,2$}
\label{ChordsFig}
\end{figure}

Suppose that $p$ is an inessential rank 2 vertex on $\ov\Sig$.
Let $E_1,E_2$ denote the two edges of $\ov\Sigma$ adjacent to $p$,
and $W_1,W_2$ the corresponding walls ($W_i$ is the unique wall of $\D$ which cuts across
the edge $E_i$). Since $\ov\Sig$ is geodesic at $p$ (in fact super-geodesic) 
 and also a simple closed curve (by Lemma \ref{sigalign}(i)),
the point $p$ is an isolated point of the set $\ov\Sig\cap\partial W_1^+$. 
Also, since $W_1$ is separating and $\ov\Sig$ is a simple closed curve,
the set $\ov\Sig\cap\partial W_1^+$ must contain points
other than $p$. Any shortest length path in $\partial W_1^+$ 
joining $p$ to another component of $\partial W_1^+\cap \ov\Sig$ now defines
a chord of $\ov\Sig$, which we shall denote $\al_1$. Note that $\al_1$ forms a right
angle with $\sig_1$ at $p$.
Similarly, there exists a second chord $\al_2$ which lies in $\partial W_2^+$ and 
and forms a right angle at $p$ with $\sig_2$. Let $q_1$, $q_2$ denote the endpoints of 
$\al_1$, $\al_2$ respectively, which are different from $p$, and write
$\ov\Sig$ as a union of paths $\sig_1$ from $p$ to $q_1$, 
$\sig_2$ from $p$ to $q_2$, and $\sig_3$ from $q_1$ to $q_2$.
By Lemma \ref{sigalign}(ii), $q_1$ and $q_2$ are both 
inessential rank 2 vertices, and $\al_i$ is necessarily aligned
with $\sig_i$, for each $i=1,2$.
The situation is illustrated in \figref{ChordsFig}.
We shall show that the union of the chords $\al_1$ and $\al_2$ 
represents a short-circuit of $\Sig$, and thereby obtain a contradiction.

Write $P$, $Q_1$ and $Q_2$ for the vertices of type $\F$ in $\Sigma$
which correspond to fixed sets containing the points $p,q_1$ and $q_2$ respectively.
By Lemma \ref{supergeodesic}, $\ov\Sig$ is super-geodesic at $p$, $q_1$ and $q_2$.
It follows that the path $\al_1\cup\sigma_3\cup\al_2$ is geodesic at both 
points $q_1$ and $q_2$, and also that $\al_1\cap\al_2=\{ p\}$, or rather that $\al_1$ and
$\al_2$ intersect orthogonally at $p$. 

For $i=1,2$, let $A_i$ denote the $\V$--path such that $\ov A_i=\al_i$. 
We may suppose that, for each $i=1,2$, the path $A_i$ is written
\[
A_i= (V_1^{(i)}, F_1^{(i)},V_2^{(i)}, F_2^{(i)},\ldots , V_{m_i}^{(i)})\,, \ m_i\in\N
\]
where $V_1^{(1)}=V_1^{(2)}=\{ p\}\subset P$ and
$V_{m_i}^{(i)}=\{q_i\}\subset Q_i$, for $i=1,2$.  

Let $\Sig_3$ denote the longest $\V$--subpath of $\Sigma$ such 
that $\ov\Sig_3$ is a subpath of $\sig_3$, and write $U_1\subset Q_1$ 
and $U_2\subset Q_2$ for the endpoints of $\Sig_3$.

We now define $\V$--paths $A'_i$ for $i=1,2$ thus:
\[
A'_i= (V_2^{(i)}, F_2^{(i)},\ldots , V_{m_i}^{(i)},Q_i,U_i)\,.
\]
By the preceding remark on geodicity at the points $q_i$, the path
$\ov A'_i$ is a geodesic in $\D$ and $L(A'_i)=L(A_i)$, for each $i=1,2$.

We now make the following claim:

\begin{claim}
The first walls crossed by $\ov A'_1$ and $\ov A'_2$ are
disjoint and separate  $U_1$ from $U_2$.
\end{claim}

\begin{proof}
To see this, let $\ga_i$, for $i=1,2$, denote the geodesic from $\{p\}$ to $U_i$ 
obtained by combining $\ov{A_i}$ and $\ov{A_i'}$.
Either the first walls crossed by $\ga_1$ and $\ga_2$ have the desired
properties, and hence also the first walls of $\ov{A_1'}$ and $\ov{A_2'}$,
or the first walls crossed by $\ga_1$ and $\ga_2$ intersect in a square
of $\D$ with rank 0 vertex $v_0$ and rank 2 vertex $p$. The vertex $v_0$
is the centre of a region $K_0=g_0(K)$ for some $g_0\in G$, and the 
geodesics $\ga_1,\ga_2$ each intersect $K_0$ in a subpath of edge length $2$, 
and leave $K_0$ at the rank 2 vertex $p_1$, respectively $p_2$.
Let $W_i$ denote the third wall crossed by $\ga_i$ (for $i=1,2$), 
namely the first wall immediately after the vertex $p_i$.
Then $W_i$ is either equal to or precedes the first wall crossed by 
$\ov{A_i'}$. Now, if $W_1\cap W_2\neq\emptyset$ then $\ga_1$, $\ga_2$, $W_1$,
and $W_2$ must together bound a flat rectangle isometrically embedded in $\D$
(just as in the proof of Lemma \ref{supergeodesic}). Moreover, this rectangle
must contain the vertex $v_0$ and at least one rank 2 vertex of $K_0$ in its interior.
This contradicts the fact that the shortest closed circuit in the link of
any rank 2 vertex of $\D$ has length at least $3\pi$ (cf, the proof of Lemma 
\ref{supergeodesic} and Lemma \ref{ShortCycles} to follow).
Thus $W_1$ and $W_2$ are disjoint
with $K_0$ lying between them. Moreover $W_i$ separates $U_i$ from $K_0$, for
each $i=1,2$. So each $W_i$ separates $U_1$ from $U_2$.
\end{proof}

It follows from the claim that $\W(A'_1)\cup\W(A'_2)$ is a 
set of mutually disjoint walls all of
which cross $\ov\Sig_3$. Moreover, we can describe $\Sig_3$ as a union of two
$\V$--paths $B_1$ and $B_2$, where $\W(A'_i)\subset \W(B_i)$ for each $i=1,2$, and
such that $B_1$ and $B_2$ overlap in at most 2 edges (ie, at most 2 type $\V$ 
vertices and one type $\F$ vertex in common). Applying Lemma \ref{Vpaths}(ii), we have
$L(A'_i)\leq L(B_i)-2$, for each $i=1,2$, and therefore  
\[
L(A_1\cup A_2)=L(A'_1)+L(A'_2)\leq L(B_1)+L(B_2)-4 \leq L(\Sig_3) - 2\,.
\]
Note that $\Sig$ may be written $(P,\Sig_1,Q_1,\Sig_3,Q_2,\Sig_2,P)$, where 
$\Sig_i$ denotes the longest $\V$--subpath of $\Sig$ such that $\ov\Sig_i$ is
a subpath of $\sig_i$. The above argument shows that $L(A_1\cup A_2) < L(\Sig_3)$.
On the other hand, by Lemma \ref{Vpaths}(i), we have $L(A_i)\leq L(\Sig_i)$ for $i=1,2$, and
therefore, $L(A_1\cup A_2)< L(\Sig_2,P,\Sig_1)$. It follows that $(Q_1,A_1,A_2,Q_2)$
is a short circuit for $\Sigma$. This contradicts the hypothesis that $\Sig$ is 
minimal, and completes the proof of Proposition \ref{first}.

\subsection{Proof of Proposition \ref{second}}

Recall that, for $s\in V(\Delta)$, we denote by $W_s$ the wall in $\D$ perpendicular 
to the edge $[V_\emptyset , V_s]$, and by $W_s^+$ the half-space bounded by $W_s$
and containing the vertex $V_s$. 
We observe that any rank 2 vertex $p$ lies in at least one of the
half-spaces $W_s^+$, for $s\in V(\Delta)$,
and never more than two. If $e=\{ s,t\}$, then $V_e\in W_s^+\cap W_t^+$.
We shall say that a rank 2 vertex \emph{projects to the edge $e=\{s,t\}$} 
if it lies in the region  $W_s^+\cap W_t^+$. Note that if a rank 2 vertex of some
$F_s$, for $s\in V(\Delta)$, projects onto an edge then 
that edge contains $s$, simply because $F_s\subset W_s^+$. 

Let $\ga$ be any path in $\mF$ which starts and ends at rank 2 vertices 
which project to edges, say $e,f\in E(\Delta)$. From the sequence of rank 2 
vertices visited by $\ga$ choose a subsequence $(p_1,p_2,..,p_n)$ 
such that $p_1$ projects to $e$,  $p_n$ projects to $f$ and, 
for $i=1,..,n$, the vertex $p_i$ projects to an edge $e_i=\{s_{i-1},s_i\}$, for 
some sequence $s_0,s_1,\dots ,s_n\in V(\Delta)$.
By passing to a subsequence if necessary, 
we may suppose that the sequence of edges $(e_1,..,e_n)$ describes a
simple path in $\Delta$ (from $s_0$ to $s_n$), where $e_1=e$, $e_n=f$.
This leads to the definition of the
following simple $\V$--path in $\what\Delta\subset \Theta$:
\[
P(\ga)=(V_e=V_{e_1},F_{s_1},V_{e_2},F_{s_2},..,F_{s_{n-1}},V_{e_n}=V_f)\,.
\]
Now suppose that the original path $\ga$ was given as $\ga=\overline A$ for 
some $\V$--path $A$ in $\Theta$. Consider an arbitrary edge $\{s,t\}\in E(\Delta)$. 
If $\ga$ ever enters the region $W_s^+\cap W_t^+$ then we claim that 
$A$ has a type $\V$ vertex which lies in $W_s^+\cap W_t^+$. 
If not, then some type $\F$ vertex of $A$, say $Q\in\F$, must contain 
a geodesic segment which enters and leaves the region. 
But this is impossible since $W_s$ and $W_t$ intersect orthogonally, 
and $Q$ intersects orthogonally with any wall that it encounters.  
It follows that, for each $i=1,..,n$, there is at
least one vertex of type $\V$ in $A$ which projects to $e_i$. 
Therefore $L(P(\ov A))=L(P(\ga))\leq L(A)$.

Finally we claim that if a basic circuit $\Sig$ in $\what\Delta$
admits a short-circuit $B$ in $\Theta$ 
then it admits a short-circuit $B'$ in $\what\Delta$, for if $A$ denotes the longest
subpath of $B$ which is a $\V$--path then replacing $A$ with $P(\ov A)$ gives the
desired path. Proposition \ref{second} now follows.


\section{Minimal circuits, chunks, and isomorphisms of $\Theta$}
\label{sect:ChunkEquiv}

Throughout this section we suppose that $\Delta,\Delta'\in\Cal G$ (CLTTF defining graphs)
 and write $\Theta=\Theta(\Delta)$ and $\Theta'=\Theta(\Delta')$.
We shall be concerned with graph isomorphisms  $\Theta\to\Theta'$ which
map $\F$ to $\F'$ and $\V$ to $\V'$. We call such an isomorphism a
\emph{$\V\F$--isomorphism} of $\Theta$. 

Our objective is to describe a 
decomposition of the graph $\Theta$ which is closely related to the structure of the
Deligne complex $\D$ and yet is canonical in the sense that it is respected
by any $\V\F$--isomorphism of $\Theta$.
The main technical idea is that the pieces of the decomposition (called ``chunks'')
may be characterised by studying the minimal circuits of $\Theta$ introduced 
in Section \ref{sect:Circuits}.

\begin{defn}[{\rm(}Chunks of $\Theta$\/{\rm)}]
Let $A$ be a connected full subgraph of $\Delta$. 
We shall say that $A$ is \emph{indecomposable} if, for every
decomposition, $\Delta=\Delta_1\cup_T\Delta_2$, of $\Delta$ over a
separating edge or vertex $T$, either $A\subset\Delta_1$ or $A\subset\Delta_2$.
By a \emph{chunk} of $\Delta$ we mean a maximal indecomposable 
(connected and full) subgraph of $\Delta$.
Clearly, any two chunks of $\Delta$ intersect, if at all, along a single
separating edge or vertex. A chunk of $\Delta$ shall be said to be \emph{solid}
if it contains a simple closed circuit of $\Delta$. It is easy to see that
any chunk is either solid or consists of just a single edge of the graph.

Recall that the graph $\what\Delta$ may be viewed both as a subset of $\Delta$
and as a subgraph of $\Theta$. 
By a \emph{(solid) chunk} of $\what\Delta$ we mean the intersection 
of $\what\Delta$ with a (solid) chunk of $\Delta$. 
This defines a subgraph of $\what\Delta$
which shall be thought of as lying inside $\Theta$.
Finally, we define a \emph{(solid) chunk} of $\Theta$ to be any translate 
of a (solid) chunk of $\what\Delta$ by an element of $G(\Delta)$.
A chunk of $\Theta$ is said to be \emph{fundamental} if it lies in the
fundamental subgraph $\what\Delta$. 

Note that a chunk of $\Theta$ is solid if and only if it contains a 
simple closed circuit of $\Theta$.
Any non-solid chunk consists of a path in $\Theta$  of the form $(gF_s,gV_e,gF_t)$
for $g\in G(\Delta)$ and $e=\{s,t\}\in E(\Delta)$ 
where $s$ and $t$ are both separating in $\Delta$,
or of the form $(gF_s,gV_e)$ if $s$ is separating but $t$ is terminal.
%
\end{defn}

\begin{lemma}\label{union-of-circuits}
Every minimal circuit of $\Theta$ is contained in a unique solid
chunk and each solid chunk of $\Theta$ is the union of the minimal circuits
which it contains.
\end{lemma}

\begin{proof}
It is a straightforward exercise to prove the corresponding statements for 
minimal circuits of $\what\Delta$ (or $\Delta$).
The Lemma is then a consequence of Propositions \ref{first} and \ref{second},
which establish that the minimal circuits of $\Theta$ are precisely the translates
of the minimal circuits of $\what\Delta$,  together with the fact that
each basic circuit lies in a unique translate 
of $\what{\Delta}$ -- cf Lemma \ref{two-circs}(i) below.
\end{proof}




\begin{defnnum}[{\rm(}Equivalence of minimal circuits, $\sim$\/{\rm)}]\label{circequiv}
We say that minimal circuits $\Sig$ and $\Sig'$ of $\Theta$ are \emph{equivalent},
written $\Sig\sim\Sig'$, if they are related by a finite sequence of 
the following type of elementary equivalence:
\begin{description}
\item [$\bullet$] $\Sig\sim\Sig'$ {\it in one step} if $\Sig$ and $\Sig'$
share a common subpath $(S,V,T)$, with $V\in\V$ and $S,T\in\F$, 
and there exists a sequence of minimal circuits 
\[
\Sig =\Sig_0,\,\Sig_1,\,\Sig_2,\dots ,\,\Sig_k = \Sig'
\] 
such that, for $i=1,..,k$, the circuits $\Sigma_i$ and $\Sig_{i-1}$ share 
a common subpath $(S,V_i,T_i)$ where $V_i\neq V$.
\end{description}
\end{defnnum}

\noindent It is easily seen that the equivalence relation just defined is 
necessarily respected by any $\V\F$--isomorphism of $\Theta$. Our interest in
this equivalence relation lies in the following Proposition.

\begin{prop}\label{chunkequiv}
Let $\Theta,\Theta'$ denote fixed set graphs of CLTTF type.
Suppose that $\Sig_1,\Sig_2$ are minimal circuits in $\Theta$. Then $\Sig_1\sim\Sig_2$ 
if and only if $\Sig_1$ and $\Sig_2$ belong to the same chunk of $\Theta$. 
Consequently, any $\V\F$--isomorphism $\Theta\to\Theta'$ maps each solid chunk of
$\Theta$ onto a solid chunk of $\Theta'$ (inducing a bijection between the solid chunks 
of $\Theta$ and those of $\Theta'$).
\end{prop}

\begin{remark}
Note that, by combining Lemma \ref{union-of-circuits} with the above Proposition, any solid
chunk of $\Theta$ may be described as just the union of a certain equivalence class of
minimal circuits.  
\end{remark}

The remainder of this section is devoted to the proof of Proposition \ref{chunkequiv}.
The Proposition is a consequence of Lemmas \ref{necessity} and \ref{sufficiency} proved 
in the following two Subsections.

\subsection{Equivalent circuits are in the same chunk}

Recall that if $e=\{ s,t\}\in E(\Delta)$ we set $x_e$ to be the group
element expressed by the word $sts..$ of length $m_e$.
Thus $z_e=x_e$ if $m_e$ is even, and $z_e=x_e^2$ if
$m_e$ is odd. In the latter case, conjugation by $x_e$ exchanges the two generators 
$s$ and $t$. The element $x_e$ generates the \emph{quasi-centre} of $G(e)$, the
subgroup of elements which respect the set $\{s,t\}$ by conjugation.
Note also that, when viewing the action of $G$ on the Deligne complex $\D$,
we have $\Fix(x_e)=\Fix(z_e)=V_e$ 
-- cf Lemma \ref{fixedsets}(ii).  

\begin{defn}[{\rm(}Basic solid subset\/{\rm)}]
We define a \emph{basic solid subset} of $\Theta$ to be any translate, by an
element of $G(\Delta)$, of a subgraph of $\what\Delta$ which contains at least one
simple circuit. It follows from Lemma \ref{two-circs}(i) below that there is a 
well-defined function 
\[
\be \co  \{ \text{ basic solid subsets of }\Theta\, \} \to G(\Delta)
\]
such that $\be(\Sig)$ is the unique group element for which $\be(\Sig)\what\Delta$
contains $\Sig$. 
\end{defn}

\begin{lemma}\label{two-circs}
Let $\Sig_1,\Sig_2$ denote basic solid subsets of $\Theta$, and suppose that
these may be written $\Sig_1=\be_1A_1$ and $\Sig_2=\be_2A_2$ 
for $\be_1,\be_2\in G(\Delta)$ and $A_1,A_2\subset \what\Delta$.
\begin{enumerate}
\item[\rm(i)] If $\Sig_1$ and $\Sig_2$ share a common path of length $\geq 3$
(so at least a path $(V',F,V,F')$ with $V,V'\in\V$ and $F,F'=\F$)
then $\beta_1=\be_2$.
\item[\rm(ii)] If $\Sig_1$ and $\Sig_2$ share a common path $(V,F,V')$ where
$V,V'\in\V$ and $F=\be_1F_s$ for some $s\in V(\Delta)$,  then
$\be_2=\be_1h$ where $h\in \< s\>$.
\item[\rm(iii)] If $\Sig_1$ and $\Sig_2$ share a common path $(F,V,F')$ where
$F,F'\in\F$ and $V=\be_1V_e$ for some $e\in E(\Delta)$,  then
$\be_2=\be_1h$ where $h\in \< x_e\>$.
\end{enumerate}
\end{lemma}

\begin{proof}
Note that (i) is a consequence of (ii) and (iii) together with the fact that the cyclic
groups $\< s\>$, for $s\in V(\Delta)$, and $\<x_e\>$, for $e\in E(\Delta)$, intersect
trivially. (Consider a common path $(V',F,V,F')$).

We assume for simplicity, in each case, that $\be_1=1$.
Suppose that, as in case (ii), $\Sig_1$ and $\Sig_2$ share a common path 
$(V_e,F_s,V_f)$ where $e,f\in E(\Delta)$ and $s\in V(\Delta)$ with $e\cap f=\{s\}$.
Then, since $V_e$ and $V_f$ lie in distinct $G$--orbits of the action on $\D$,
we must have $\be_2(V_e)=V_e$ and $\be_2(V_f)=V_f$. 
But then $\be_2\in G(e)\cap G(f)=\< s\>$.

If, as in case (iii), $\Sig_1$ and $\Sig_2$ share a common path 
$(F_s,V_e,F_t)$ where $e=\{s,t\}\in E(\Delta)$, then $\be_2(V_e)=V_e$ and 
$\be_2(\{F_s,F_t\})=\{F_s, F_t\}$. Thus $\be_2\in G(e)$ and conjugation by $\be_2$ 
preserves the set $\{s,t\}$. In other words, $\be_2$ lies in the ``quasi-centre'' 
of $G(e)$ which is generated by $x_e$.
\end{proof}
 
\begin{lemma}\label{freegroup}
Let $G=G(\Delta)$ be a large-type triangle-free Artin group.
Then the set $\{x_e\,:\,e\in E(\Delta)\}$ freely generates
a free subgroup of $G(\Delta)$.
\end{lemma}

\begin{proof}
Let $\Cal X$ denote the complex of groups with underlying complex $K$, a vertex
group $\< x_e\>$ at each vertex $V_e$, for $e\in E(\Delta)$, and all other vertex groups
trivial. Then $\pi_1(\Cal X)$ is simply the free product of the cyclic groups $\< x_e\>$.
(In fact the the universal covering complex $\wtil{\Cal X}$ admits an equivariant deformation
retraction onto a Bass-Serre tree for this free product). The obvious map from 
$\pi_1(\Cal X)$ onto the subgroup of $G$ generated by the set $\{x_e\,:\,e\in E(\Delta)\}$ 
 is associated with
an equivariant map  $\psi\co  \wtil{\Cal X}\to\D$ with image the union of  
translates $gK$ of the fundamental region $K$, for $g\in\< x_e:e\in E(\Delta)\>$.
We equip both spaces $\wtil{\Cal X}$ and $\D$ with the natural cubical metric $d_C$,
and observe that, for each edge $e\in E(\Delta)$, the union $\bigcup\limits_{k\in\Z} x_e^k K$ 
 is a convex subset of $(\D,d_C)$ -- one simply
needs to check local convexity at the vertex $V_e$, using the hypothesis that $G$ be large type. 
It follows that the map $\psi$ is locally isometric, and hence a globally isometric embedding
(since the image lies in a CAT(0) space). The map $\pi_1(\Cal X)\to G$ is therefore injective.
\end{proof}

\begin{remark}
Note that the convexity statement made in the above proof is not true with respect 
to the Moussong metric, which explains the triangle-free hypothesis. 
Nevertheless, the above result seems likely to be true for an arbitrary large type Artin group, 
although we do not have an obvious proof to hand. 
\end{remark}

\begin{lemma}\label{betas} 
Let $\Sig,\Sig'$ denote minimal ciruits in $\Theta$.
If $\Sig\sim \Sig'$ then $\be(\Sig)=\be(\Sig')$. Moreover, if $\Sig\sim\Sig'$
in one step then the sequence $\{\Sig_i\}_i$ in Definition \ref{circequiv}
may be chosen so that $\be(\Sig_i)=\be(\Sig)$ for all $i$.
\end{lemma}

\begin{proof} 
It suffices to verify the statement concerning one step equivalence.
We suppose therefore that, as in the definition of one step equivalence 
(Definition \ref{circequiv}), there exists a sequence
$\{\Sig_i\}$, $i=0,1,..,k$, of minimal circuits with $\Sig_0=\Sig$ and $\Sig_k=\Sig'$,
a subpath $(S,V,T)$ common to $\Sig_0$ and $\Sig_k$, 
and a sequence of subpaths $(S,V_i,T_i)\subset \Sig_i\cap\Sig_{i-1}$
with $V_i\neq V$  for $i=1,..,k$. 
We shall write $\be_i=\be(\Sig_i)$ for $i=0,1,..,k$. We also let $e,e_1,..,e_k$
denote the edges in $\Delta$ such that $V=\be_0V_e$ and $V_i=\be_iV_{e_i}$ for $i=1,..,k$.
Note that these edges are not necessarily distinct. 
Finally, we suppose that the sequence $\{\Sig_i\}_i$ is chosen so as to minimise 
the length $k$. With this assumption we make the following claim:

\begin{claim} 
If $e_i=e_j$ for some $1\leq i<j\leq k$ then $\be_i\neq \be_{j-1}$.
\end{claim}

\begin{proof}
Fix $j\in\{1,..,k\}$, and set $e_j=\{ s,t\}$. Then 
$(S,V_j,T_j)=\be_j(F_s,V_{e_j},F_t)$. Since this subpath is common 
to both $\Sig_j$ and $\Sig_{j-1}$ we have, by Lemma \ref{two-circs}(iii),
that $\be_{j-1}\inv\be_j\in\< x_{e_j}\>$, and therefore that 
$\be_{j-1}(V_{e_j})=\be_j(V_{e_j})=V_j$ and 
$\be_{j-1}(\{F_s,F_t\})=\be_j(\{F_s,F_t\})=\{T_j,S\}$.
Now suppose, by way of contradiction, that $e_i=e_j$ and 
$\be_i=\be_{j-1}$, for some $i<j$. Then, by the previous 
observation, we have  $V_i=\be_i(V_{e_i})=\be_{j-1}(V_{e_j})=V_j$, 
and  $\{ T_i,S\}=\{ T_j,S\}$ similarly. That is $T_i=T_j$ and $V_i=V_j$.
But then one can simply remove the circuits $\Sig_i,..,\Sig_{j-1}$ from
the sequence to obtain a shorter sequence satisfying the one 
step equivalence of Definition \ref{circequiv}, contrary to our 
choice of a shortest length sequence. 
\end{proof}

Applying Lemma \ref{two-circs} (iii), we have a sequence of elements
$h_i=\be_{i-1}^{-1}\be_i$, for $i=1,\dots,k$, such that $h_i\in\<x_{e_i}\>$ 
for each $i$. Also, writing $h=\be_0^{-1}\be_k$ we have that $h\in\< x_e\>$
(again by Lemma \ref{two-circs} (iii)\,) and
\[
h=h_1h_2h_3\cdots h_k\,.
\]
Note that by Lemma \ref{freegroup} the elements $x_d$ for $d\in E(\Delta)$ 
generate a free group $L$.
The above Claim implies that deleting all trivial syllables $h_i=1$ from the 
expression $h_1h_2\dots h_k$ yields a reduced form for the element $h$ with 
respect to the structure
of $L$ as a free product of cyclic groups $\<x_d\>$ for $d\in E(\Delta)$.
(For, if $h_{i+1}= h_{i+2}=\cdots =h_{j-1}=1$ then 
$\be_i^{-1}\be_{j-1}=h_{i+1}\cdots h_{j-1}=1$ and so, by the Claim,
$e_i\neq e_j$ which implies that $h_i,h_j$ do not belong to the same free factor).
However, since $h$ lies in the free factor $\<x_e\>$, it follows 
that either $h=1$ or there is a unique $r\in\{1,..,k\}$ for
which $h_r$ is nontrivial with $h_r=h$ and $e_r=e$.
In the latter case we would have $\be_0^{-1}\be_r=h\in\< x_e\>$
which would imply that $V_r=\be_r(V_e)=\be_0(V_e)=V$, a contradiction. 
Thus  $h=1$, and consequently $h_i=1$, and so $\be_i=\be_0$, 
for all $i=1,..,k$.
\end{proof}

\begin{lemma}\label{necessity}
Let $\Sig,\Sig'$ denote minimal circuits in $\Theta$. If $\Sig\sim\Sig'$ 
then $\Sig$ and $\Sig'$ belong to the same chunk of $\Theta$. 
\end{lemma}

\begin{proof}
Suppose that $\Sig\sim\Sig'$ by a one step equivalence as described in
Definition \ref{circequiv}. Without loss of generality we may suppose that
$\Sig\subset\what\Delta$. Then
we have a sequence $\Sig=\Sig_0,\Sig_1,..,\Sig_k=\Sig'$
of minimal circuits as in the definition where, by Lemma \ref{betas},
we may suppose that all circuits $\Sig_i$ lie in $\what\Delta$.
Considering, for simplicity, all circuits as
circuits in $\Delta$, we have that $\Sig,\Sig'$ have an edge $e$ in common while, for 
each $i=1,..,k$, the circuits $\Sig_i$ and $\Sig_{i-1}$ meet along an edge $e_i$
different from $e$. 
If $\Sig$ and $\Sig'$ were to lie in distinct chunks then 
the edge $e$ would have to separate $\Delta$ into two pieces $A$ and $B$ containing
$\Sig$ and $\Sig'$ respectively. However, the fact that $e_i\neq e$ would then 
imply that $\Sig_i\subset A$ if and only if $\Sig_{i-1}\subset A$, and therefore
that $\Sig$ and  $\Sig'$ both lie in $A$, a contradiction. 
Thus $\Sig$ and $\Sig'$ lie in the same chunk.
\end{proof}

\subsection{Circuits in the same chunk are equivalent}

For this part of the proof it will be convenient to only consider  circuits  in
the original defining graph $\Delta$ (before subdivision of the edges). We shall say that 
an edge or vertex $T$ of $\Delta$ \emph{separates} subsets $A$ and $B$ if $\Delta$ can be
written $\Delta_1\cup_T\Delta_2$ with $A\subset\Delta_1$ and $B\subset \Delta_2$.
Note that this terminology implies that $T$ is a separating edge or vertex of $\Delta$.

\begin{lemma}\label{stepI}
If $\Delta$ is a connected graph and  $\Sig,\Sig'$ are minimal 
circuits then either there is a vertex of $\Delta$ which separates them, or
they are joined by a sequence $\Sig=\Sig_0,\Sig_1,\Sig_2,..,\Sig_n=\Sig'$
of minimal circuits of $\Delta$ such that $\Sig_i$ and 
$\Sig_{i-1}$ have an edge in common, for each $i=1,..,n$.
\end{lemma}

\begin{proof}
We first observe that there is necessarily a simple edge path in $\Delta$
whose first edge lies in $\Sig$ and whose last edge lies in $\Sig'$. 
Moreover, since no vertex of this path can separate $\Sig$ and $\Sig'$,
we may suppose that no two subsequent edges along this path are 
separated in $\Delta$ by their common vertex.
It will now suffice to show that if $e,e'$ are edges of $\Delta$ with a common
vertex that does not separate them 
then there exists a sequence of minimal circuits $\Sig_0,\Sig_1,..,\Sig_n$ 
such that $e\subset\Sig_0$, $e'\subset\Sig_n$ 
and $\Sig_i\cap\Sig_{i-1}$ contains an edge for each $i=1,..,n$.
Suppose we have $e=\{t,s\}$ and $e'=\{s,t'\}$.
Since $s$ does not separate the two edges, there must exist a simple path $\al$
from $t$ to $t'$ which does not pass through $s$. 
Concatenating $\al$ with the edges $e,e'$ yields a simple circuit
$C$ through $e,e'$. If $C$ is minimal then there is nothing left to prove. 
Otherwise, we may find a short-circuit  which decomposes $C$ into
circuits $C_1$ and $C_2$ of strictly shorter length. 
Either $e,e'$ both still lie in the same 
circuit, $C_1$ say,  in which case we replace $C$ with $C_1$, or $e\subset C_1$ 
and $e'\subset C_2$ and there is an edge $e''=\{s,t''\}$ common to $C_1$ and $C_2$ 
and adjacent to $s$, in which case we replace $C$ with the sequence $C_1,C_2$. 
In either case, the desired result follows by induction on the length of $C$.
\end{proof} 

\begin{lemma}\label{stepII}
If minimal circuits $\Sig,\Sig'$ of a graph $\Delta$ meet along an edge $e$ then
either $e$ separates $\Sig$ from $\Sig'$, or there exists a sequence
$\Sig=\Sig_0,\Sig_1,\Sig_2,..,\Sig_k=\Sig'$ of minimal circuits such that, for each
$i=1,..,k$, the circuits $\Sig_i$ and $\Sig_{i-1}$ have a common edge $e_i$ which
is adjacent to but not equal to $e$.
\end{lemma}

\begin{proof}
Let $e\in E(\Delta)$, and write $S(e)$ for the set of vertices of $\Delta$
which are not themselves endpoints of $e$, but which are adjacent
along an edge to one or other endpoint of $e$. For $u,v$ distinct
vertices in $S(e)$, we define a \emph{joining arc} from $u$ to $v$, to
be any simple edge-path from $u$ to $v$ in $\Delta$ which does not pass
through either of the endpoints of $e$. If $\al$ is a joining arc from $u$ to $v$,
then write $\Sig(\al)$ for the simple circuit formed from $\al$ and the unique
shortest edge-path (of length 2 or 3) from $u$ to $v$ passing through one
or both of the endpoints of $e$. Note that joining arcs are thus in bijective
correspondence with the simple closed circuits in $\Delta$ which intersect the
edge $e$ (at one or both of its endpoints). The joining arc $\al$ shall be said to be
\emph{minimal} if $\Sig(\al)$ is a minimal circuit.

\begin{sublemma}
Let $e\in E(\Delta)$, $u,v$ distinct vertices in $S(e)$, and $\al$ a joining arc
from $u$ to $v$. Then there exists a sequence $u=u_1,u_2,..,u_k=v$
of elements of $S(e)$ and minimal joining arcs $\al_i$
from $u_i$ to $u_{i+1}$, for $i=1,..,k-1$.
\end{sublemma} 

\begin{proof}
Suppose $\al$ is not minimal. Then $\Sig(\al)$ admits a short circuit $\sigma$
which decomposes $\Sig(\al)$ into two simple circuits $C_1$ and $C_2$ each of
strictly smaller length than $\Sig(\al)$. Either $\sigma$ is
disjoint from $e$, or one of its endpoints is also an endpoint of $e$.
In the first case, one of the two smaller circuits, $C_1$ say,
contains $u,v$ and at least one vertex of $e$, and is thus equal to $\Sig(\al')$
where $\al'$ is a joining arc from $u$ to $v$.
In the second case, both circuits $C_1$ and $C_2$ pass though a vertex of
$e$ and we
may suppose that $C_1$ contains $u$ and $C_2$ contains $v$.
In fact, if $w$ denotes the interior vertex of $\sigma$ closest to $e$
(so that $w\in S(e)$), then $C_1=\Sig(\al_1)$ where $\al_1$
is a joining arc from $u$ to $w$,
and $C_2=\Sig(\al_2)$ where $\al_2$ is a joining arc from $w$ to $v$.
The Sublemma now follows by induction on the length of $\Sig(\al)$.
\end{proof}

We shall now complete the proof of Lemma \ref{stepII}.
We have an edge $e$ which is common to minimal circuits $\Sig$ and $\Sig'$.
Suppose that $e$ does not separate $\Sig$ and $\Sig'$.
Take distinct vertices $u,v\in S(e)$ such that $u\in \Sig$ and $v\in \Sig'$.
since $e$ does not separate the two circuits, there is a joining arc from $u$ to $v$ and, by
the Sublemma, a sequence $u=u_1,..,u_k=v$ of
vertices in $S(e)$, and minimal joining arcs $\al_i$ from $u_i$ to $u_{i+1}$,
for $i=1,..,k-1$. Write $\Sig_0=\Sig$, $\Sig_i=\Sig(\al_i)$,
for $i=1,..,k-1$, and $\Sig_k=\Sig'$.
Now, for each $i=1,..,k$, $\Sig_i\cap \Sig_{i-1}$ contains at least one edge
$e_i$ which is adjacent to $e$ and has the vertex $u_i$ as one
of its endpoints. In particular $e_i\neq e$. 
\end{proof}

\begin{lemma}\label{sufficiency}
Let $\Sig,\Sig'$ denote minimal circuits in $\Theta$.
If $\Sig$ and $\Sig'$ belong to the same chunk of $\Theta$ then $\Sig\sim\Sig'$. 
\end{lemma}

\begin{proof} 
Without loss of generality we may suppose that $\Sig$ and $\Sig'$ 
are fundamental minimal circuits lying in the same chunk of $\what\Delta$.
We may also, for simplicity, consider these circuits as circuits 
of the graph $\Delta$. (Note that by Proposition \ref{second}
circuits which are minimal in $\Delta$, equivalently $\what\Delta$,
are also minimal as circuits of $\Theta$).

Suppose firstly that $\Sig$ and $\Sig'$ share an edge 
$e\in E(\Delta)$. Since $\Sig$, $\Sig'$ lie in the same chunk of $\Delta$, 
$e$ cannot separate them and so Lemma \ref{stepII} gives a 
sequence of minimal circuits
$\Sig=\Sig_0,\Sig_1,\Sig_2,\dots,\Sig_n=\Sig'$ such that each $\Sig_i$ shares 
an edge $e_i$ with the previous and $e_i$ is adjacent to but distinct from $e$.
Let $e=\{s,t\}$. If  $e_i\cap e=\{s\}$ for all $i=1,..,n$, then 
$\Sig_0\sim\Sig_n$ in one step. Suppose then that $e_i\cap e=\{s\}$, for $i=1,..,k$,
and $e_{k+1}\cap e=\{t\}$, for some $k<n$. Then  $\Sig_k$ contains
the path $e_k,e,e_{k+1}$. It follows, since $e\subset\Sig_k$,
that $\Sig_0\sim\Sig_k$ in one step. 
By a straightforward induction we conclude that $\Sig\sim\Sig'$.

More generally, if $\Sig$ and $\Sig'$ lie in the same chunk of $\Delta$ then,
since no vertex can separate them, Lemma \ref{stepI} gives a sequence of minimal circuits
$\Sig=\Sig_0,\Sig_1,\Sig_2,\dots,\Sig_n=\Sig'$ in $\Delta$ where each shares 
an edge with the next. Moreover, since any two adjacent chunks meet along at most a
single edge, we may suppose (by passing to a subsequence) that all circuits $\Sig_i$ 
lie in the same chunk of $\Delta$.  But, by the preceding argument, this means that
$\Sig_{i-1}\sim\Sig_i$ for each $i$, and so $\Sig\sim\Sig'$.  
\end{proof}


\section{Chunk rigidity in $\Theta_W$ the graph of fixed sets in $\D_W$}

We obtain results exactly analogous to those of the last two sections in the 
context of the action of a CLTTF Coxeter group $W=W(\Delta)$ on its Davis complex $\D_W$.
Throughout this section we shall view $\D_W$ as a subcomplex of 
$\D$ via the map induced by the Tits section. We note that $\D_W$ is geodesically 
convex in $\D$ with respect to either the Moussong metric, or the cubical metric.
We shall  write $\ov{x}$ for the image in $W$ of an element $x\in G$,
and  $i_W\co W\to G$ for the Tits section.

For $s\in V(\Delta)$ we set $H_s:= F_s\cap\D_W$. Then $H_s$ is the fixed set in $\D_W$
of the standard reflection $\ov s$, and may be thought of as a ``hyperplane" in the 
Davis complex. Note that $H_s$ is infinite if and only if $F_s$ is unbounded.

\begin{defn}
We define the following sets 
\[
\begin{aligned}
\V_W &=\{ \text{ singletons }\{wV_e\} : w\in W\,,\ e\in E(\Delta)\}\,,\text{ and }\\
\F_W &= \{ \text{ infinite hyperplanes } wH_s :  w\in W\,,\ s\in V(\Delta)\}\,.
\end{aligned}
\]
We define the graph $\Theta_W$ to be the bipartite graph with vertex set
$\V_W\cup\F_W$ and edges $(V,H)$ whenever $V\in\V_W$, $H\in \F_W$ and $V\subset H$.
\end{defn}
 
The elements of $\V_W$ are characterized as the fixed sets of the maximal finite 
subgroups of $W$, namely the conjugates of dihedral groups 
$\<\ov{s},\ov{t}\>=\Stab(V_e)$ for all edges $e=\{s,t\}\in E(\Delta)$.
Any two maximal finite subgroups $D_1$ and $D_2$ in $W$ intersect nontrivially 
if and only if their fixed points lie on a common hyperplane $wH_s$ for some 
$s\in V(\Delta)$ and $w\in W$, in which case $D_1\cap D_2=\{1, w\ov{s}w\inv \}$.
This condition is also equivalent to the statement that $\Fix(D_1)$ and $\Fix(D_2)$
both lie in $gF_s$ where $g=i_W(w)$ is the image of $w$ under the 
Tits section $W\hookrightarrow G$. Thus we have the following.

\begin{lemma}\label{ThetaW}$\phantom{99}$
\begin{itemize}
\item[\rm(i)] Any isomorphism $\phi\co W\to W'$ induces a graph isomorphism 
$\Phi\co  \Theta_W\to\Theta_W'$ such that, for all $X\in\V_W\cup\F_W$, we have
$\Stab(\Phi(X))=\phi(\Stab(X))$. 
Moreover, this is a ``$\V\F_W$--isomorphism": $\Phi(\V_W)=\V'_W$ 
and $\Phi(\F_W)=\F'_W$. 
\item[\rm(ii)] The inclusion of $\D_W$ into $\D$ induces a natural inclusion
$\Theta_W\hookrightarrow \Theta$ (such that $\{wV_e\}\mapsto \{i_W(w)V_e\}$ and
$wH_s\mapsto i_W(w)F_s$). 
\end{itemize}
\end{lemma}

We shall identify $\Theta_W$ with its image in $\Theta$ under the map of
part (ii) of the Lemma. Note that the fundamental subgraph $\what\Delta$ lies inside
$\Theta_W$, ie, $\what\Delta\subset\Theta_W\subset \Theta$.
Note also that any graph isomorphism $\Phi\co  \Theta_W\to\Theta_W'$ will respect the
family of circuits of $\Theta_W$ which are minimal \emph{as circuits of $\Theta_W$}. 

\begin{lemma}
Let $\Sigma$ denote a simple circuit in $\Theta_W$. Then $\Sigma$ is minimal 
as a circuit of $\Theta_W$ if and only if it is minimal as a circuit of $\Theta$.
In particular, the minimal circuits of $\Theta_W$ are precisely the translates
of minimal circuits of the fundamental subgraph $\what\Delta$ by elements of $W$.
\end{lemma}

\begin{proof}
In Section \ref{sect:TwoProps},  Proposition \ref{first} was proved by showing
that if a circuit $\Sigma$ fails to lie in a fundamental region of $\Theta$
then it admits a short circuit, and the only short circuits exhibited throughout
the proof were constructed from one or more chords of the polygon $\ov\Sigma$
(cf Lemma \ref{sigalign} and \figref{ChordsFig}).
If we are given $\Sigma_W$ a simple circuit in $\Theta_W$ then,
since $\D_W$ is a geodesically convex subcomplex of $\D$, the simple closed
curve $\ov\Sigma_W$ and any chord of $\ov\Sigma_W$  are contained in $\D_W$.
It follows that any short circuit produced in the proof of Proposition 
\ref{first} is a short circuit in $\Theta_W$.
Thus any circuit which is minimal as a circuit of $\Theta_W$ is basic.
Note that a circuit in $\what\Delta$ which is minimal in $\Theta_W$ is necessarily
minimal in $\what\Delta$ and therefore, by Proposition \ref{second}, minimal in $\Theta$.
It follows that any circuit which is minimal in $\Theta_W$ is minimal in $\Theta$. 
The converse is obvious.
\end{proof}

We may now define an equivalence relation $\sim_W$ on the set of minimal 
circuits of $\Theta_W$ exactly as per Definition \ref{circequiv}, 
but with reference only to the minimal circuits of $\Theta_W$. 
Clearly, this equivalence relation is respected by any isomorphism of $\Theta_W$.
Let $\Sigma,\Sigma'$ be minimal circuits of $\Theta_W$. Then it follows
immediately from Lemma \ref{betas} that $\Sigma\sim_W\Sigma'$ if and only if 
$\Sigma\sim\Sigma'$ (as circuits in $\Theta$). Thus, by Proposition \ref{chunkequiv},
we have that $\Sigma\sim_W\Sigma'$ if and only if 
$\Sigma$ and $\Sigma'$ lie in the same chunk of $\Theta_W$ where, by \emph{chunk of $\Theta_W$},
we understand any translate of a chunk of the fundamental subgraph $\what\Delta$ by an 
element of $W$. The following is an immediate consequence of this statement
and Lemma \ref{ThetaW}(i) above:

\begin{prop}\label{CoxeterChunks}
Let $W=W(\Delta)$ and $W'=W(\Delta')$ be CLTTF Coxeter groups.
Any isomorphism $\phi\co W\to W'$ naturally induces a $\V\F_W$--isomorphism
$\Theta_W\to \Theta'_W$ which maps the solid chunks of $\Theta_W$ bijectively onto the 
solid chunks of $\Theta'_W$.
\end{prop}


\section{Automorphisms -- Theorems \ref{Thm1}, \ref{Thm2}, and \ref{ThmCoxeter}}

Let $\Delta,\Delta'\in\Cal G$ be CLTTF defining graphs, and write $G=G(\Delta)$,
$G'=G(\Delta')$. We shall also simply write $\Theta$ and $\Theta'$ for the 
corresponding fixed set graphs.

\begin{lemma}\label{induced-from-phi}
Let $\varphi\co  G\to G'$ be an isomorphism. Then:
\begin{itemize}
\item[\rm(i)] 
$\varphi$ induces a graph isomorphism $\Phi\co \Theta\to\Theta'$ 
such that $\Phi(\Fix(C)) = \Fix(\varphi(C))$ for any  CNVA  
subgroup $C<G$. The isomorphism $\Phi$ is a $\V\F$--isomorphism.
\item[\rm(ii)] 
$\varphi$ induces a  label preserving bijection 
$\ov\varphi\co  E(\Delta)\to E(\Delta')$ which is defined uniquely 
such that $G(\ov\varphi(e))$ and $\varphi(G(e))$
are conjugate  subgroups of $G'$ for each $e\in E(\Delta)$.
\end{itemize}
\end{lemma} 

\begin{proof}
{\bf (i)}\qua
The fact that $\Phi$ is well-defined is just a restatement of Proposition \ref{ThetaIsom}.
To see that $\Phi$ is necessarily a $\V\F$--isomorphism, 
observe that maximal  CNVA  subgroups conjugate to $\<s\>$ for 
$s\in V(\Delta)$ are distinguished from those conjugate to $\< z_e\>$ for 
$e\in E(\Delta)$ by the fact that a standard generator $s$ is a primitive element
of $G$, while $z_e$ is not (if $e=\{s,t\}$ then $z_e=(st)^k$ where 
$k=\text{lcm}(m_e,2)/2$). 
\medskip

\noindent{\bf (ii)}\qua
For each $e\in E(\Delta)$ we have that $\Phi(V_e)= g_eV_{e'}$
for some $g_e\in G'$ and $e'\in E(\Delta')$, and therefore 
(since $\Stab(V_e)=G(e)$, etc) $\varphi(G(e))=g_eG(e')g_e\inv$.
It follows from the description of the Deligne complex that
rank 2 vertices $V_e$ and $V_f$ lie in distinct orbits of the group action 
unless $e=f$, and therefore that the stabilizers $G(e)$ and $G(f)$ are non-conjugate
unless $e=f$. Therefore, setting $\ov\varphi(e)=e'$ gives a well-defined 
bijection $\ov\varphi\co  E(\Delta)\to E(\Delta')$. Moreover, since groups 
$G(e)$ and $G(e')$ (for any edges $e$ and $e'$) are non-isomorphic 
unless $m_e = m_{e'}$ the bijection $\ov\varphi$ must be label preserving.
\end{proof}

Note that the proof of the above Lemma  depends heavily on the fact that
$\varphi$ is an \emph{isomorphism}, rather than an abstract commensurator.

Rcall that any Artin group $G(\Delta)$ admits a length 
homomorphism  $\ell\co G(\Delta)\to \Z$ such that $\ell(s)=1$ for each
generator $s\in V(\Delta)$.

\begin{lemma}\label{leafgroup}
Let $e=\{s,t\}\in E(\Delta)$, $m_e\geq 3$, and suppose that $\al\in\Aut(G(e))$ 
such that $\ell(\al(s))= \ell(\al(t))=1$. 
Then $\al$ differs by an inner automorphism of $G(e)$ from either the identity or
the graph automorphism $\tau \co  s\leftrightarrow t$.
\end{lemma}

\begin{proof}
This is an easy consequence of the computation of the automorphism group of
a dihedral type Artin group first performed in \cite{GHMR} (see also 
\cite{ChCr} for a description of the same automorphism group). 
\end{proof}

\begin{lemma}[Chunk Invariance]\label{chunk-invariance}
Let $\varphi\co  G\to G'$ be an isomorphism and suppose that $\ell(\varphi(s))=1$ 
for every generator $s\in V(\Delta)$ of $G$. Then, for each chunk $A$ of $\Delta$
(solid or otherwise), 
there exists a chunk $A'$ of $\Delta'$, a labelled graph 
isomorphism $\tau_A\co A\to A'$, and an element $g_A\in G'$ such that 
the restriction $\varphi_A$ of $\varphi$ to the subgroup $G(A)$
is given by
\[ 
\varphi_A= g_A\circ \tau_A
\]
(where $g_A$ denotes conjugation by $g_A$ and, by abuse of notation, $\tau_A$
denotes the group isomorphism induced by the graph isomorphism $\tau_A$).  
Note that the mapping $A\mapsto A'$ defines a bijection between the
chunks of $\Delta$ and those of $\Delta'$.
\end{lemma}

\begin{proof}
First suppose that $A$ is a solid chunk of $\Delta$. We may equally view $A$ as a solid
chunk of $\what\Delta$. Proposition \ref{chunkequiv} states that $\Phi$, the induced
$\V\F$--isomorphism of Lemma \ref{induced-from-phi}(i), carries solid chunks of
$\Theta$ to solid chunks of $\Theta'$. That is to say that $\Phi(A)=g_A(A')$ for 
some $g_A\in G'$ and some solid fundamental chunk $A'$ of $\what\Delta'$. 
In fact, the element $g_A$ and the chunk $A'$ are uniquely determined. 
Now, for each edge $e\subset A$, we have $\varphi(G(e))=g_A G(e')g_A^{-1}$ 
(and $\Phi(V_e)=g_A V_{e'}$) for some $e'\subset A'$.
Moreover, if $e$ and $f$ are edges of $A$ and $e\cap f=\{ s\}$ then
$\varphi(\< s\>)=g_A\<s'\>g_A\inv$  where 
$e'\cap f' = \{ s'\}$ (since $\< s\>=G(e)\cap G(f)$).
In fact $\varphi(s)=g_A s'g_A\inv$, since we suppose that $\ell(\varphi(s))=1$. 
Note also that $e'=\ov\varphi(e)$, where
$\ov\varphi \co  E(\Delta)\to E(\Delta')$ is the induced bijection 
of Lemma \ref{induced-from-phi}(ii). 
It now follows that the restriction of $\ov\varphi$
to the edges of $A$ determines a graph isomorphism $\tau_A$ and 
that $\varphi_A= g_A\circ \tau_A$.
 
The case where $A=e$ is a non-solid chunk follows easily by taking $A'=\ov\varphi(e)$
and applying Lemma \ref{leafgroup}.
\end{proof}

\subsection{Proof of Theorem \ref{Thm1}}

We suppose that we are given a group isomorphism $\varphi\co  G\to G'$. 
Our approach will be to compose this isomorphism with known isomorphisms of types
(2)--(4) (inversions, Dehn twist and inner automorphisms and edge twist isomorphisms)
until it is reduced to a graph automorphism (type (1)).

\paragraph{Applying the inversion automorphisms}
We shall use the existence of the standard length homomorphism $\ell\co G'\to \Z$
(such that $\ell(s)= 1$ for all $s\in V(\Delta')$).
For each $e\in E(\Delta)$, we have 
$\varphi(\< z_e\>)=g\< z_{e'}\>g\inv$, where $e'=\ov\varphi(e)$ and $g\in G'$.
 Also, since $\ell(z_{e'})\neq 0$, 
the element $z_{e'}$ is not conjugate to its inverse.
Thus we have a well-defined function $\nu\co  E(\Delta)\to\{\pm 1\}$ such that  
$\varphi(z_e)\sim z_{e'}^{\nu(e)}$.
If $e=\{s,t\}$ where both $\<s\>$ and $\< t\>$ are  CNVA  then both
$\varphi(s)$ and $\varphi(t)$ have absolute length $1$ 
(since $s$ and $t$ are mapped to generators of maximal  CNVA  subgroups).
Since $\varphi$ respects the relation $(st)^k=z_e$ ($k=\text{lcm}(m_e,2)/2$), 
we must therefore have $\ell(\varphi(s))=\ell(\varphi(t))=\nu(e)$. 
This argument applies to
every edge of $\Delta$ with the exception of the even labelled terminal edges,
where the terminal generator does not generate a  CNVA  subgroup. 
By connectedness of $\Delta$ it follows that $\nu$ is constant
on the set $E(\Delta)\setminus\{ \text{ even labelled terminal edges }\}$.
By precomposing $\varphi$  with leaf inversions and a global inversion, as necessary,
we may now suppose that $\nu(e)=1$ for all $e\in E(\Delta)$ and that 
$\ell(\varphi(s))=1$ for every CNVA generator $s\in V(\Delta)$.
(Recall that if $\mu_e$ is a leaf inversion then it fixes all CNVA generators
and $\mu_e(z_e)=z_e^{-1}$).  
Note that, for any edge $e=\{ s,t\}$ the relation $(st)^k=z_e$
implies that if $\ell(\varphi(s))=\nu(e)=1$ then $\ell(\varphi(t))=1$.
Thus, it actually follows that $\ell(\varphi(s))=1$ for every single
generator $s\in V(\Delta)$. 

(In the last statement we are implicitly using the assumption that 
$\Delta$ has at least 3 vertices and is connected (C), 
and so has at least one CNVA generator).

\paragraph{Applying edge twists and Dehn twists}
We suppose from now on that $\ell(\varphi(s))=1$ for all $s\in V(\Delta)$.
Let $B$ denote any connected subgraph of $\Delta$ which is a union
of chunks. We shall show, by induction on the number of chunks in $B$,
that we may arrange (by composing with isomorphisms of type (3) and (4))
that the restriction $\varphi_B$ of $\varphi$ to 
$G(B)$ is induced by a labelled graph isomorphism $\tau_B\co B\to B'$,
for some connected subgraph $B'$ of $\Delta'$ (which is also necessarily
a union of chunks of $\Delta'$). In the case where $B=\Delta$ we have 
$\Delta=\Delta'$ and $\varphi=\tau_B$ a graph automorphism of $G(\Delta)$, completing 
the proof of Theorem \ref{Thm1}.

If we take $B$ to be any single chunk then, by Lemma \ref{chunk-invariance}, 
we may suppose, up to an inner automorphism of $G'$, that $g_B=1$ and $\varphi_B=\tau_B$.

Suppose now that the statement is already proven for some subgraph $B\neq\Delta$, 
and that $A$ is a chunk of $\Delta$ which does not lie in $B$ but which intersects
$B$ nontrivially. We shall prove the statement for $A\cup B$. 
We have $\varphi_B=\tau_B$ for some  graph isomorphism 
$\tau_B \co  B\to B'$ and, by Lemma \ref{chunk-invariance},
we may suppose that $\varphi_A=g_A\circ\tau_A$, 
for some $g_A\in G'$ and $\tau_A\co A\to A'$ where $A'$ denotes a chunk of $\Delta'$.
There are two cases to consider.

Suppose firstly that $A$ and $B$ intersect along an edge $e=\{s,t\}$. 
Then $\varphi_A=g_A\circ\tau_A$ and $\varphi_B=\tau_B$ must agree on the subgroup $G(e)$.
In particular $\tau_A(e)=\tau_B(e)=e'=\{s',t'\}$, and $g_A(\{s',t'\})=\{s',t'\}$.
It follows that $g_A$ lies in $G(e')$ 
(since $g_AG(e')g_A^{-1}=G(e')$ implies that $g_A(V_{e'})=V_{e'}$ and so $g_A\in G(e')$).
More precisely $g_A$ lies in the quasi-centre
$\< x_{e'}\>$ of $G(e')$. By composing with a sequence of edge 
twist isomorphisms  we may now suppose that $g_A=1$, and consequently that
$\tau_A$ and $\tau_B$ agree on $e$ (ie, $\tau_A(s)=\tau_B(s)=s'$ etc). 
But then $\varphi_{A\cup B}$ is induced from a labelled graph isomorphism, as required.
(Note that each of the above edge twists will change $\Delta'$ by a twist 
equivalence. However, they will compose to give a genuine Dehn twist in the
case where the given $g_A$ is central in $G(e')$).

Now suppose that $A\cap B=\{s\}$, for some $s\in V(\Delta)$. 
Write $s',t'$ for the elements of $V(\Delta')$ such that 
$\tau_A(s)=s'$ and $\tau_B(s)=t'$. Since $g_A\circ\tau_A(s)=\tau_B(s)$ we have that
$g_As'g_A\inv=t'$. We may now find a simple edge path $\ga$ from $s'$ to $t'$
which consists of only odd label edges. (If not $s'$ and $t'$ would map to 
distinct cyclic factors of the abelianisation of $G'$ and hence could not 
be conjugate in $G'$).

Note that two chunks of the graph $\Delta$ are separated by a vertex if and only
if they cannot be joined by a sequence of chunks where each has a edge in common
with the next. Since $\ov\varphi\co E(\Delta)\to E(\Delta')$ is a bijection which
preserves chunks it also preserves the above property. Therefore, since $A$
is separated from any chunk in $B$ by the vertex $s$, it follows that 
$A'=\tau_A(A)$ is separated from at least one chunk of $B'=\tau_B(B)$ by some 
vertex $v'\in V(\Delta')$.
(Note that both graph isomorphisms $\tau_A$ and $\tau_B$ are induced by $\ov\varphi$).
Since $v'$ separates $A'$ from one chunk in $B'$, and $B'$ is connected, 
it must separate $A'$ from every chunk of $B'$, and hence from $B'$ itself.
Moreover, the path $\ga$ must pass through the vertex $v'$ separating 
$A'$ from $B'$.  We may thus write 
\[
\Delta'=\Delta'_1\bigcup_{v'_1=v'=v'_2}\Delta'_2\,,
\]
where $A'\subset\Delta'_1$ and $B'\subset\Delta'_2$.
We also decompose $\ga$ into a union of subpaths $\ga_1$ 
from $s'$ to $v_1'$ in $\Delta_1$ and $\ga_2$ from $v_2'$ 
to $t'$ in $\Delta_2$ where the endpoints $v_1'$, $v_2'$ are identified
to the vertex $v'$.
Now, by applying a sequence of edge twist isomorphisms to $G'$ (along the
edges of the subpath $\ga_2$ in $\Delta'_2$), we may
modify $\Delta'$ to the graph 
\[
\Delta'_1\bigcup_{v'_1=t'}\Delta'_2\,.
\]
By a similar sequence of edge twists (following the subpath $\ga_1$ 
from $v_1'$ back to $s'$ in $\Delta_1$), 
we may now further modify this graph to the graph
\[
\Delta'_1\bigcup_{s'=t'}\Delta'_2\,.
\]
Note that the above edge twists may be chosen to be the identity on $G(B)$. 
Composing with these edge twists therefore alters the isomorphism 
$\varphi \co G\to G'$ so that $\varphi_B=\tau_B$ is unchanged and 
$\varphi_A=g_A\circ\tau_A$ for a (possibly different) $g_A$ and $\tau_A$
satisfying $\tau_A(s)=\tau_B(s)=s'$ and $g_As' g_A\inv=s'$. 
But that is to say that $g_A\in C_{G'}(\< s'\>)$, and so, 
by applying a Dehn twist automorphism to $G'$, we may suppose that $g_A=1$. 
At this point the restriction of $\varphi$ to 
$G(A\cup B)$ is simply induced by a labelled graph isomorphism, as required.
This completes the proof of Theorem \ref{Thm1}.

\subsection{Proof of Theorem \ref{Thm2}}
The nontrivial content of Theorem \ref{Thm2} is contained in the statement that
\[
\ker(\pi,\Delta)=\Pure(\Delta)\rtimes\Inv(\Delta)\,.
\] 
To prove this we simply repeat the proof of Theorem \ref{Thm1} above 
with the added assumptions that $\Delta'=\Delta$ and that the map $\ov\varphi$ induced
on edges is the identity. One observes that only inversions and Dehn twist isomorphisms
are needed to complete the proof, and the statement of Theorem \ref{Thm2} follows.

\subsection{Proof of Theorem \ref{ThmCoxeter}}

To establish Theorem \ref{ThmCoxeter} it suffices to repeat once again the arguments
of Theorems \ref{Thm1} and \ref{Thm2}, replacing the isomorphism $\varphi\co G\to G'$ with 
an isomorphism $\phi\co  W\to W'$ between the corresponding Coxeter groups,
 and using Proposition \ref{CoxeterChunks} in the place of 
Propositions \ref{ThetaIsom} and \ref{chunkequiv}. We also make the following observations:
 
\begin{itemize}
\item The fact that $\phi$ induces a well defined bijection 
$\ov\phi\co E(\Delta)\to E(\Delta)$ follows 
by consideration of the action on the Davis complex, and the 
fact that $\ov\phi$ respects the labelling follows from the fact
that $W(e)$, being a dihedral group of order $2m_e$, is distinguished up to 
isomorphism by the label $m_e$.

\item The consideration of inversion automorphisms does not appear in the 
Coxeter group situation, but is replaced with a consideration of dihedral twist 
automorphisms.
In place of Lemma \ref{leafgroup} we make the following observations.
Let $e=\{s,t\}\in E(\Delta)$ and $m=m_e\geq 3$. 
Then the dihedral group $W(e)\cong D_{2m}$
has presentation $\<t,\rho|t^2=\rho^m=(\rho t)^2=1\>$, where $\rho=st$.
The reflections (conjugates of $s$ and $t$) are characterized as the 
primitive involutions of $W(e)$. Therefore, any automorphism of $W(e)$ 
can be modified by an inner automorphism and, if necessay, 
by the graph automorphism exchanging $s$ and $t$
so that it fixes the generator $t$ say.
The cyclic subgroup generated by $\rho$ is characteristic. 
So any automorphism of $W(e)$ which fixes the generator 
$t$ is given by $t\mapsto t$ and $\rho\mapsto \rho^k$ where 
$k$ represents a unit of the ring $\Z/m\Z$. 
Equivalently $t\mapsto t$ and $s\mapsto (st)^r s(st)^{-r}$ where 
$2r+1\equiv k(m)$. That is to say that any automorphism of $W(e)$
differs from a dihedral twist automorphism by a composition of inner
and graph automorphisms. This establishes the analogue of 
Lemma \ref{chunk-invariance} in the case of a non-solid chunk $A=e$.
\end{itemize}


\section{``Vertex links" in $\Theta$ and rank 2 vertices of the Deligne complex}
\label{sect:LinkTheta}

In order to arrive finally at a proof of Theorem \ref{Thm3} we need to pursue a little
further our study of the graph $\Theta$ so as to establish a rigidity property 
which is closely associated with the structure of links of rank 2 vertices
in the Deligne complex. We suppose throughout this section that $\Delta$ is a 
CLTTF defining graph and $\Theta$ the associated graph of fixed sets in the 
Deligne complex $\D$.

Let $\textsl{Chk}(\Theta):=\{ \text{ solid chunks of }\Theta\,\}$, and 
$\textsl{Chk}(\what\Delta):=\{ \text{ solid chunks of }\what\Delta\,\}$.
By virtue of Proposition \ref{chunkequiv} and Lemma \ref{betas} there exist
well-defined maps
\[
\be\co \textsl{Chk}(\Theta)\to G\hskip3mm\text{ and } \hskip3mm
\textsl{fund}\co \textsl{Chk}(\Theta)\to\textsl{Chk}(\what\Delta)
\]
such that $X=\be(X).\textsl{fund}(X)$, for each solid chunk $X$ in $\Theta$.

\begin{defn}[{\rm(}Oriented solid chunks\/{\rm)}] 
Let $G=G(\Delta)$ be a CLTTF Artin group, and $\Theta$ the associated graph 
of fixed sets. 
Fix $V\in\V$ and let $X$ be a solid chunk of $\Theta$ containing $V$.
By an \emph{orientation} of $X$ we mean a choice of vertex $F\in X$ adjacent to $V$
(necessarily, $F\in\F$). We write $\X=(X,V,F)$ for the \emph{oriented chunk based at $V$}
with orientation given by $F$. Note that there are always exactly two choices of 
orientation for a  chunk $X$ based at $V$. Namely, writing $V=\beta(X)V_e$ for some edge 
$e=\{ s,t\}\in E(\Delta)$, we have orientations given by $\beta(X)F_s$ and $\beta(X)F_t$.
We shall denote by $-\X$ the chunk $\X$ with the opposite orientation.

We say that two oriented solid chunks $\X_1=(X_1,V,F_1)$ and $\X_2=(X_2,V,F_2)$ 
based at $V$ in $\Theta$ are \emph{equivalent}, written $\X_1\simeq\X_2$, 
if $\X_2=g\X_1$ for some $g$ in the pointwise stabilizer of $F_1$ under the action of $G$. 
In other words, 
$(X_1,V,F_1)\simeq (X_2,V,F_2)$ if and only if, $\be(X_1)\inv\be(X_2)\in\< s\>$ and 
$F_1=F_2=\be F_s$, for some $s\in V(\Delta)$.
(Here $\be$ may be any element of the coset $\be(X_1)\Stab(F_s)$ of the 
\emph{setwise} stabilizer of $F_s$ under the action of $G$). 
\end{defn}

Note that the pointwise stabilizer of a fixed tree $F\in\F$ is strictly smaller
than its setwise stabilizer. Thus the equivalence class of an oriented solid chunk 
$(X,V,F)$ is not determined just by the pair $(V,F)$. We can however characterize the above
equivalence relation purely in terms of the structure of $\Theta$. Given an oriented 
solid chunk $\X=(X,V,F)$ we define the set
\[
N(\X)=\{ U\in\V\, :\, U\in X\,,\ U\neq V \text{ and $U$ is adjacent to } F\,\}\,.
\]
Note that $N(\X)$ and $N(-\X)$ are disjoint nonempty sets. 

\begin{lemma}\label{LVstruct}
Let $\X_i=(X_i,V,F_i)$ denote oriented chunks based at $V$, for $i=1,2$. 
Then the following are equivalent:
\begin{itemize}
\item[\rm(1)] $N(\X_1)= N(\X_2)$;
\item[\rm(2)] $N(\X_1)\cap N(\X_2)\neq\emptyset$;
\item[\rm(3)] $\X_1\simeq\X_2$.
\end{itemize} 
\end{lemma} 

\begin{proof}
Clearly (1)$\implies$(2), and (3)$\implies$(1) because any element of the 
pointwise stabilizer of $F_1$ fixes every vertex of $\Theta$ adjacent to $F_1$ 
(since $U\in\V$ is adjacent to $F_1$ if and only if $U\subset F_1$). We shall
show (2)$\implies$(3). Suppose $U \in N(\X_1)\cap N(\X_2)$. Then $F_1=F_2$ and 
$X_1$ and $X_2$ share the path $(V,F_1,U)$. 
Statement (3) now follows from Lemma \ref{two-circs}(ii). 
\end{proof}


\begin{defn}[{\rm(}The link graph $L(V,\Theta)\,$\/{\rm)}]
For each $V\in\V$ we define $L(V,\Theta)$ to be the graph with vertices  the
equivalence classes of oriented solid chunks based at $V$, and  an edge 
for each solid chunk $X$ containing $V$, the endpoints of which are determined by the
two possible orientations of $X$ based at $V$.
Thus, for each pair $\{\X,-\X\}$ of oriented solid chunks at $V$, 
there is a single edge in $L(V,\Theta)$ whose vertices are the 
just the equivalence classes of $\X$ and $-\X$ respectively.

Let $A$ denote a solid chunk of $\what\Delta$.
We shall write $L_A(V,\Theta)$ for the subgraph of $L(V,\Theta)$ consisting of those
edges associated to solid chunks  in the same $G$--orbit as $A$ (ie, solid chunks $X$
such that $\textsl{fund}(X)=A$).
\end{defn}

Note that the graph $L_A(V,\Theta)$ may often be empty, and will be non-empty if 
and only if some translate of the fundamental chunk $A$ contains $V$.
It is also possible that $V$ lies in no solid chunk whatsoever, in which case
the whole graph $L(V,\Theta)$ is empty. This happens at every vertex $V\in\V$ whenever
$\Delta$ is a tree.

The following Proposition shows that the vertex link graph just defined is,
on the one hand, canonical (with respect to $\V\F$--isomorphism) and, on the other hand,
strongly tied to the structure of the Deligne complex.

\begin{prop}\label{linkV}
Let $\Delta,\Delta'$ denote CLTTF defining graphs and $\Theta=\Theta(\Delta)$
and $\Theta'=\Theta(\Delta')$ the associated fixed set graphs. Let $V\in\Theta$
denote a type $\V$ vertex.
\begin{itemize} 
\item[\rm(i)] Any $\V\F$--isomorphism $\Phi\co  \Theta\to\Theta'$ 
induces a well-defined graph isomorphism
\[
\Phi_V\co L(V,\Theta)\to L(\Phi(V), \Theta')
\]
such that the vertex of $L(V,\Theta)$ represented by a based oriented solid 
chunk $(X,V,F)$ is mapped under $\Phi_V$ to the vertex of $L(\Phi(V),\Theta')$
represented by $(\Phi(X),\Phi(V),\Phi(F))$.

\item[\rm(ii)] The graph $L(V,\Theta)$ is a disjoint
union of the subgraphs $L_A(V,\Theta)$, for solid chunks $A\subset\what\Delta$,
and each nontrivial component $L_A(V,\Theta)$ is naturally isomorphic 
to $\Lk(V,\D)$, the link of 
$V$ in the Deligne complex. 
\end{itemize}
\end{prop}

\begin{proof}
{\bf (i)}\qua Since, by Proposition \ref{chunkequiv}, the $\V\F$--isomorphism $\Phi$
respects solid chunks, it will also respect \emph{oriented} solid chunks mapping the 
set $N(X,V,F)$ onto the set $N(\Phi(X),\Phi(V),\Phi(F))$. 
By Lemma \ref{LVstruct}, it follows that the equivalence relation 
$\simeq$ is preserved under the $\V\F$--isomorphism. 
The map $\Phi_V$ is therefore well-defined, and clearly a graph isomorphism.
\medskip

\noindent{\bf (ii)}\qua
Let $e=\{ s,t\}\in E(\Delta)$ denote the edge of $\Delta$ such that $V$ is a translate 
of $V_e$. Then the link $\Lk(V,\D)$ of the rank 2 vertex in $\D$ may be described as
follows. For each $g\in G$ such that $gV_e=V$ (ie, such that $V\subset gK$), there is 
exactly one edge in $\Lk(V,\D)$ contributed by the translate $gK$ of the fundamental 
region $K$. This edge has endpoints corresponding to the edges $[gV_e,gV_s]$ and $[gV_e,gV_t]$
emanating from $V=gV_e$. We shall denote this edge by $gE$ and the its endpoints by
$gS$ and $gT$ respectively. Clearly we have 
$gS=hS \iff g\inv h\in\< s\>$ and $gT=hT \iff g\inv h\in\< t\>$, while $gS\neq hT$, for all
$g,h\in G$. (The vertices of $\Lk(V,\D)$ may in fact be thought of as cosets in $G$ of  
subgroups $S=\<s\>$ and $T=\< t\>$ which lie in a common coset of $G(e)$).

Now suppose that $A$ is a solid chunk of $\Delta$ such that $L_A(V,\Theta)$ is non-empty.
Equivalently, $A$ contains the vertex $V_e$. Then each oriented solid chunk which
contributes to $L_A(V,\Theta)$ is either of the form $g\X_A=(gA,gV_e,gF_s)$ or 
of the form $-g\X_A=(gA,gV_e,gF_t)$
for some $g\in G$ such that $gV_e=V$.  It now follows, from the definition of $L(V,\Theta)$
and the above discussion, that mapping the edge $(g\X_A, -g\X_A)$ to the 
edge $gE=(gS, gT)$, for each $g$, defines a graph isomorphism. 

Finally we note that the graph $\Lk(V,\D)$ is connected (essentially because the group
$G(e)$ is generated by $s$ and $t$) and that it is clear from the definitions that 
the different subgraphs $L_A(V,\Theta)$ lie in different connected components of $L(V,\Theta)$.
\end{proof}

\begin{remark}
Note that the isomorphism of Proposition \ref{linkV}(ii) is natural in the sense that 
it is equivariant with respect to the obvious $\Stab(V)$ action on each graph.
\end{remark}

\begin{notation}[{\rm(}Generic rank 2 vertex link\/{\rm)}]
We shall adopt the notation suggested in the above proof in order to describe 
the link $\Lk(V,\D)$ of a generic rank 2 vertex $V$ of the Deligne
complex $\D$. For simplicity we shall suppose that $V=V_e$ where $e=\{ s,t\}\in E(\Delta)$
and we shall write $m=m_e$. Recall that the stabilizer, $\Stab(V_e)$, of this vertex 
under the action of $G$ on $\D$ is the group 
\[
G(e)=\< s,t\mid \pro(s,t;m)=\pro(t,s;m)\,\>\,,
\]
where $\pro(s,t;m)$ denotes the word $sts...$ of length $m$. 
Let $S=\< s\>$ and $T=\< t\>$.
The vertices of the graph $\Lk(V_e,\D)$ shall be represented by the cosets
of the subgroups $S$ and $T$ in $G(e)$ and, for each $g\in G(e)$,
there is a single edge with vertices $gS$ and $gT$, written $gE=\{ gS,gT\}$. 
(The symbol $E$ may be thought of as representing the trivial subgroup
$E=\{ 1\}$). 
The action of $G(e)$ on $\Lk(V_e,\D)$ (coming from the action of $G$ on $\D$) 
is defined in the obvious way, be left multiplication of cosets. 
Thus $\Stab_{G(e)}(gS)=g\< s\>g\inv$, $\Stab_{G(e)}(gT)=g\< t\>g\inv$ 
and $\Stab_{G(e)}(gE)=1$.
\end{notation}

We view $G(e)$ as the quotient of the free product $\<s\>\star\<t\>$ by the single
relation shown in the above presentation, and we make the following observation.
Locally geodesic circuits in $\Lk(V_e,\D)$ which pass through the fundamental 
edge $E$ correspond bijectively to cyclically reduced expressions over 
$\<s\>\star\<t\>$ for the identity in $G(e)$, ie, expressions
$w=a_1a_2\ldots a_n$ where the $a_i$ belong alternately to $\< s\>\setminus\{ 1\}$
and $\< t\>\setminus\{ 1\}$ ($a_1$ and $a_n$ belonging to distinct subgroups).
To be precise, a cyclically reduced expression for the identity in $G(e)$
which is written $w=a_1a_2\ldots a_n$ 
corresponds to the circuit $W=(E,a_1E,a_1a_2E,\dots , a_1..a_{n-1}E, wE=E)$ 
of the same length in $\Lk(V_e,\D)$. Note, also,
that any circuit in $\Lk(V_e,\D)$ may be translated by a graph automorphism 
(action by an element of $G(e)$) to a circuit passing through $E$.

\begin{lemma}\label{ShortCycles}
Suppose that $\Delta$ is the defining graph for a 2--dimensional Artin group.
Let $e=\{ s,t\}$ be an edge of $\Delta$ with label $m=m_e\geq 3$, 
and let $L=\Lk(V_e,\D)$ denote the link in the Deligne complex
of the rank 2 vertex $V_e$ fixed by $G(e)$.

Let $w$ denote a nonempty cyclically reduced expression over $\<s\>\ast\<t\>$
which represents the identity in $G(e)$, and write $\text{len}(w)$ for the (syllable)
length of $w$. Then $\text{len}(w)\geq 2m$ and if $\text{len}(w)=2m$ then, 
up to a cyclic permutation and inversion ($w\leftrightarrow w\inv$),
the word $w$ is one of the following balanced expressions, 
for some $n\in\Z\setminus\{ 0\}$,
\[
\begin{aligned}
s^nt\dots st(ts\dots ts^n)^{-1}\hskip3mm \text{ or }\hskip3mm 
t^ns\dots ts(st\dots st^n)^{-1}
\hskip3mm &\text{ if $m$ even, and}\\
s^nt\dots ts(ts\dots st^n)^{-1}\hskip3mm \text{ or }\hskip3mm 
t^ns\dots st(st\dots ts^n)^{-1}
\hskip3mm &\text{ if $m$ odd.}
\end{aligned}
\]
Equivalenty, every simple circuit in $L$ has edge length at least $2m$ and 
if it has edge length precisely $2m$ then it is a translate 
(by some element of $G(e)$ acting on $L$)
of one of the circuits through $E$ corresponding to the above expressions.
\end{lemma}

\begin{proof}
Recall that $G(e)$ acts by isometries on a regular $m$--valent tree $T$ in such a 
way that the generators $s$ and $t$ are each hyperbolic on $T$, the stabilizer 
of the midpoint of any edge is conjugate to $\< x\>$ where $x=\pro(s,t;m)$, 
and the kernel of the action is  the centre of $G(e)$, generated by the element $x$ 
if $m$ is even and $x^2$ if $m$ is odd. The tree may be embedded in the plane $\R^2$ 
(and the action extended non-isometrically) in such a way that the axis for 
each generator $s$ and $t$ (and each of their conjugates) bounds a connected 
component of $\R^2\setminus T$. Moreover, the action is such that the axes for $s$ 
and $t$ in $T$ intersect along a single edge $A$, but are oriented 
in opposite directions along this edge. 

We let $M$ denote the graph dual to $T$
in the plane. Observe that there is a natural $G(e)$--equivariant map $p\co L\to M$
which sends the edge $E$ to the edge of $M$ dual to $A$, and vertices $S$ and $T$ to
the vertices of $M$ lying in the regions bounded by the axes for $s$ and $t$ 
respectively. This map $p$ is in fact a covering projection. 
We also observe that any simple closed path $\rho$ in $M$
which starts at a vertex and runs exactly once around the boundary of a single
region of $\R^2\setminus M$ (thus a path of length $m$ surrounding a 
single vertex of the tree $T$) always lifts to a path
in $L$ of the form $g(E,sE, stE,stsE,..,\al E)$ for $g\in G$ and
$\al=\pro(s,t;m-1)$, or of similar form using the word $\be = \pro(t,s;m-1)$ 
or one of $\al^{-1}$ or $\be^{-1}$ in place of $\al$. 
Also, by choice of orientation of $\R^2$, we may suppose
that the lift of $\rho$ is associated with a positive word ($\al$ or $\be$) if and 
only if $\rho$ runs in a clockwise direction. 

More generally, any simple circuit $\rho$ in $M$ encloses a region 
containing a finite number of vertices of $T$, and if the simple circuit
$\rho$ surrounds exactly $N$ vertices of $T$ then it has length $(m-2)N +2\geq m$.
Moreover, the circuit $\rho$ lifts to a path in $L$ associated with a strictly positive 
or strictly negative word in the generators $s,t$ depending on whether it is
oriented in the clockwise or anti-clockwise direction.

Now consider a simple circuit $W$ in the graph $L$, corresponding to a reduced
expression $w$ over $\< s\>\star\< t\>$ for the identity in $G(e)$.
This projects to a locally geodesic (but not necessarily 
simple) circuit $\ov W$ in $M$. One may easily find a subpath of
$\ov W$ which describes a simple circuit in $M$. Thus we may decompose
$\ov W$ into the concatenation of paths $\rho.\rho'$ where $\rho$ is 
a simple circuit. In particular, the length of $\ov W$ is at least $m$.
Also, $\rho$ lifts to a path in $L$ associated to a word
$u$ in $s,t$ which is strictly positive or negative. Such a word
cannot represent the identity in $G(e)$, so is not equal to $w$. 
It follows that $\rho'$ is a nontrivial circuit in $M$. 
Repeating the above argument we have $l(\rho')\geq m$, and so
$l(\ov W)=l(\rho)+l(\rho')\geq 2m$. 
Moreover, $l(\ov W)=2m$ only if both $\rho$ and $\rho'$ 
are simple circuits, each surrounding a single vertex of $T$. 
Finally, in this case, since $w=1$ in $G(e)$ one of 
$\rho,\rho'$ is oriented clockwise, the other anti-clockwise,
and $w$ is necessarily given by one of the words 
listed in the statement of the Lemma.
\end{proof}

\begin{prop}\label{linkautos}
Suppose that $\Delta$ is the defining graph for a 2--dimensional Artin group.
Let $e=\{ s,t\}$ be an edge of $\Delta$ with label $m=m_e\geq 3$, 
and let $L=\Lk(V_e,\D)$ denote the link in the Deligne complex
of the rank 2 vertex $V_e$.\hfill\break
If $\tau$ is a graph automorphism of $L$ which 
fixes the fundamental edge $E=\{ S,T\}$ (ie, $\tau(S)=S$ and $\tau(T)=T$) 
then, either $\tau$ is the identity on $L$, 
or it is induced by the group inversion such that $s\mapsto s\inv$, and $t\mapsto t\inv$.
\end{prop}

\begin{proof}
Consider the circuits of minimal length $2m$ in $L$ which pass through the 
fundamental edge $E$. We observe that the edge pair $(E,sE)$ appears in infinitely 
many minimal length circuits, while $(E,s^kE)$, with $|k|>1$, 
appears in at most one or two minimal length circuits (depending on whether 
$m$ is odd or even). This implies
that the natural total order on the set $\{s^nE:n\in\Z\}$ (coming from the natural
ordering of the integers) is determined up to a reversal of order by graph theoretic
information.
It follows that, by composing $\tau$ with an inversion automorphism if
necessary, we may suppose that $\tau$ is the identity on the neighbourhood of $S$
(ie, the union of edges $s^kE$ for $k\in\Z$). Also, $\tau(tE)= tE$ or $t\inv E$.
Note also that, since the canonical cyclic ordering on any vertex of $L$ is
respected (up to reversal) by any graph automorphism, the family of minimal
length circuits associated to words of total \emph{word} length $2m$ ($n=\pm 1$ 
in Lemma \ref{ShortCycles}) is respected by $\tau$. These are the circuits associated
to the following cyclic words and their inverses
\[
\begin{aligned}
st\dots st(ts\dots ts)^{-1}
\hskip6mm &\text{ if $m$ even, and}\\
st\dots ts(ts\dots st)^{-1} 
\hskip6mm &\text{ if $m$ odd.}
\end{aligned}
\]
However, we observe that, in each of these cyclic words, 
the word $st\inv$ (or its inverse) appears exactly
once as a subword, while the word $st$ (or its inverse)
appears a total of $m-1$ times. Since $m\geq 3$, the paths
$(sE,E,tE)$ and $(sE,E,t\inv E)$ are thus differentiated by the number of minimal 
length circuits of this type which contain them. Therefore $\tau(tE)=tE$ and in fact 
$\tau$ must fix the whole neighbourhood of $T$. Since a similar argument may
be applied at each edge of $L$, and the graph is connected, it now 
follows that $\tau$ is the identity on the whole graph.
\end{proof}

Propositions \ref{linkV} and \ref{linkautos} together give the ``rigidity in the neighbourhood 
of a vertex'' property that will be needed in the following Section to complete the proof
of Theorem \ref{Thm3}.


\section{Abstract commensurators -- Theorem \ref{Thm3}}\label{sect:Comms}

Throughout this section we suppose that $\Delta,\Delta'$ denote CLTTF defining graphs.
For simplicity we shall write $G=G(\Delta)$, $G'=G(\Delta')$, and $\Theta$, $\Theta'$
for the associated fixed set graphs respectively.

\begin{prop}\label{VFisom}
If the defining graph $\Delta$ is not a tree (ie $\Delta$ contains at least one
simple circuit) then any graph isomorphism $\Theta\to\Theta'$ is a $\V\F$--isomorphism. 
\end{prop}

\begin{proof}
Let $\Phi \co \Theta\to\Theta'$ be a graph isomorphism. Recall that $\Theta$ 
and $\Theta'$ are connected bi-partite graphs. If $\Phi$ is not a $\V\F$--isomorphism
then we may suppose that $\Phi(\F)=\V'$ and $\Phi(\V)=\F'$. 
Since $\Delta$ is not a tree, we may choose some $V\in\V$ 
which lies in a solid chunk of $\Theta$ (so that $L(V,\Theta)\neq\emptyset$). 
Note that $L(V,\Theta)$ contains many simple closed
circuits (cf. Proposition \ref{linkV} and Lemma \ref{ShortCycles}).  Let 
$X_1,X_2,\dots,X_n$ denote a sequence of solid chunks in $\Theta$ which represents
a simple closed edge path in $L(V,\Theta)$ (each $X_i$ contains $V$).
The image $\Phi(X_i)$, $i=1,..,n$, of this sequence is a sequence of solid chunks
of $\Theta'$ all of which contain the vertex $\Phi(V)\in\F'$. 
By an application of Lemma \ref{two-circs}(iii), the conditions controlling 
adjacency of edges in $L(V,\Theta)$ (see Section \ref{sect:LinkTheta})
translate under $\Phi$ to the following condition. For each $i=1,..,n$,
there exists $e_i\in E(\Delta')$ such that 
\[
\be(X_i)\inv\be(X_{i+1}) = x_{e_i}^{m_i} \hskip3mm \text{ for some nonzero } m_i\in\Z\,. 
\] 
Moreover, we must have $e_i\neq e_{i+1}$, for all $i$. (Here indices are taken mod $n$).
But then we have that 
\[
x_{e_1}^{m_1}x_{e_2}^{m_2}\dots x_{e_n}^{m_n} =1\,,
\]
which contradicts the fact that the elements $\{ x_e : e\in E(\Delta')\}$ freely 
generate a free group, by Lemma \ref{freegroup}.
\end{proof}

\subsection{Proof of Theorem \ref{Thm3}}

We recall that, in the statement of the Theorem, $\Delta$ denotes a CLTTF defining graph
with no separating edge or vertex.

\paragraph{(i)} 
Suppose $G(\Delta)$ is abstractly commensurable to $G(\Delta')$
for some CLTTF defining graph $\Delta'$. We wish to show that $\Delta$ and $\Delta'$
are label isomorphic.

The condition that $\Delta$ has no separating edge
or vertex simply means that $\what\Delta$ is itself a solid chunk of $\Theta$
(the unique fundamental chunk in this case). 
In particular, $\Delta$ is not a tree and Proposition
\ref{VFisom} applies. Let $\varphi\in\Comm(G(\Delta),G(\Delta'))$. 
By Proposition \ref{VFisom} and Proposition \ref{ThetaIsom}, $\varphi$ induces
a $\V\F$--isomorphism $\Phi\co \Theta\to\Theta'$ which, by Proposition \ref{chunkequiv},
maps solid chunks of $\Theta$ to solid chunks of $\Theta'$.  
But then, up to modification of $\varphi$ by an inner automorphism, we may suppose
that $A:=\Phi(\what\Delta)$ is a solid chunk of $\what\Delta'$. 
Thus $\Delta$ is isomorphic to a subgraph of $\Delta'$. 
Moreover, the isomorphism respects labels because
the label $m_e$ is determined by the structure of the link graph $L(V_e,\Theta)$.
Namely, the shortest simple closed path in $L(V_e,\Theta)$ has length $2m_e$ (cf.
Lemma \ref{ShortCycles} and Proposition \ref{linkV}).

Recall that the Deligne complex $\D$ (of type $\Delta$) may be described as the 
universal cover of a complex of groups structure over the fundamental region $K$.
We write $\D'$ and $K'$ for the Deligne complex of type $\Delta'$
and its fundamental region. There is a naturally defined subcomplex $K_A\subset K'$ 
associated to the fundamental chunk $A$ (which is spanned by those vertices 
corresponding to standard parabolic subgroups lying in $G(A)$), and 
we define the following subcomplex of $\D'$:
\[
\D_A=\bigcup_{g\in G(A)}\, gK_A\,.
\]
Clearly, $\D_A$ is an isometric copy of the Deligne complex associated to $G(A)$ 
sitting inside $\D'$. We claim that the map $\Phi\co \Theta\to\Theta'$ induces an 
isometry $\D\to\D_A$. 
 
The fact that $\Phi$ maps chunks of $\Theta$ to solid chunks of $\Theta'$ means
that there is a naturally induced family of isomorphisms 
\[
gK\to \phi(g)K_{B(g)} \hskip3mm\text{ for each }g\in G(\Delta)\,.
\]
where $\phi\co G(\Delta)\to G(\Delta')$ and $B\co G(\Delta)\to\textsl{Chk}(\Delta')$
are simply functions. By the discussion in the opening paragraph, 
We have that $\phi(1)=1$, $B(1)=A$, and the map $K\to K_A$ is induced by a
label isomorphism between the graphs $\Delta$ and $A$.

Let $e\in E(\Delta)$, and $e'$ its image in $A$. Then, 
by Proposition \ref{linkV}, the $\V\F$--isomorphism $\Phi$
must induce an isomorphism $L(V_e,\Theta)\cong L(V_{e'},\Theta')$, and 
there is a naturally induced isomorphism $\Lk(V_e,\D)\cong \Lk(V_{e'},\D')$. 
In particular, the function $\phi$ restricts to 
an isomorphism $G(e)\to G(e')<G(A)$ and $B(g)=A$ for all $g\in G(e)$. 
This naturally induces a well-defined isometric embedding 
of the neighbourhood of a rank 2 vertex of $\D$ into $\D_A$. 
Applying the same argument at every rank 2 vertex of $\D$ and 
using the fact that $\D$ is connected we obtain a map $\Phi_\D\co \D\to\D_A$ which is 
locally isometric, so globally isometric since $\D$ is CAT(0).
The isomorphism $\Phi_\D$ is natural in the sense
that if $H<G$ is the domain of $\varphi$, then 
$\Stab_{\varphi(H)}(\Phi_\D(p))=\varphi(\Stab_H(p))$ for all vertices $p\in\D$.

It follows from the above argument that the image of the abstract 
commensurator $\varphi$ is a finite index subgroup of $G(A)$. However, since $G(A)$
is infinite index in $G(\Delta')$ unless $A=\Delta'$, 
statement (i) of Theorem \ref{Thm3} follows.

\paragraph{(ii)} 
Suppose now that $\Delta$ satisfies the vertex rigidity condition:
\begin{description}
\item[\rm(VR)] Any label preserving automorphism of $\Delta$ which
fixes the neighbourhood of a vertex is the identity automorphism,
\end{description}
and write $G=G(\Delta)$. We wish to show that $\Comm(G)\cong\Aut(G)$.

As in the proof of part (i) (with $\Delta'=\Delta$), 
we may suppose, up to modification of $\varphi$ by
an inner automorphism and a graph automorphism of $G$, that $\varphi$ naturally induces
an isometry $\Phi_\D$ of $\D$ which is the identity on the fundamental region $K$. 
Moreover, by Lemma \ref{linkautos}, $\Phi_\D$ induces either the identity
or the ``inversion'' automorphism on the link of each rank 2 vertex in $K$. 
Fixing $e\in E(\Delta)$, we may suppose, up to modification of $\varphi$ by a 
global inversion if necessary, that $\Phi_\D$ induces the identity on 
$\Lk(V_e,\D)$, and hence restricts to the identity on a small 
open neighbourhood of $V_e$ in $\D$.

We now use the hypothesis that $\Delta$ satisfies the vertex rigidity condition (VR) 
to show that $\Phi_\D$ is the identity on the whole of $\D$. 
On the one hand, if $\Phi_\D$ is the identity on a small open
neighbourhood of any rank 2 vertex $p$ of $\D$ then by (VR) it is the identity on every 
translate of $K$ adjacent to this vertex. On the other hand, 
if $p$ and $q$ are rank 2 vertices joined
by a path $(p,r,q)$ in $\mF$, 
where $r$ is a rank 1 vertex, then we  observe that 
every translate of $K$ adjacent to $r$ is also adjacent to both $p$ and $q$. 
If $\Phi_\D$ is the identity on every translate of $K$ adjacent to $p$ then, since
at least two, in fact infinitely many, of these translates are also adjacent to 
$r$ and $q$, it follows, by Lemma \ref{linkautos}, 
that $\Phi_\D$ induces the identity on the link of $q$.
This argument, together with the statement involving the (VR) hypothesis,
shows that if $\Phi_\D$ is the identity on the neighbourhood of $p$ 
then it is the identity on every translate of $K$ which is adjacent to 
a rank 2 vertex $q$ within a ball of radius 2 about $p$ in the Deligne complex.
Since, from the previous paragraph, we may suppose that $\Phi_\D$ is 
the identity on a neighbourhood of
the vertex $V_e$, by applying this 
argument inductively and appealing to the connectedness of $\D$, 
we show that the map $\Phi_\D$ must be the identity on 
the neighbourhood of every rank 2 vertex in $\D$, 
and therefore equal to the identity everywhere.

Finally, if $\Phi_\D$ is the identity then $\varphi$ must also be the identity (since
for each $h\in H$, $\Phi_\D(hK)=\varphi(h)K$). This completes the proof of 
statement (ii) of Theorem \ref{Thm3}.

\section{Examples of nontrivial abstract commensurators}\label{sect:Examples}

We conclude by giving some examples of abstract commensurators which illustrate the
situations one might need to consider in order to extend Theorem \ref{Thm3}.
We begin with the necessity of the (VR) hypothesis in part (ii) of the Theorem.

\begin{example}
Let $\Delta$ be the CLTTF defining graph shown in
\figref{ExampleFig}(i).  Note that $\Delta$ has no separating edge or
vertex, but does \emph{not} satisfy the vertex rigidity condition
(VR).  The standard generators of $G(\Delta)$ are labelled $u,v,x,y,z$
as shown in the Figure. Let $X$ denote the presentation 2--complex of
the standard presentation of $G=G(\Delta)$ given in the
introduction. Thus $X$ has a single vertex, an oriented labelled
1-cell for each of the generators, and a 2--cell corresponding to each
relator in the presentation, and $\pi_1(X)=G$.  

\begin{figure}[ht!]\small\anchor{ExampleFig}
\psfrag {4}{$4$}
\psfrag {x}{$x$}
\psfrag {y}{$y$}
\psfraga <-1pt,0pt> {z}{$z$}
\psfrag {A}{$A$}
\psfrag {B}{$B$}
\psfrag {v}{$v$}
\psfrag {A}{$A$}
\psfrag {i}{(i)}
\psfraga <6pt,6pt>  {*}{$\ast$}
\psfrag {ii}{(ii)}
\cl{\includegraphics[width=12cm]{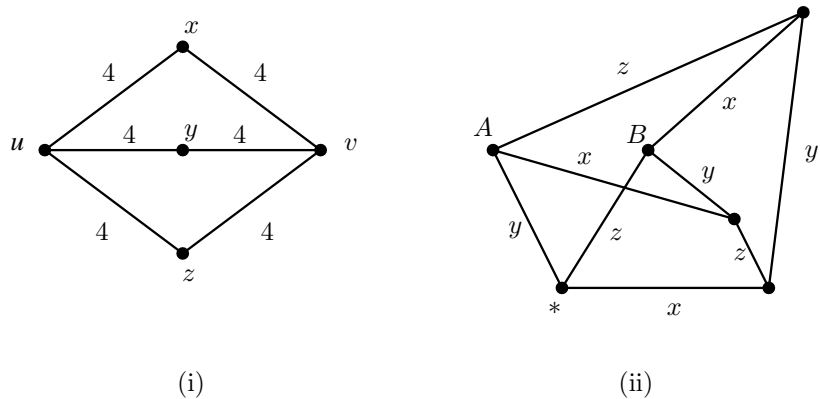}}
\vspace{-10mm}
\caption{(i)\qua A non-``vertex rigid'' defining graph $\Delta$ and
(ii) a recipe $L$ 
for an index 6 subgroup of $G=G(\Delta)$ which admits an automorphism
not induced from an element of $\Aut(G)$}
\label{ExampleFig}
\end{figure}
 
The labelled graph $L$ shown in \figref{ExampleFig}(ii) is a recipe for 
building a finite index cover $\wtil X$ of $X$, and so represents a finite index 
subgroup of $G$, as follows. Let the vertices of $\wtil X$ be in bijection
with the vertices of $L$. For each edge $(P,Q)$ in $L$ labelled with a generator $w$ 
of $G$ there are two oriented 1-cells between $P$ and $Q$ in $\wtil X$, 
each labelled $w$, one oriented from $P$ to $Q$, and the other in the opposite sense. 
At each vertex $P$, $\wtil X$ has an oriented 1-cell (a loop from $P$ back to $P$)
labelled $u$, and another labelled $v$. This defines the 1-skeleton of $\wtil X$.
Note that the labelling and orientation define a 6-fold 
covering map $\wtil{X}^{(1)}\to X^{(1)}$. We now define the 2--cells of $\wtil X$ 
in the unique way that will enable us to extend this covering map to a 6-fold covering
map $\wtil X\to X$. Note that there is no obstruction to doing this
because every edge label in $\Delta$ is even.  Finally, choose a
basepoint $\ast$ for $\wtil X$ as indicated in \figref{ExampleFig}(ii).

Let $H=\pi_1(\wtil X, \ast)$ denote the index 6 subgroup of $G$ associated with this covering
map and observe that any automorphism of the underlying graph of $L$ 
induces an automorphism of the group $H$. Let $\varphi\co H\to H$ denote the automorphism 
induced by exchanging the vertices labelled $A$ and
$B$ in \figref{ExampleFig}(ii), and leaving all other vertices of $L$ fixed. 

We remark that $\varphi$ is not induced by any automorphism of $G$.
To see this, we note that $x^2, y^2, z^2, u$ and $v$ are all elements of $H$ and 
$\varphi$ exchanges $y^2$ and $z^2$ while fixing $x^2$, $u$, and $v$.
Therefore the only candidate for an element of $\Aut(G)$ which induces $\varphi$
would be the graph automorphism $\tau$ which exchanges generators $y$ and $z$, leaving
all other generators fixed. However $\varphi$ also fixes the element $xz^2x\inv$
while $\tau(xz^2x\inv)=xy^2x\inv$. Thus $\varphi$ and $\tau$ are inequivalent as
elements of $\Comm(G)$.
\end{example}

We state the next example in the form of a lemma:

\begin{lemma}
Let $\Delta,\Delta'$ denote Artin defining graphs. Suppose that
$e=\{s,t\}\in E(\Delta)$ is a cut edge of $\Delta$, equivalently, 
$e$ is itself a (non-solid) chunk of $\Delta$.
Suppose moreover that $\Delta$ and $\Delta'$ differ only in the label on the edge $e$,
but that this label is at least 3 in each case.
Then the Artin groups $G(\Delta)$ and $G(\Delta')$ are abstractly commensurable.
\end{lemma}

\begin{proof}
If $A\subset\Delta$ is any full subgraph, and $n\geq 1$, then we write
$H(A;n)$ for the index $n$ subgroup of $G(A)$ which is the kernel of the
mod $n$ length function (the group of elements $x$ such that 
$\ell(x)\equiv 0 \mod n$).
 
Let $e=\{s,t\}$ be an edge with label $m_e\geq 3$.
Let $k=\textsl{lcm}(m_e,2)$ and let $n$ be any positive multiple of $k$.
Then, since the order of every torsion element of $G(e)/Z$ divides $k$, it follows 
that $H(e;n)\cong F\times\Z$ where $F$ is a finitely generated nonabelian free group
(nonabelian since $m_e\geq 3$).
Moreover, up to isomorphism of $F\times\Z$, we may suppose that the subgroups 
$H(s;n)=H(e;n)\cap\< s\>=\< s^n\>$ and  $H(t;n)=H(e;n)\cap\< t\>=\< t^n\>$ 
are free factors of the subgroup $F$. 
That is $H(e;n)= (F'\star \< s^n\>\star\< t^n\>)\times\Z$. 
Since all finite rank free groups are abstractly commensurable we may suppose,
up to an abstract commensurator which fixes the subgroups $\<s^n\>$ and $\< t^n\>$,
that the rank of $F'$ is any given integer. It follows that, if $e_i=\{s_i,t_i\}$ are
edges, for $i=1,2$, with labels $m_i\geq 3$ respectively, then $G(e_1)$ and $G(e_2)$ are 
abstractly commensurable by a commensurator which maps $s_1^n\mapsto s_2^n$
and $t_1^n\mapsto t_2^n$, for sufficiently large $n$ (we may take $n=k_1k_2$ where
$k_i=\textsl{lcm}(m_i,2)$). 

Now suppose that $e=\{ s,t\}$ is a cut edge of the defining graph $\Delta$, and write
$\Delta=\Delta_1\cup_s e\cup_t\Delta_2$. Then $G(\Delta)$ is an amalgmated product
\[
G(\Delta)=G(\Delta_1)\star_{\< s\>} G(e) \star_{\<t\>} G(\Delta_2)\,.
\]
We consider two choices of the edge label $m_e$, writing 
$\Delta,\Delta'$ for the two defining graphs thus obtained, and $k, k'$ for the
corresponding values of $\textsl{lcm}(m_e,2)$. Let $n=kk'$. Then the subgroup 
$H(\Delta;n)$ is written as an amalgamated product as follows
\[
H(\Delta;n)=H(\Delta_1;n)\star_{\< s^n\>} H(e;n) \star_{\<t^n\>} H(\Delta_2;n)\,.
\]
It follows from remarks in the previous paragraph that $H(\Delta;n)$ and 
$H(\Delta';n)$ are abstractly commensurable. Thus $G(\Delta)$ and $G(\Delta')$ are
abstractly commensurable.
\end{proof}

It would be interesting to give a classification of all CLTTF Artin groups up to abstract
commensurability. Theorem \ref{Thm3} gives a partial result in this direction. 
The above Lemma shows that, in order to give a complete treatment of the question,
it suffices to consider only those CLTTF defining graphs where every cut edge 
(equivalently, every edge that does not lie in a circuit) is labelled 3. 
Moreover, by applying twist isomorphisms we may 
further restrict our attention to the case where every such edge contains a 
terminal vertex. The following example suggests that even amongst these defining 
graphs there may be many non-obvious commensurations.

\begin{lemma}
Let $\Delta$ denote an arbitrary defining graph, and let $s\in V(\Delta)$.
Let $(\Delta_i,s_i)$ denote a label isomorphic copy of $(\Delta,s)$, for each $i\in\N$,
 and let $E$ denote the graph
consisting of a single edge $E=\{s_0,t\}$ with label $m_E=3$. 
For $n\in\N$ we write $\Delta^{(n)}$ for the union of the labelled graphs 
$E,\Delta_1,\dots ,\Delta_n$ with the vertices $s_0,s_1,..,s_n$ identified to
a single vertex $s$.

Then the Artin groups $G(\Delta^{(n)})$ and $G(\Delta^{(m)})$ 
are abstractly commensurable for all $m,n\geq 1$.
\end{lemma}

\begin{proof}
We note that, for a sufficiently large $k\in\N$, $G(\Delta^{(1)})$ is abstractly 
commensurable to $G_1:= H(\Delta;k)\star_C((C\star F)\times\Z)$ where $C$ denotes the 
cyclic subgroup of $G(\Delta^{(1)})$ generated by $s^k$, and $F$ denotes a  
free group of finite rank at least 2. More generally $G(\Delta^{(n)})$ is abstractly
commensurable to the amalgamated product
\[
G_n:= H(\Delta_1;k) \star_C \dots\star_C H(\Delta_n;k) \star_C ((C\star F)\times \Z)\,.
\]
We now observe that the group $G_n$ is isomorphic to an index $n$ subgroup of $G_1$. 
Let $\lambda \co  G_1\to\Z$ denote the surjective homomorphism  whose kernel contains 
the free factor $H(\Delta;k)$ as well as the subgroup $C\star F$ of the remaining factor. 
Then $G_n$ is isomorphic to the kernel of the quotient map $G_1\to \Z/n\Z$ which 
factors through $\lambda$.
Thus, each group $G_n$ is abstractly commensurable to $G_1$, completing the proof.
\end{proof}

\end{document}

%% file: gtoutput.tex
\def\ifplaintex{\expandafter\ifx\csname documentclass\endcsname\relax}

\ifplaintex 
\else
\fi

\expandafter\ifx\csname beginpicture\endcsname\relax
\expandafter\ifx\csname documentclass\endcsname\relax
\input pictex \else
\input prepictex \input pictex \input postpictex \fi\fi

\def\gt{{\mathsurround=0pt\it $\cal G\mskip-2mu$eometry \&\ 
$\cal T\!\!$opology}}        

\def\gtp{{\mathsurround=0pt\it $\cal G\mskip-2mu$eometry \&\ 
$\cal T\!\!$opology $\cal P\!$ublications}}  

\def\lognumber#1{\def\thelognumber{#1}}
\def\volumenumber#1{\def\thevolumenumber{#1}}
\def\papernumber#1{\def\thepapernumber{#1}}
\def\volumeyear#1{\def\thevolumeyear{#1}}

\def\pagenumbers#1#2{\def\startpage{#1}\def\finishpage{#2}}
\def\published#1{\def\publishdate{#1}}
\def\proposed#1{\def\theproposer{#1}}
\def\seconded#1{\def\theseconders{#1}}
\def\received#1{\def\receiveddate{#1}}
\def\revised#1{\def\reviseddate{#1}}
\def\accepted#1{\def\accepteddate{#1}}
\def\asciititle#1{\def\theasciititle{#1}}

\def\asciiaddress#1{\def\theasciiaddress{#1}}

\long\def\asciiabstract#1{\long\def\theasciiabstract{#1}}
\def\asciikeywords#1{\def\theasciikeywords{#1}}

\def\shorttitle#1{\def\theshorttitle{#1}}

\let\\\par\let\thelognumber\relax
\let\thevolumenumber\relax\let\thepapernumber\relax
\let\thevolumeyear\relax\let\thesamplenumber\relax\let\startpage\relax
\let\finishpage\relax\let\publishdate\relax\let\receiveddate\relax
\let\reviseddate\relax\let\accepteddate\relax\let\theasciititle\relax
\let\theasciiauthors\relax\let\theasciiaddress\relax
\let\theasciiabstract\relax\let\theasciikeywords\relax
\let\theasciiemail\relax\let\theshortauthors\relax\let\theshorttitle\relax

\long\def\maketitlep{   

\count0=\startpage

\gt\hfill      
\beginpicture
\setcoordinatesystem units <0.33truein, 0.33truein> point at 2.2 0.9
\setplotsymbol ({$\cal G$})
\plotsymbolspacing=9truept
\circulararc 315 degrees from 0 1 center at 0 0
\setplotsymbol ({$\cal T$})
\circulararc 315 degrees from 1 -1 center at 1 0
\endpicture
\break
{\small\ifx\thesamplenumber\relax 
Volume \else Sample
\fi\thevolumenumber\ (\thevolumeyear)
\startpage--\finishpage\nl
Published: \publishdate}
\vglue 0.5truein plus 0.4fil minus 0.1truein

{\parskip=0pt\leftskip 0pt plus 1fil\def\\{\par\smallskip}{\ifplaintex\large
\else\Large\fi\bf\thetitle}\par\medskip}   

\vglue 0pt plus 0.1fil 

{\parskip=0pt\leftskip 0pt plus 1fil\def\\{\par}{\sc\theauthors}
\par\medskip}

\vglue 0pt plus 0.1fil 

{\small\parskip=0pt\let\newline\\
{\leftskip 0pt plus 1fil\def\\{\par}{\sl\theaddress}\par}
\expandafter\ifx\theemail\relax    
\relax\else\vglue 5pt plus 0.02fil minus 2pt\def\\{\stdspace{\rm 
and}\stdspace} 
\cl{Email:\stdspace\tt\theemail}\fi
\ifx\theurl\relax                  
\relax\else\vglue 5pt plus 0.02fil minus 2pt\def\\{\stdspace{\rm 
and}\stdspace}
\cl{URL:\stdspace\tt\theurl}\fi\par}

\vglue 7pt plus 0.3fil minus 3pt

{\bf Abstract}
\vglue 5pt plus 0.1fil minus 2pt

\theabstract

\vglue 7pt plus 0.3fil minus 3pt

{\bf AMS Classification numbers}\quad Primary:\quad \theprimaryclass

Secondary:\quad \thesecondaryclass

\vglue 5pt plus 0.3fil minus 2pt

{\bf Keywords:}\quad \thekeywords

\vglue 10pt plus 0.5fil minus 5pt

{\small  Proposed: \theproposer\hfill Received: \receiveddate\nl
Seconded: \theseconders\hfill 
\ifx\reviseddate\relax                         
Accepted: \accepteddate                        
\else
Revised: \reviseddate                          
\fi}
\eject
}       

\let\maketitlepage\maketitlep
\let\maketitle\maketitlepage

\font\phead=cmsl9 scaled 950
\font=cmsl9 scaled 1050
\font\pnum=cmbx10 scaled 913
\font\lnum=cmbx10 
\font\pfoot=cmsl9 scaled 950
\font=cmsl9 scaled 1050
\ifplaintex
\headline{\vbox to 0pt{\vskip -4.5mm\line{\small\phead\ifnum
\count0=\startpage ISSN 1364-0380 (on line)
1465-3060 (printed) \hfill {\pnum\folio}\else\ifodd\count0\def\\{ }%
\ifx\theshorttitle\relax\thetitle\else\theshorttitle\fi\hfill{\pnum\folio}
\else\def\\{ and }{\pnum\folio}\hfill\ifx\theshortauthors\relax\theauthors
\else\theshortauthors\fi\fi\fi}\vss}}
\footline{\vbox to 0pt{\vglue 0mm\line{\small\pfoot\ifnum\count0=\startpage
\copyright\ \gtp\hfill\else
\gt, Volume \thevolumenumber\ (\thevolumeyear)\hfill\fi}\vss
}}
\else
\makeatletter
\def\@oddhead{{\small\lhead\ifnum\count0=\startpage ISSN 1364-0380 (on line)
1465-3060 (printed) \hfill {\lnum\number\count0}\else\ifodd\count0
\def\\{ }\ifx\theshorttitle\relax \thetitle \else\theshorttitle\fi\hfill
{\lnum\number\count0}\else\def\\{ and }{\lnum\number\count0}
\hfill\ifx\theshortauthors\relax 
\theauthors\else\theshortauthors\fi\fi\fi}}\def\@evenhead{\@oddhead}
\def\@oddfoot{\small\lfoot\ifnum\count0=\startpage\copyright\ \gtp\hfill\else
\gt, Volume \thevolumenumber\ (\thevolumeyear)\hfill\fi}
\def\@evenfoot{\@oddfoot}
\makeatother
\fi

\newwrite\gtoutfile
\long\gdef\makeheadfile{  
{\def\\{, }\def\s{ }
\immediate\openout\gtoutfile head.xxx
\immediate\write\gtoutfile{Proxy-for: \ifx\theasciiauthors\relax
\theauthors\else\theasciiauthors\fi\s<\ifx\theasciiemail\relax\theemail\else\theasciiemail\fi>}
\immediate\write\gtoutfile{\noexpand\\}
\immediate\write\gtoutfile{Authors: \ifx\theasciiauthors\relax
\theauthors\else\theasciiauthors\fi}
{\def\\{ }\immediate\write\gtoutfile{Title: \ifx\theasciititle\relax
\thetitle\else\theasciititle\fi}}
\immediate\write\gtoutfile{Subj-class: GT or SG or MG etc}
\immediate\write\gtoutfile{MSC-class: \theprimaryclass\ifx\thesecondaryclass\relax\else, \thesecondaryclass\fi}
\immediate\write\gtoutfile{Journal-ref: Geom. Topol. \thevolumenumber
(\thevolumeyear) \startpage-\finishpage}
\immediate\write\gtoutfile{Comments: Published by Geometry and Topology at}
\immediate\write\gtoutfile{\s\s http://www.maths.warwick.ac.uk/gt/GTVol\thevolumenumber/paper\thepapernumber.abs.html}
\immediate\write\gtoutfile{\noexpand\\}
\immediate\write\gtoutfile{}
\ifx\theasciiabstract\relax
\immediate\write\gtoutfile{\theabstract}\else
\immediate\write\gtoutfile{\theasciiabstract}\fi
\immediate\write\gtoutfile{}
\immediate\write\gtoutfile{\noexpand\\}
\immediate\write\gtoutfile{}
\immediate\closeout\gtoutfile}}  

\def\maketitlepage{\maketitlep\makeheadfile}
\let\maketitle\maketitlepage